\documentclass[preprint, 10pt]{elsarticle}

\usepackage{palatino}
\usepackage[usenames]{color}
\usepackage{epsfig}
\usepackage[fleqn,reqno]{amsmath}
\usepackage{amsfonts,amsthm,bm}
\usepackage{amssymb}
\usepackage{mathrsfs}
\usepackage{stmaryrd}
\usepackage{subfigure}
\usepackage[titletoc]{appendix}
\usepackage{tabularx,booktabs}
\usepackage{graphics,caption}
\usepackage{fancybox}
\usepackage{lineno}
\usepackage[top=1.2in,bottom=1.2in,left=1in, right=1in]{geometry}
\usepackage{array}
\usepackage{multirow}
\usepackage{wrapfig}
\usepackage{algorithmic,algorithm}
\usepackage{paralist}
\usepackage{ifthen}
\usepackage{enumitem}
\usepackage{tikz,pgfplots,filecontents}
\usepackage{mdframed}

\usepackage[pagebackref=false,bookmarks=false]{hyperref} 

\definecolor{dblue}{rgb}{0.03,0.3,0.62}
\definecolor{dorange}{rgb}{1,0.55,0}

\hypersetup{
  bookmarksnumbered=true,
  bookmarksopen=false,
  hypertexnames=false,      
  breaklinks=true,          
  unicode=false,
  pdffitwindow=true,        
  pdfnewwindow=true,        
  colorlinks=true,         
  linkcolor=dblue,
  anchorcolor=red,
  citecolor=dorange,
  filecolor=magenta,
  urlcolor=dblue,
  pdfstartview = FitH,
  pdfkeywords = {},
  pdfcreator = {LaTeX with hyperref package}
}

\definecolor{sblue}{cmyk}{0.98,0.13,0,0.43} 
\definecolor{sblue}{cmyk}{0.98,0.13,0,0.43} 

\newcommand{\bigO}{\mathcal{O}}
\newcommand{\Div}{\mathrm{Div}}

\newcommand{\xx}{{\mathbf{x}}}
\newcommand{\yy}{{\mathbf{y}}}
\newcommand{\zz}{{\mathbf{z}}}
\newcommand{\RR}{{\mathbb{R}}}
\newcommand{\rr}{{\mathbf{r}}}
\newcommand{\ee}{{\mathbf{e}}}
\newcommand{\txx}[1]{{\tilde{\xx}_{#1}}}
\newcommand{\tee}[1]{{\tilde{\ee}_{#1}}}
\newcommand{\exx}[1]{{\mathbf{e}_{\xx_{#1}}}}
\newcommand{\tsigma}[1]{{\tilde{\sigma}_{#1}}}
\newcommand{\esigma}[1]{{e_{\sigma_{#1}}}}
\newcommand{\ff}{{\mathbf{f}}}
\renewcommand{\SS}{{\mathcal{S}}}
\newcommand{\BB}{{\mathcal{B}}}
\newcommand{\DD}{{\mathcal{D}}}
\newcommand{\TT}{{\mathcal{T}}}

\newcommand{\vv}{{\mathbf{v}}}
\newcommand{\eeta}{{\boldsymbol\eta}}
\newcommand{\new}{{\mathrm{new}}}
\newcommand{\opt}{{\mathrm{opt}}}
\newcommand{\sdc}{{\mathrm{sdc}}}
\newcommand{\gmres}{{\mathrm{gmres}}}
\newcommand{\up}{{\mathrm{up}}}
\newcommand{\down}{{\mathrm{down}}}
\newcommand{\wall}{{\mathrm{wall}}}
\newcommand{\boldnsdc}[1]{{\boldsymbol{n_{\sdc}=#1}}}
\newcommand{\mcaption}[2]{\caption{\small \em #1}\label{#2}}
\newtheorem{theorem}{Theorem}[section]

\newif\ifTikz
\Tikztrue

\begin{document}

\title{Adaptive Time Stepping for Vesicle Suspensions}

\author[ut]{Bryan Quaife}\ead{quaife@ices.utexas.edu}
\author[ut]{George Biros}\ead{gbiros@acm.org}

\address[ut]{Institute of Computational Engineering and Sciences,\\
  The University of Texas at Austin, Austin, TX, 78712.}

\begin{abstract} 
We present an adaptive arbitrary-order accurate time-stepping numerical
scheme for the flow of vesicles suspended in Stokesian fluids. Our
scheme can be summarized as an approximate implicit spectral deferred
correction (SDC) method.  Applying a textbook fully implicit SDC scheme
to vesicle flows is prohibitively expensive. For this reason we
introduce several approximations.  Our scheme is based on a
semi-implicit linearized low-order time stepping method. (Our
discretization is spectrally accurate in space.)  We also use invariant
properties of vesicle flows, constant area and boundary length in two
dimensions, to reduce the computational cost of error estimation for
adaptive time stepping. We present results in two dimensions for
single-vesicle flows, constricted geometry flows, converging flows, and
flows in a Couette apparatus. We experimentally demonstrate that the
proposed scheme enables automatic selection of the step size and
high-order accuracy.
\end{abstract}


\maketitle

\section{Introduction\label{s:intro}}
Vesicles are deformable capsules filled with a viscous fluid. Vesicle
flows refer to flow of vesicles that are suspended in a Stokesian fluid.
Vesicles are differentiated from other capsules in the balance of forces
in their interface.  In particular, their boundary is locally
inextensible (in 3D is locally incompressible).  Vesicle flows find many
applications many biological applications such as the simulation of
biomembranes~\cite{sac1996} and red blood
cells~\cite{ghi:rah:bir:mis2011,kao:tah:bib:ezz:ben:bir:mis2011,mis2006,nog:gom2005,poz1990}.

The dynamics of a vesicle is governed by bending, tension (enforces
inextensibility), hydrodynamics forces from other vesicles in the
suspension, and possibly hydrodynamic forces from flow confinement
boundary walls.  Vesicle suspensions are modeled using the Stokes
equations, a jump in the stress across each vesicle to match the
interfacial forces on the membrane the vesicle, and a no-slip boundary
condition on the vesicles and the confinement walls.

A significant challenge in simulating vesicle flows is that their
governing equations are stiff. One stiffness source is associated with
the interfacial forces, bending forces, and inextensibility related
tension forces (depends mainly on the curvature of the vesicle).
Another stiffness source is the hydrodynamic interaction between
vesicles (depends on the minimum distance between vesicles). A third
stiffness source is the hydrodynamic interaction between vesicles and
external boundary walls (depends on the minimum distance between the
vesicles and the wall).  Therefore, to resolve the dynamics of a
vesicle, it may be necessary to take a small time step when the vesicle
has regions of large curvature, or when a vesicle approaches a boundary
wall, or approaches another vesicle. However, a large time step may be
taken when a vesicle has a smooth boundary and is separated from all
other vesicles and the boundary walls.  These considerations
necessitate that we use an adaptive time-stepping scheme mostly for
robustness of the code and to remove the need to manually select a time
step size.

Another challenge in simulating vesicle flows is that maintaining good
accuracy for long horizons requires a great number of time steps when
using a low-order time-stepping method. Higher-order methods in time
can mitigate error accumulation over long simulation times.  However,
implicit higher-order methods are quite challenging to combine with
adaptive schemes. The most common methods are implicit multistep
schemes, but those are problematic especially for dynamics that involve
moving interfaces.

In summary, the design goals for a time-stepping scheme for vesicle
flows is to address stiffness with reasonable computational costs,
allow for adaptivity, and enable high-order time marching (at least
second-order). In the past (in other groups as well as our group),
methodologies have addressed parts of the design goals above, but no
method, to our knowledge, address all of them.  In this paper, we
propose a scheme that provides this capability.

\paragraph{Summary of the method and contributions}
To address stiffness, we use our recent work in which wall-vesicle,
vesicle-vesicle, and bending-tension self interactions are treated
implicitly~\cite{qua:bir2013b}.  There we used backward difference
formulas (BDFs) with fixed time step sizes.  As we discussed, it is
possible to extend BDF with adaptive time stepping, but, since BDFs
require the solution from multiple previous time step sizes but higher
order methods become quite difficult to use in
practice~\cite{dah:lin:nev1983}. We develop a new high-order method
that uses spectral deferred correction to iteratively increase the
order of a first-order time integrator.  By introducing SDC, we are
able to iteratively construct high-order solutions that only require
one previous time step.

In addition, we propose a time-stepping error estimation specific to
vesicle flows.  The incompressibility and inextensibility conditions
require that the vesicle preserves both its enclosed area and total
length. Thus, errors in area enclosed by a vesicle and errors in its
perimeter can be used estimate the local truncation error. In this way
we avoid forming multiple expensive solutions to estimate the local
truncation error.

Our contributions are summarized below.
\begin{itemize}
  \item We propose an SDC formulation to construct high-order solutions of an
  integro-differential-algebraic equation that governs vesicle
  suspensions.

  \item We propose an adaptive time stepping method that uses conserved quantities
  to estimate the local truncation error, therefore not requiring
  multiple numerical solutions to form this estimate.

  \item We conduct numerical experiments for several different
  flows involving multiple vesicles and confined flows that
  demonstrate the behavior of our scheme.
\end{itemize}

\paragraph{Related work}
There is a rich literature on numerical methods for Stokesian
particulate flows. We only discuss some representative time-stepping
methods for vesicle flows and we omit details on the spatial
discretization. Most methods for vesicle flows are based on integral
equation
formulations~\cite{ram:poz1998,soh:tse:li:voi:low2010,suk:sei2001,zha:sha2013b,rah:las:vee:cha:mal:moo:sam:shr:vet:vud:zor:bir2010,
rah:vee:bir2010,vee:gue:zor:bir2009,vee:rah:bir:zor2011,zha:sha2011a},
but other formulations based on stencil discretizations also
exist~\cite{laa:sar:mis2014,bib:kas:mis2005,du:zha2007,kim:lai2010}.

The time-stepping schemes used for vesicle flows can be grouped in
three classes of methods. The first class is fully explicit methods
with a penalty formulation for the inextensibility constraint. These
are relatively easy to implement but inefficient due to bending
stiffness. Another class of methods treats self-interactions using
linear semi-implicit time stepping. This addresses one main source of
stiffness but they are fragile. No adaptive or higher than second-order
scheme have been reported for these schemes. Finally, a third class of
methods treat all the vesicle interactions using linearization and
semi-implicit stepping. We are aware of two papers in this direction.
In \cite{zha:sha2013b} a first order backward Euler scheme is used
(with no adaptivity). As we mentioned in our own
work~\cite{qua:bir2013b,rah:vee:zor:bir2012}, we used a BDF scheme for vesicle
flows that treats all interactions implicitly, but again, it is not
adaptive and is accurate only up to second-order.

Finally, we briefly review the literature on SDC methods.  SDC was
first introduced by Dutt et al.~\cite{dut:gre:rok2000}.  SDC has been
applied to partial differential
equations~\cite{jia:hua2008,spe:rup:emm:bol:kra2013}, differential
algebraic equations~\cite{bu:hua:min2012,hua:jia:min2007}, and
integro-differential equations~\cite{hua:lai:xia2006}, but not vesicle
flows.  In~\cite{min2003}, Minion extended SDC to implicit-explicit
(IMEX) methods by studying the problem $\dot{\xx}(t) = F_{E}(\xx,t) +
F_{I}(\xx,t)$ where $F_{E}$ is non-stiff and $F_{I}$ is stiff.
Unfortunately, our governing equations do not exhibit an additive
splitting as non-stiff and stiff terms.

Algorithmically, SDC has been accelerated using
preconditioners~\cite{hua:jia:min2006,hua:jia:min2007} and parallel
algorithms~\cite{chr:ong2010,chr:mac:ong2010,chr:hay:ong2012,ong:mel:chr2012,spe:rup:emm:bol:kra2013}.
In our work, the SDC iterations are formed with first-order time
methods.  Higher-order methods can be used to form the iterations, and
this is analyzed in~\cite{chr:ong:qui2009}.  However, such an analysis
for integro-differential-algebraic equations is not the focus of this work.

\paragraph{Limitations}
The main limitations of our approach are the following.
\begin{itemize}
  \item We cannot provide a proof on the convergence order of our
  scheme.  We will see in Section~\ref{s:results} that each SDC
  correction reduces the error significantly. However, when large
  numbers of SDC corrections are used, the asymptotic convergence
  rates are unclear.  Minion observes similar behavior for stiff
  problems in~\cite{min2003}.  Understanding the asymptotic rates of
  convergence could be used to further optimize the procedure for
  choosing optimal time step sizes.

  \item We have not explored the scheme for vesicle suspensions
  where the viscosity inside and outside each vesicles differs
  (viscosity contrast).  The difficulty arises because flows with a
  viscosity contrast include a double-layer potential of the fluid
  velocity. 

  \item Currently, each vesicle must use the same time step size.
  Therefore, in flows with multiple vesicles, the time step size is
  controlled by the vesicle requiring the smallest time step.  This can
  result in the actual global error being much less than the desired
  error (see {\em Couette} example).  Simulations could be accelerated
  if each vesicle could have its own time step size.

\end{itemize}

Let us remark that we do not use adaptivity in space.

\paragraph{Outline of the paper}
In Section~\ref{s:formulation}, we briefly discuss SDC for initial
value problems (IVPs), and then discuss how to extend SDC to vesicle
suspensions.  In Section~\ref{s:numerics}, we discretize the
differential and integral equations arising in the SDC framework,
discuss preconditioning, and provide complexity estimates.
Section~\ref{s:adaptive} discusses our adaptive time stepping strategy,
and numerical results are presented in Section~\ref{s:results}.
Finally, we make concluding remarks in Section~\ref{s:conclusions}, and
we reserve a more in depth formulation of our governing equations in
the SDC framework Appendix~\ref{a:appendix1}.

\paragraph{Notation}
In Table~\ref{t:notation}, we summarize the main notation used in this
paper.
\begin{table}[htps]
\colorbox{gray!20}{
\begin{tabular}{|l|l|}
\hline
Symbol & Definition \\
\hline
$\gamma_{k}$                & Boundary of vesicle $k$ \\
$\xx_{k}(s,t)$              & Parameterization of $\gamma_{k}$ at time
                              $t$ and parameterized in arclength $s$ \\
$\sigma_{k}(s,t)$           & Tension of vesicle $k$ at time $t$ and
                              parameterized in arclength $s$ \\
$\Gamma$                    & Boundary of the confined geometry \\
$\BB(\xx_{k})$              & Bending operator due to vesicle $k$ \\
$\TT(\xx_{k})$              & Tension operator due to vesicle $k$ \\
$\Div(\xx_{k})$             & Surface divergence operator due to 
                              vesicle $k$ \\
$\SS(\xx_{j},\xx_{k})$      & Single-layer potential due to vesicle $k$ 
                              and evaluated on vesicle $j$ \\
$\DD(\xx_{j},\Gamma)$       & Double-layer potential due to $\Gamma$ 
                              and evaluated on vesicle $j$ \\
$\vv_{\infty}(\xx_{j})$     & Background velocity due to a far-field 
                              condition \\
$\vv(\xx_{j};\xx_{k})$      & Velocity of vesicle $j$ due to
                              hydrodynamic forces from vesicle $k$ \\
$\rr(t;\txx{})$             & Residual of equation~\eqref{e:picardEqn} 
                              due to the provisional solution 
                              $\txx{}$ \\
$\exx{j}$                   & Error between the exact vesicle position and
                              the provisional position $\txx{j}$  \\
$\esigma{j}$                & Error between the exact vesicle tension and
                              the provisional tension $\tsigma{j}$  \\
$A(t),L(t)$                 & Area and length of the vesicles at time 
                              $t$ \\
$e_{A},e_{L}$               & Error in area and length \\
$\epsilon$                  & Desired final tolerance for the area and 
                              length of the vesicles \\
$\beta_{\up} \geq 1$        & Maximum amount the time step is allowed
                              to be scaled up per time step \\
$\beta_{\down} \leq 1$      & Minimum amount the time step is allowed
                              to be scaled down per time step \\
$\alpha \leq 1$             & Multiplicative safety factor for the new
                              time step size \\
\hline
\end{tabular}
}
\mcaption{Index of frequently used notation.}{t:notation}
\end{table}

\section{Formulation\label{s:formulation}} 
In this section, we briefly summarize the SDC formulation for IVPs and
then extend SDC to an integro-differential equation that governs
vesicle dynamics.  We only discuss the theory of SDC and outline its
numerical implementation in Section~\ref{s:numerics}.

\subsection{Spectral Deferred Correction}
In its original development~\cite{dut:gre:rok2000}, SDC iteratively
constructed a high-order solution of the IVP
\begin{align}
  \frac{d\xx}{dt} = f(\xx,t), \quad t \in [0,T],
  \quad \xx(0) = \xx_{0}. 
\label{e:ivp}
\end{align}
While classical deferred correction
methods discretize the time derivative in~\eqref{e:ivp}, SDC uses a
Picard integral to avoid unstable numerical differentiation.
Equation~\eqref{e:ivp} is reformulated as
\begin{align}
  \xx(t) = \xx_{0} + \int_{0}^{t} f(\xx,\tau) d\tau, \quad t \in [0,T].
\label{e:picardIVP}
\end{align}
Any time stepping scheme can be used to generate a provisional solution
$\txx{}(t)$.  Given this provisional solution $\txx{}(t)$
of~\eqref{e:picardIVP}, the residual is defined as
\begin{align}
  \rr(t;\txx{}) = \xx_{0} - \txx{}(t) + \int_{0}^{t} f(\txx{},\tau)
  d\tau, \quad t \in [0,T].
  \label{e:picardResidual}
\end{align}
The integral in~\eqref{e:picardResidual} is approximated by a
$p^{th}$-order accurate quadrature rule.  The error $\ee = \xx -
\txx{}$ satisfies
\begin{align}
  \ee(t) = \rr(t;\txx{}) + \int_{0}^{t} 
      (f(\txx{} + \ee,\tau) - f(\txx,\tau)) d\tau, \quad
      t \in [0,T],
  \label{e:picardCorrect}
\end{align}
and an approximation $\tilde{\ee}$ of $\ee$ is computed, generally
using the same numerical method used to compute $\txx{}$.  Finally, the
provisional solution is updated to $\txx{} + \tilde{\ee}$ and the
procedure is repeated as many times as desired.  The process of forming
the residual $\rr$, approximating the error $\ee$, and updating the
provisional solution is referred to as an SDC correction or iteration.


The accuracy of SDC depends on the discretization of~\eqref{e:picardIVP}
and~\eqref{e:picardCorrect} and the accuracy of the quadrature rule
required to evaluate $\rr(t;\txx{})$.  An abstract error analysis for
applying deferred correction methods to the operator equation $F(y)=0$
has been preformed in~\cite{boh:ste1984,lin1980,ske1982,ste1973}.
In~\cite{han:str2011}, this abstract framework was applied
to~\eqref{e:ivp} by defining $F(\xx) = (-\xx'(t) + f(\xx,t), -\xx(0) +
\xx_{0})$, and the main result is state in Theorem 4.2.  Here we only
summarize the result to avoid introducing additional notation.
\begin{theorem} 
\label{thm:convergence}
Let $f \in C^{\infty}$ have bounded derivatives, and $\xx$ be the unique
solution of~\eqref{e:ivp}.  Suppose that the time integrator $\phi$ is
stable.  That is, there exists some $S > 0$ which only depends on $f$
and $T$, such that for all $\yy, \zz \in \RR^{m}$,
\begin{align*}
  \|\yy - \zz\| \leq S\|\phi_{m}(\yy) - \phi_{m}(\zz)\|.
\end{align*}
Moreover, suppose that $\phi_{m}$ is of order $k$ meaning that
$\phi_{m}(\xx) = \mathcal{O}(\Delta t^{k})$, where $\Delta t = T/m$.
Suppose that a numerical solution $\txx{}$ of $\xx$ satisfies
\begin{align*}
  \|\txx{}(T) - \xx(T)\| \leq C \Delta t^{\ell},
\end{align*}
where $C$ depends only on the derivatives of $f$ and on $T$, but not on
$m$.  Assuming the residual $\rr$ is computed exactly, if $\tee{}$ is
formed with the order $p$ time integrator $\phi$ to approximate $\ee$,
then
\begin{align*}
  \|\tee{}(T) + \txx{}(T) - \xx(T)\| \leq C \Delta t^{\ell+k}.
\end{align*}
However, since $\rr$ is approximated with a $p^{th}$-order quadrature,
the asymptotic error is
\begin{align*}
  \|\tee{}(T) + \txx{}(T) - \xx(T)\| = \bigO(\Delta t^{\min(\ell+k,p)}).
\end{align*}

\end{theorem}

Theorem~\ref{thm:convergence} tells us that by estimating the error
$\ee$ with a first-order method, which is the only order we consider in
the SDC framework, the order of accuracy is increased by one, with the
constraint that this convergence is limited by the accuracy of the
quadrature rule for approximating~\eqref{e:picardResidual}.  However,
the theorem states nothing about the stability of SDC.
In~\cite{dut:gre:rok2000}, the authors consider three discretizations
of~\eqref{e:picardIVP} and~\eqref{e:picardCorrect}: fully explicit,
fully implicit, and linear combinations of the two.  We do not consider
explicit methods since the governing equations are stiff.  The other two
methods could be used for vesicle suspensions, but both would require
solving non-linear equations.

We have successfully used a variant of IMEX methods for vesicle
suspensions~\cite{qua:bir2013b,rah:vee:bir2010,vee:gue:zor:bir2009}, and
we couple these methods with SDC in this work.  IMEX
methods~\cite{asc:ruu:wet1995} are a family of time integrators that
treat some terms (generally, linear) implicitly and other terms
explicitly.  IMEX methods for additive splittings $\dot{\xx}(t) =
F_{\mathrm{EXPLICIT}}(\xx,t) + F_{\mathrm{IMPLICIT}}(\xx,t) =
F_{E}(\xx,t) + F_{I}(\xx,t)$ of~\eqref{e:ivp} were first applied to SDC
by Minion in~\cite{min2003}.  We summarize the numerical results
from~\cite{min2003} since their behaviour resembles results that we will
observe.  First, the Van der Pol oscillator is considered in a
non-stiff, mildly stiff, and stiff regime, and the number of SDC
iterations ranges from two to six.  In the non-stiff regime, the error
behaves according to Theorem~\ref{thm:convergence}.  In the mildly stiff
case, the correct asymptotic result is observed, but not until $\Delta
t$ is much smaller.  In the stiff regime, the convergence behaviour
differs considerably from the formal error.  This behaviour is
attributed to order reduction which is further analyzed.  The author
proceeds to claim that ``the correct asymptotic convergence rates would
be observed given sufficiently small $\Delta t$; however, this is not
the relevant issue in most applications".  Two other examples considered
are a system of differential equations resulting from a spatial
discretization of an advection-diffusion equation, and Burgers equation.
Convergence rates up to fifth-order are achieved.  As before, as the
system becomes stiffer, smaller time steps are required before the
expected convergence behaviour is observed.

\subsection{SDC for Vesicle Suspensions}

Let $\{\gamma_{j}\}_{j=1}^{M}$ be a collection of vesicles
parameterized by $\xx_{j}$ and with tension $\sigma_{j}$
(Figure~\ref{f:schematics}).  We use the integro-differential equation
from~\cite{vee:gue:zor:bir2009} to balance the bending and tension
forces of the vesicle with the stress jump of the fluid across the
vesicle membrane, and an algebraic constraint to enforce the
inextensibility condition.  We start by defining the velocity of
vesicle $j$ due to the hydrodynamic forces from vesicle $k$
\begin{align*}
  \vv(\xx_{j};\xx_{k}) = \vv_{\infty}(\xx_{j})\delta_{j,k} + \SS(\xx_{j},\xx_{k})
    (-\BB(\xx_{k})(\xx_{k}) + \TT(\xx_{k})(\sigma_{k})),
\end{align*}
where $\delta_{j,k}$ is the Kronecker delta function,
\begin{align*}
  &s = \|\xx'\|, \quad \quad \rho = \|\xx - \yy\|, \\
  &\SS(\xx_{j},\xx_{k})(\ff) = \frac{1}{4\pi}\int_{\gamma_{k}} \left(
    -\log \rho + \frac{(\xx - \yy) \otimes (\xx - \yy)}{\rho^{2}} \right)
    \ff(\yy) ds_{\yy}, \quad \xx \in \gamma_{j}, \\
  &\BB(\xx)(\ff) = \frac{d^{4}\ff}{ds^{4}}, \\ 
  &\TT(\xx)(\sigma) =
    \frac{d}{ds}\left(\sigma\frac{d\xx}{ds}\right),
\end{align*}
and $\vv_{\infty}$ is the background velocity (unconfined flows) or the
velocity due to solid walls (confined flows).  In the case of confined
flows, we use the double-layer potential of an unknown density function
$\eeta$ defined on the boundary of the solid walls.  The extra equation
comes from a non-slip boundary condition on the solid walls and the
details are presented in~\cite{rah:vee:bir2010}.  We point out that
$\vv$ is not symmetric, meaning that $\vv(\xx_{j};\xx_{k}) \neq
\vv(\xx_{k};\xx_{j})$ for all $j \neq k$, and $\SS$, $\BB$, $\TT$ are
all linear in their second argument.

The notation we are using is chosen so that terms such as
$\BB(\xx^{N})(\xx^{N+1})$, which approximates
$\BB(\xx^{N+1})(\xx^{N+1})$, are understood to mean
\begin{align*}
  \BB(\xx^{N})(\xx^{N+1}) = \frac{d^{4}}{ds^{4}} \xx^{N+1},
    \quad s = \|\xx'^{N}\|,
\end{align*}
where $\xx'$ is the derivative of $\xx$ with respect to its
parameterization variable.  The tension $\sigma_{j}$ acts as a Lagrange
multiplier to satisfy the inextensibility constraint
\begin{align}
  \Div(\xx_{j})\left(\sum_{k=1}^{M}\vv(\xx_{j};\xx_{k})\right) = 0,
  \label{e:inextens}
\end{align}
where
\begin{align*}
  \Div(\xx)(\ff) = \frac{d\xx}{ds} \cdot \frac{d\ff}{ds}, \quad s =
  \|\xx'\|,
\end{align*}
which is also linear in its second argument.
Equation~\eqref{e:inextens} can be eliminated using the Schur complement
of the tension to write $\sigma_{j}$ in terms of the positions
$\xx_{k}$, $k=1,\ldots,M$.  Then, $\vv(\xx_{j};\xx_{k})$ can be written
entirely in terms of $\xx_{j}$ and $\xx_{k}$, and the no-slip boundary
condition of vesicle $j$ gives
\begin{align}
  \frac{d\xx_{j}}{dt} = \sum_{k=1}^{M}\vv(\xx_{j};\xx_{k}), \quad
  j=1,\ldots,M.
  \label{e:diffEqn}
\end{align}
The formulation~\eqref{e:diffEqn} is easiest to present our numerical
methods, but we in fact do not eliminate the tension in our
implementation.  The resulting changes to the SDC formulation are
presented in Appendix~\ref{a:appendix1}.

\begin{figure}[htps]
\begin{center}
  \begin{tabular}{cc}
  \ifTikz
  \scalebox{0.9}{
%
%
%
\begin{tikzpicture}[scale=1.0]

\begin{axis}[
  xmin = -3,
  xmax = 18,
  ymin = -6,
  ymax = 6,
  axis equal = true,
  hide axis,
  ]

\addplot [mark=none,red,line width=1.5] table{
2.4625e+00 2.4625e+00
2.3710e+00 2.5302e+00
2.2567e+00 2.5736e+00
2.1207e+00 2.5922e+00
1.9642e+00 2.5859e+00
1.7889e+00 2.5546e+00
1.5963e+00 2.4987e+00
1.3883e+00 2.4188e+00
1.1669e+00 2.3155e+00
9.3436e-01 2.1900e+00
6.9278e-01 2.0434e+00
4.4453e-01 1.8771e+00
1.9199e-01 1.6927e+00
-6.2390e-02 1.4920e+00
-3.1617e-01 1.2770e+00
-5.6691e-01 1.0496e+00
-8.1219e-01 8.1219e-01
-1.0496e+00 5.6691e-01
-1.2770e+00 3.1617e-01
-1.4920e+00 6.2390e-02
-1.6927e+00 -1.9199e-01
-1.8771e+00 -4.4453e-01
-2.0434e+00 -6.9278e-01
-2.1900e+00 -9.3436e-01
-2.3155e+00 -1.1669e+00
-2.4188e+00 -1.3883e+00
-2.4987e+00 -1.5963e+00
-2.5546e+00 -1.7889e+00
-2.5859e+00 -1.9642e+00
-2.5922e+00 -2.1207e+00
-2.5736e+00 -2.2567e+00
-2.5302e+00 -2.3710e+00
-2.4625e+00 -2.4625e+00
-2.3710e+00 -2.5302e+00
-2.2567e+00 -2.5736e+00
-2.1207e+00 -2.5922e+00
-1.9642e+00 -2.5859e+00
-1.7889e+00 -2.5546e+00
-1.5963e+00 -2.4987e+00
-1.3883e+00 -2.4188e+00
-1.1669e+00 -2.3155e+00
-9.3436e-01 -2.1900e+00
-6.9278e-01 -2.0434e+00
-4.4453e-01 -1.8771e+00
-1.9199e-01 -1.6927e+00
6.2390e-02 -1.4920e+00
3.1617e-01 -1.2770e+00
5.6691e-01 -1.0496e+00
8.1219e-01 -8.1219e-01
1.0496e+00 -5.6691e-01
1.2770e+00 -3.1617e-01
1.4920e+00 -6.2390e-02
1.6927e+00 1.9199e-01
1.8771e+00 4.4453e-01
2.0434e+00 6.9278e-01
2.1900e+00 9.3436e-01
2.3155e+00 1.1669e+00
2.4188e+00 1.3883e+00
2.4987e+00 1.5963e+00
2.5546e+00 1.7889e+00
2.5859e+00 1.9642e+00
2.5922e+00 2.1207e+00
2.5736e+00 2.2567e+00
2.5302e+00 2.3710e+00
2.4625e+00 2.4625e+00
};

\addplot [mark=none,red,line width=1.5] table{
5.9625e+00 2.4625e+00
5.8710e+00 2.5302e+00
5.7567e+00 2.5736e+00
5.6207e+00 2.5922e+00
5.4642e+00 2.5859e+00
5.2889e+00 2.5546e+00
5.0963e+00 2.4987e+00
4.8883e+00 2.4188e+00
4.6669e+00 2.3155e+00
4.4344e+00 2.1900e+00
4.1928e+00 2.0434e+00
3.9445e+00 1.8771e+00
3.6920e+00 1.6927e+00
3.4376e+00 1.4920e+00
3.1838e+00 1.2770e+00
2.9331e+00 1.0496e+00
2.6878e+00 8.1219e-01
2.4504e+00 5.6691e-01
2.2230e+00 3.1617e-01
2.0080e+00 6.2390e-02
1.8073e+00 -1.9199e-01
1.6229e+00 -4.4453e-01
1.4566e+00 -6.9278e-01
1.3100e+00 -9.3436e-01
1.1845e+00 -1.1669e+00
1.0812e+00 -1.3883e+00
1.0013e+00 -1.5963e+00
9.4541e-01 -1.7889e+00
9.1414e-01 -1.9642e+00
9.0778e-01 -2.1207e+00
9.2638e-01 -2.2567e+00
9.6976e-01 -2.3710e+00
1.0375e+00 -2.4625e+00
1.1290e+00 -2.5302e+00
1.2433e+00 -2.5736e+00
1.3793e+00 -2.5922e+00
1.5358e+00 -2.5859e+00
1.7111e+00 -2.5546e+00
1.9037e+00 -2.4987e+00
2.1117e+00 -2.4188e+00
2.3331e+00 -2.3155e+00
2.5656e+00 -2.1900e+00
2.8072e+00 -2.0434e+00
3.0555e+00 -1.8771e+00
3.3080e+00 -1.6927e+00
3.5624e+00 -1.4920e+00
3.8162e+00 -1.2770e+00
4.0669e+00 -1.0496e+00
4.3122e+00 -8.1219e-01
4.5496e+00 -5.6691e-01
4.7770e+00 -3.1617e-01
4.9920e+00 -6.2390e-02
5.1927e+00 1.9199e-01
5.3771e+00 4.4453e-01
5.5434e+00 6.9278e-01
5.6900e+00 9.3436e-01
5.8155e+00 1.1669e+00
5.9188e+00 1.3883e+00
5.9987e+00 1.5963e+00
6.0546e+00 1.7889e+00
6.0859e+00 1.9642e+00
6.0922e+00 2.1207e+00
6.0736e+00 2.2567e+00
6.0302e+00 2.3710e+00
5.9625e+00 2.4625e+00
};

\addplot [mark=none,red,line width=1.5] table{
1.2062e+01 2.4625e+00
1.1971e+01 2.5302e+00
1.1857e+01 2.5736e+00
1.1721e+01 2.5922e+00
1.1564e+01 2.5859e+00
1.1389e+01 2.5546e+00
1.1196e+01 2.4987e+00
1.0988e+01 2.4188e+00
1.0767e+01 2.3155e+00
1.0534e+01 2.1900e+00
1.0293e+01 2.0434e+00
1.0045e+01 1.8771e+00
9.7920e+00 1.6927e+00
9.5376e+00 1.4920e+00
9.2838e+00 1.2770e+00
9.0331e+00 1.0496e+00
8.7878e+00 8.1219e-01
8.5504e+00 5.6691e-01
8.3230e+00 3.1617e-01
8.1080e+00 6.2390e-02
7.9073e+00 -1.9199e-01
7.7229e+00 -4.4453e-01
7.5566e+00 -6.9278e-01
7.4100e+00 -9.3436e-01
7.2845e+00 -1.1669e+00
7.1812e+00 -1.3883e+00
7.1013e+00 -1.5963e+00
7.0454e+00 -1.7889e+00
7.0141e+00 -1.9642e+00
7.0078e+00 -2.1207e+00
7.0264e+00 -2.2567e+00
7.0698e+00 -2.3710e+00
7.1375e+00 -2.4625e+00
7.2290e+00 -2.5302e+00
7.3433e+00 -2.5736e+00
7.4793e+00 -2.5922e+00
7.6358e+00 -2.5859e+00
7.8111e+00 -2.5546e+00
8.0037e+00 -2.4987e+00
8.2117e+00 -2.4188e+00
8.4331e+00 -2.3155e+00
8.6656e+00 -2.1900e+00
8.9072e+00 -2.0434e+00
9.1555e+00 -1.8771e+00
9.4080e+00 -1.6927e+00
9.6624e+00 -1.4920e+00
9.9162e+00 -1.2770e+00
1.0167e+01 -1.0496e+00
1.0412e+01 -8.1219e-01
1.0650e+01 -5.6691e-01
1.0877e+01 -3.1617e-01
1.1092e+01 -6.2390e-02
1.1293e+01 1.9199e-01
1.1477e+01 4.4453e-01
1.1643e+01 6.9278e-01
1.1790e+01 9.3436e-01
1.1916e+01 1.1669e+00
1.2019e+01 1.3883e+00
1.2099e+01 1.5963e+00
1.2155e+01 1.7889e+00
1.2186e+01 1.9642e+00
1.2192e+01 2.1207e+00
1.2174e+01 2.2567e+00
1.2130e+01 2.3710e+00
1.2062e+01 2.4625e+00
};

  \foreach \y in {-3,-2,-1,1,2,3}
     \addplot[color=black,line width = 1.0pt,solid,->]
       plot coordinates{(15,\y)(15+\y,\y)};

\end{axis}

\node[font = \Huge] at (1.4,4.3) {$\gamma_{1}$};
\node[font = \Huge] at (2.8,4.3) {$\gamma_{2}$};
\node[font = \Huge] at (3.9,4.3) {$\cdots$};
\node[font = \Huge] at (5.0,4.3) {$\gamma_{_{M}}$};

\node[font = \Large] at (7.8,2.7) {$\mathbf{v}_{\infty} = (y,0)$};

\end{tikzpicture}

} & 
  \scalebox{0.9}{\input{boundedGeom.tikz}} 
  \fi
  \end{tabular}
\end{center}
\mcaption{Two typical vesicle suspensions.  Left: $M$ vesicles are
submerged in an unbounded shear flow.  Right: $13$ vesicles are in a
bounded domain.  In the right configuration, the flow is driven by
Dirichlet boundary conditions on $\Gamma$ (The internal cylinders are
rotating.)}{f:schematics}
\end{figure}

Following the SDC method we reformulate~\eqref{e:diffEqn} as
\begin{align}
  \xx_{j}(t) = \xx_{j}(0) + \int_{0}^{t}
    \sum_{k=1}^{M}\vv(\xx_{j};\xx_{k})d\tau,
    \quad t \in [0,T].
  \label{e:picardEqn}
\end{align}
We form a provisional solution $\txx{}$ at a set of quadrature nodes in
$[0,\Delta t]$ using a method that we describe in the next section
(Section~\ref{s:vesiclePicard}), and then evaluate
\begin{align}
  \rr_{j}(t;\txx{}) = \xx_{j}(0) - \txx{j}(t) + \int_{0}^{t}
    \sum_{k=1}^{M}\vv(\txx{j};\txx{k})d\tau,
    \quad t \in [0,T],
  \label{e:picardRes}
\end{align}
with a quadrature rule that we also define in the next section
(Section~\ref{s:vesicleQuadrature}).  Then, the error $\exx{j} = \xx_{j}
- \txx{j}$ satisfies
\begin{align*}
  \exx{j}(t) = \xx_{j}(0) - \txx{j}(t) + \int_{0}^{t}
    \sum_{k=1}^{M}\vv(\txx{j} + \exx{j};\txx{k} + \exx{k})
    d\tau, \quad t \in [0,T],
\end{align*}
which we write using the residual $\rr$ as
\begin{align}
  \exx{j}(t) = \rr_{j}(t;\txx{}) + \int_{0}^{t}
  \sum_{k=1}^{M}(\vv(\txx{j}+\exx{j};\txx{k} + \exx{k}) - 
      \vv(\txx{j};\txx{k})) d\tau,
    \quad t \in [0,T].
  \label{e:sdcUpdate}
\end{align}
We define the new provisional solution as $\txx{} +
\tilde{\ee}_{\xx}$, where $\tilde{\ee}_{\xx}$ is an approximate
solution of~\eqref{e:sdcUpdate}.  Again, this procedure of computing the
residual~\eqref{e:picardRes}, numerically solving~\eqref{e:sdcUpdate}
for the error, and updating the provisional solution is what we call an
SDC iteration.  Assuming that we take $n_{\sdc}$ first-order SDC
iterations and $p$ is the quadrature error for computing the residual
$\rr$, from Theorem~\ref{thm:convergence}, we expect that the asymptotic
rate of convergence is $\min(n_{\sdc},p)$.

\section{Numerical Scheme\label{s:numerics}}
In this section, we present our numerical scheme for
solving~\eqref{e:picardEqn},~\eqref{e:picardRes},
and~\eqref{e:sdcUpdate}.  Evaluating $\vv$ involves layer potentials
and Fourier differentiation, and for concentrated suspensions, this can
be costly.  This cost will be further amplified in three-dimensional
simulations (which are not discussed in this work).  Therefore, we
discretize~\eqref{e:sdcUpdate} differently from standard methods of
discretizing~\eqref{e:picardCorrect}.  We also discuss how the
block-diagonal preconditioner outlined in~\cite{qua:bir2013b} is
applied.  We conclude with estimates of the overall work per time step
as a function of the number of vesicles and the number of points per
vesicle.

\subsection{Spatial Discretization}

Following our previous
work~\cite{qua:bir2013b,rah:vee:bir2010,vee:gue:zor:bir2009}, we use a
Lagrangian formulation by marking each vesicle with $N$ tracker
points.  Since vesicles are modeled as smooth closed curves, all the
derivatives are computed spectrally using the fast Fourier transform
(FFT).  Near-singular integrals are handled using the near-singular
integration scheme outlined in~\cite{qua:bir2013b}.  Finally, we use
Alpert's high-order Gauss-trapezoid quadrature rule~\cite{alp1999} with
accuracy $\bigO(h^{8} \log h)$ to evaluate the single-layer potential,
and the trapezoid rule for the double-layer potential which has
spectral accuracy since its kernel is smooth.

\subsection{Temporal Discretization}

A fully explicit discretization of~\eqref{e:diffEqn} results in a stiff
system for multiple reasons.  A stiffness analysis
in~\cite{vee:gue:zor:bir2009} reveals that the leading sources of
stiffness, and the corresponding time step restrictions, are:
\begin{itemize}
  \item $\SS(\xx_{j},\xx_{j})\BB(\xx_{j})(\xx_{j})$ (self-hydrodynamic
  bending force) -- $\Delta t \sim \Delta s^{3}$;

  \item $\SS(\xx_{j},\xx_{j})\TT(\xx_{j})(\xx_{j})$ (self-hydrodynamic
  tension force) -- $\Delta t \sim \Delta s$;

  \item $\Div(\xx_{j})(\xx_{j})$ (self-inextensibility force) -- $\Delta
  t \sim \Delta s$;

  \item $\SS(\xx_{j},\xx_{k})(\BB(\xx_{k})(\xx_{k}) +
  \TT(\xx_{k})(\sigma_{k}))$, $j \neq k$ (inter-vesicle hydrodynamic
  forces) -- depends on the inter-vesicle distance.
\end{itemize}
The leading sources of stiffness result from the intra-vesicle
interactions; but, for concentrated suspensions, inter-vesicle
interactions become significant and introduce stiffness.

To address these multiple sources of stiffness, we use a variant of
IMEX~\cite{asc:ruu:wet1995} time integrators.  IMEX methods were
developed for the problem $\dot{\xx}(t) = F_{E}(\xx) + F_{I}(\xx)$,
where $F_{E}$ is non-stiff and is treated explicitly, and $F_{I}$ is
stiff and is treated implicitly.  For problems with this additive
splitting, the family of time integrators is
\begin{align*}
  \frac{\beta \xx^{n+1} - \xx^{0}}{\Delta t} = 
    F_{E}(\xx^{e}) + F_{I}(\xx^{n+1}),
\end{align*}
where $\xx^{0}$ and $\xx^{e}$ are linear combinations of previous time
steps and $\beta > 0$.  Unfortunately,~\eqref{e:diffEqn} does not have
an additive splitting between stiff and non-stiff terms.  However, we
have observed first- and second-order convergence~\cite{qua:bir2013b}
for the time integrator
\begin{align*}
  \frac{\beta \xx^{n+1}_{j} - \xx^{0}_{j}}{\Delta t} =
    \sum_{k=1}^{M} \vv(\xx_{j}^{e};\xx_{k}^{n+1}), \quad j=1,\dots,M,
\end{align*}
where
\begin{align*}
  \vv(\xx_{j}^{e};\xx_{k}^{n+1}) = \SS(\xx^{e}_{j},\xx^{e}_{k})
    (-\BB(\xx^{e}_{k})(\xx^{n+1}_{k}) +
      \TT(\xx^{e}_{k})(\sigma^{n+1}_{k})).
\end{align*}
Note that in this formulation, it is the operators involved in $\vv$,
such as the bending $\BB(\xx_{k})$ and the single-layer potential
$\SS(\xx_{j},\xx_{k})$, that are discretized explicitly at $\xx^{e}$.
In line with our previous work, we use the first-order method given by
$\beta = 1$, $\xx^{0}=\xx^{n}$, and $\xx^{e}=\xx^{n}$, and the
second-order backward difference formula (BDF) given by $\beta = 3/2$,
$\xx^{0} = 2\xx^{n} - 1/2\xx^{n-1}$, and $\xx^{e} = 2\xx^{n} -
\xx^{n-1}$.  

\subsection{Quadrature Formula}
\label{s:vesicleQuadrature}

SDC requires a quadrature formula to approximate~\eqref{e:picardRes},
the residual of the Picard integral.  The quadrature rule effects the
accuracy of SDC (Theorem~\ref{thm:convergence}).  However, the
stability of SDC is also effected by the quadrature rule, and this has
been investigated in depth in~\cite{bou:lay:min2003}.  We do not use
equi-spaced nodes since they are vulnerable to the Runge phenomenon.
In order to achieve the highest possible accuracy, we could use
Gaussian nodes; however, in order to avoid extrapolation, we would like
both endpoints to be quadrature nodes.  Therefore, we use Gauss-Lobatto
points, $0 = t_{1} < \cdots < t_{p}=\Delta t$, since they include both
endpoints, have maximal order $2p-3$, and have successfully been used
by other groups~\cite{bou:lay:min2003,bou:min2010,min2003}.
Alternatively, we could use Radau quadrature formula which include the
left endpoint, but not the right.

For computational efficiency, the accuracy of the quadrature formula
should not exceed the expected rate of convergence.  This is especially
important in three dimensions, which is not considered in this work,
since multiple variables must be stored at each of the quadrature nodes
for future SDC iterations.

\subsection{Picard Integral Discretization}
\label{s:vesiclePicard}

Equations~\eqref{e:picardEqn} and~\eqref{e:sdcUpdate} have similar
structure, and we take advantage of this structure in their
discretizations.  For adaptivity, we want a scheme that easily allows
for variable time step sizes.  For this reason, we use only first-order
methods for~\eqref{e:picardEqn} and~\eqref{e:sdcUpdate}.  When desired,
we use SDC iterations to increase the accuracy.

The first-order provisional solution is found by
discretizing~\eqref{e:picardEqn} as
\begin{align}
  \xx_{j}^{n+1} = \xx_{j}^{n} + \Delta t_{n}
    \sum_{k=1}^{M} \vv(\xx_{j}^{n};\xx_{k}^{n+1}), \quad
    n=0,\ldots,p-1,
  \label{e:numericProvisionalImplicit}
\end{align}
where $\Delta t_{n} = t_{n+1} - t_{n}$.  This is exactly the first-order
time integrator we introduced in~\cite{qua:bir2013b}.  As we march in
time with~\eqref{e:numericProvisionalImplicit}, we save the variables
required for near-singular integration for future SDC iterations.  Then,
we evaluate the residual~\eqref{e:picardRes} using the Gauss-Lobatto
quadrature rule.

In line with~\cite{min2003}, which considers semi-implicit SDC methods,
we would like to discretize~\eqref{e:sdcUpdate} as
\begin{align*}
  \ee^{n+1}_{\xx_{j}} = \ee^{n}_{\xx_{j}} + 
  \rr_{j}^{n+1} - \rr_{j}^{n} + 
  \Delta t_{n}\sum_{k=1}^{M}(
    \vv(\txx{j}^{n}+\ee^{n}_{\xx_{j}};\txx{k}^{n+1}+\ee^{n+1}_{\xx_{k}}) - 
    \vv(\txx{j}^{n},\txx{k}^{n+1})).
\end{align*}
The issue with this formulation is that it requires additional storage
and computations to find the velocity due to the vesicle parameterized
by $\txx{k}^{n} + \ee^{n}_{\xx_{k}}$.  Again, in three dimensions
this restriction is even more prohibitive.  We have experimented with
other discretizations of~\eqref{e:sdcUpdate}.  The simplest such one is
\begin{align*}
  \ee^{n+1}_{\xx_{j}} = \ee^{n}_{\xx_{j}} + 
    \rr_{j}^{n+1} - \rr_{j}^{n}.
\end{align*}
This discretization is consistent with the governing
equations\footnote{If $\ee^{n}_{\xx_{j}}$ converges to 0, then
$\rr_{j}^{n+1}=\rr_{j}^{n}$, and by~\eqref{e:picardRes}, we have
solved~\eqref{e:picardEqn} up to quadrature error.}, but,
experimentally, SDC converges only if a very small time step is used.
To allow for larger time steps, we can include the implicit term in the
discretization
\begin{align*}
  \ee^{n+1}_{\xx_{j}} = \ee^{n}_{\xx_{j}} + 
  \rr_{j}^{n+1} - \rr_{j}^{n} + 
  \Delta t_{n} \sum_{k=1}^{M} 
  (\vv(\txx{j}^{n};\txx{k}^{n+1}+\ee_{\xx_{k}}^{n+1}) -
    \vv(\txx{j}^{n};\txx{k}^{n+1})),
\end{align*}
where we use a slight abuse of notation by defining
\begin{align*}
  \vv(\txx{j}^{n};\txx{k}^{n+1} + \ee_{\xx_{k}}^{n+1}) = 
    \SS(\txx{j}^{n},\txx{k}^{n})
    (-\BB(\txx{k}^{n})(\txx{k}^{n+1}+\ee_{\xx_{k}}^{n+1}) +
      \TT(\txx{k}^{n})(\sigma^{n+1}_{k})).
\end{align*}
(Note how none of the operators depend on the error $\ee$.)
Appendix~ref{a:appendix1} presents this same discretization without the
abuse of notation.  While this discretization allows larger time steps,
it is incompatible with the inextensibility constraint (see
Appendix~\ref{a:appendix1}).  A discretization that is is compatible
with the inextensibility constraint is
\begin{align}
  \ee^{n+1}_{\xx_{j}} = \ee^{n}_{\xx_{j}} + 
  \rr_{j}^{n+1} - \rr_{j}^{n} + 
  \Delta t_{n} \sum_{k=1}^{M} 
  (\vv(\txx{j}^{n+1};\txx{k}^{n+1}+\ee_{\xx_{k}}^{n+1}) -
    \vv(\txx{j}^{n+1};\txx{k}^{n+1})),
  \label{e:numericSDCImplicit}
\end{align}
where we are using the same abuse of notation.  Note
that~\eqref{e:numericSDCImplicit} only requires evaluating the velocity
field due to the vesicle configuration given by $\txx{k}^{n+1}$.  Since
these velocity fields are required to form residual $\rr$, no additional
velocity fields need to be formed.  However, there is a loss in
accuracy, but we have observed experimentally that if only
one or two SDC corrections are used, it does not reduce the rate of
convergence.

To summarize, the main steps for using SDC to solve~\eqref{e:picardEqn}
are
\begin{enumerate}
  \item Find a first-order provisional solution $\tilde{\xx}$
  using~\eqref{e:numericProvisionalImplicit}.

  \item \label{step:residual} Compute the residual $\rr$ by
  approximating the integral in~\eqref{e:picardRes} with the
  Gauss-Lobatto quadrature rule.

  \item Use~\eqref{e:numericSDCImplicit} to approximate the error
  $\tee{}$.

  \item Define the new provisional solution to be $\txx{} + \tee{}$.

  \item Go to step~\ref{step:residual}.
\end{enumerate}

\subsection{Preconditioning}
\label{s:preco}
Equations~\eqref{e:numericProvisionalImplicit}
and~\eqref{e:numericSDCImplicit} are ill-conditioned and require a large
number of GMRES iterations~\cite{vee:gue:zor:bir2009}.  To reduce the
cost of the linear solves, we used a block-diagonal preconditioner
in~\cite{qua:bir2013b} which is formed and factorized in matrix form at
each time step.  Using this preconditioner, the number of preconditioned
GMRES iterations depends only on the magnitude of the inter-vesicle
interactions, which in turn is a function of the proximity of the
vesicles.  For further savings, we freeze and factorize the
preconditioner at the first Gauss-Lobatto point, and this preconditioner
is used for all the subsequent Gauss-Lobatto points and SDC iterates.
By freezing the preconditioner, we significantly reduce the number of
GMRES iterations, and we only require one matrix factorization per time
step, which is the number of factorizations required when we
precondition our time integrators introduced in~\cite{qua:bir2013b}.

\subsection{Complexity Estimates}
Here we summarize the cost of the most expensive algorithms required in
our formulation.
\begin{itemize}

\item \emph{Matrix-vector multiplication:} For unbounded flows, if $M$
vesicles are each discretized with $N$ points, the bending and tension
calculations require $\bigO(MN\log N)$ operations using the FFT, and the
single-layer potential requires $\bigO(MN)$ operations using the FMM.
If the solid wall is discretized with $N_{\wall}$ points, using the FMM,
the matrix-vector multiplication requires $\bigO(MN\log N + N_{\wall})$
operations.  

\item \emph{Computing the residual:}
Given a provisional solution $\tilde{\xx}$, computing the residual
$\rr$ is equivalent to $p$ matrix-vector multiplications.  Therefore,
computing the residual requires $\bigO(p(MN\log N + N_{\wall}))$
operations.

\item \emph{Forming the provisional solution and SDC corrections:}
Equations~\eqref{e:numericProvisionalImplicit}
and~\eqref{e:numericSDCImplicit} require solving the same linear system
(only the right-hand sides are different), and our preconditioner
results in a mesh-independent number of GMRES iterations.  Therefore,
if $n_{\gmres}$ total iterations are required to find the provisional
solution, then $n_{\sdc}$ SDC iterations requires
$\bigO(n_{\gmres}p(n_{\sdc}+1)(MN\log N + N_{\wall}))$ operations.  

\item \emph{Forming the preconditioner:} The preconditioner is computed
and stored in matrix form and requires $\bigO(MN^{2}\log N)$ operations
per time step by using Fourier differentiation and the FFT.  The
preconditioner must be factorized which requires $\bigO(MN^{3})$
operations.  This is computed only once per time step and is reused at
all the additional Gauss-Lobatto quadrature points.  In two-dimensions,
this cost is acceptable since the number of unknowns on each vesicle is
relatively small.  However, in three-dimensions, this cost is
unacceptable and different preconditioners will need to be constructed.
\end{itemize}

\section{Adaptive Time Stepping\label{s:adaptive}}
During the course of a simulation, vesicles come close to each other and
to confining walls.  For instance, as a vesicle passes through a
constriction (see Figure~\ref{f:stenosisGeom}), it can come very close
to the solid wall and small time steps should be taken.  Resolving
multiple time scales is an important step towards a robust solver for
vesicle suspensions for two reasons.  First, the simulation is sped up
since the largest possible time step is always taken, and second, we
eliminate a trial and error procedure for finding a time step size that
results in the desired tolerance.

A common strategy of an adaptive time step method is to control an
estimate of the local truncation error.  One estimate is the difference
of two numerical solutions $\xx_{1}$ and $\xx_{2}$.  One choice for
$\xx_{1}$ and $\xx_{2}$ is the numerical solution formed from a single
time step of size $\Delta t$ and one formed from two time steps of size
$\Delta t/2$ (step-doubling).  Another choice is to use solutions
formed by two different numerical methods whose orders differ by one.
For vesicle suspensions, we can instead use two invariants to estimate
the local truncation error.  We use the errors in area and length,
which are invariant by the incompressibility and inextensibility
conditions, to estimate the local truncation error.  This estimate does
not require forming multiple numerical solutions.

We now outline how this estimate is used to accept or reject a time
step, and how it is used to select a new time step size.  Suppose we
have a single vesicle\footnote{If there are multiple vesicles, we
choose the minimum requested time step size over all of the vesicles.}
at time $t$ with area $A(t)$, length $L(t)$, and the desired tolerance
for the global error is $\epsilon$ at the time horizon $T=1$.  We
compute the solution at time $t+\Delta t$ using a $k^{th}$-order time
stepping scheme which results in a new area $A(t + \Delta t)$ and
length $L(t + \Delta t)$.  First, we check if the solution at time $t +
\Delta t$ satisfies 
\begin{align}
  |A(t+\Delta t) - A(t)| \leq A(t) \Delta t \epsilon,
  \text{ and }
  |L(t+\Delta t) - L(t)| \leq L(t) \Delta t \epsilon.
  \label{e:errorBounds}
\end{align}
If condition~\eqref{e:errorBounds} is not satisfied, we reject this
solution and we compute a new solution using a smaller time step size.
If condition~\eqref{e:errorBounds} is satisfied, we accept this solution
and we increase the time step size so that we are taking the largest
possible $\Delta t$ such that~\eqref{e:errorBounds} is satisfied for
future time steps.  As an alternative to fixing $\epsilon$ for the
entire simulation, we experiment with increasing it as the simulation
progresses in Section~\ref{s:couette}.  The strategy is to increase
$\epsilon$ at each time step to account for the fact that the actual
error committed will be less than the predicted error.  In this manner,
the global error is much closer to the tolerance.

Regardless of the acceptance or rejection of the solution at time $t +
\Delta t$, a new time step size must be chosen.  We first require
estimates of the asymptotic constants of proportionality for the errors
in area and length
\begin{align*}
  C_{A} = \frac{|A(t+\Delta t) - A(t)|}{A(t)\Delta t^{k+1}}, 
  \text{ and } 
  C_{L} = \frac{|L(t+\Delta t) - L(t)|}{L(t)\Delta t^{k+1}}.
\end{align*}
We only consider the error in the area since the same argument can be
applied to the error in length.  The optimal time step size $\Delta
t_{\opt}$ satisfies
\begin{align}
  |A(t+\Delta t_{\opt}) - A(t)| = A(t) \Delta t_{\opt} \epsilon,
  \label{e:optimalEstimate1}
\end{align}
and we also have the estimate
\begin{align}
  |A(t+\Delta t_{\opt}) - A(t)| = C_{A} A(t) \Delta t_{\opt}^{k+1}.
  \label{e:optimalEstimate2}
\end{align}
Equating~\eqref{e:optimalEstimate1} and~\eqref{e:optimalEstimate2},
$\Delta t_{opt}$ satisifies 
\begin{align*}
  C_{A} \Delta t_{\opt}^{k+1} = \Delta t_{\opt} \epsilon.
\end{align*}
Finally, using our estimate for $C_{A}$, we have
\begin{align*}
  \frac{|A(t+\Delta t) - A(t)|}{A(t) \Delta t^{k+1}} 
      \Delta t_{\opt}^{k+1} = \Delta t_{\opt} \epsilon, 
\end{align*}
which implies that our next time step size should be
\begin{align*}
  \Delta t_{\opt} = \left(
    \frac{A(t) \epsilon \Delta t}{|A(t+\Delta t) - A(t)|}\right)^{1/k}
    \Delta t.
\end{align*}
A similar optimal time step size is computed based on the length and
the smaller of these two time steps is selected for the next time step.

We place restrictions on the new time step size to increase the
probability of the next time step being accepted.  Since our error
estimates are local and asymptotic, it is dangerous to change the time
step size too rapidly.  Therefore, we restrict the new time step size
scaling to the interval $[\beta_{\down},\beta_{\up}]$, where
$\beta_{\down} < 1$ and $\beta_{\up} > 1$.  Next, we multiply the new
time step size by a safety factor $\alpha^{1/k}<1$ to increase the
likelihood of the next time step size being accepted.  Finally, we
never increase the time step size if the previous time step size is
rejected~\cite{hai:nor:wan1993}.  In summary, if the previous time step
is accepted, the new time step size is
\begin{align*}
  \Delta t_{\new} = \alpha^{1/k}\min(\beta_{\up}\Delta t,
    \max(\Delta t_{\opt},\beta_{\down}\Delta t)),
\end{align*}
and if the previous time step is rejected, the new time step size is
\begin{align*}
  \Delta t_{\new} = \alpha^{1/k}\min(\Delta t,
    \max(\Delta t_{\opt},\beta_{\down}\Delta t)).
\end{align*}
The scaling factor $\alpha^{1/k}$ is chosen so that
\begin{align*}
  \max\left(\frac{|A(t+\Delta t)-A(t)|}{A(t)},
            \frac{|L(t+\Delta t)-L(t)|}{L(t)}\right)
  \approx \alpha \Delta t \epsilon,
\end{align*}
regardless of the order $k$.  That is, with this scaling of $\alpha$,
the method tries to commit the same amount of error per time step,
independent of the time stepping order.

The parameters $\alpha$, $\beta_{\up}$, and $\beta_{\down}$  affect the
overall efficiency of the method.  For instance, if $\alpha$ is too
large, then $\Delta t_{\new}$ may be too large and the time steps will
be rejected too often.  However, if it is too small, the bounds
in~\eqref{e:errorBounds} will not be tight which will increase the
total number of time steps.  We have experimented with different values
and have had success with $\alpha = \sqrt{0.9}$, $\beta_{\up} = 1.5$
and $\beta_{\down} = 0.6$.  These values are used for all the numerical
examples.  We point out that when a small tolerance $\epsilon$ is
chosen, the estimate~\eqref{e:optimalEstimate2} will be much more
accurate, and $\beta_{\up}$ and $\beta_{\down}$ will not play a role
since only small changes to the time step size will be made.  However,
when larger tolerances are requested,~\eqref{e:optimalEstimate2} is
less accurate, and without $\beta_{\up}$ and $\beta_{\down}$,
unreasonable time step sizes may be chosen.

\section{Results\label{s:results}} 
We discuss the behaviour of SDC and adaptive time stepping for confined
and unconfined flows.  We are interested in problems whose dynamics
exhibit multiple time scales, and problems with long time horizons.
For such problems, it is beneficial that an adaptive time stepping
takes the largest possible time steps and does not require a trial and
error procedure to achieve a desired tolerance.  The main parameters
are:
\begin{itemize}
  \item $N$: The number of points per vesicle;
  \item $N_{\mathrm{wall}}$: The number of points per solid wall;
  \item $m$: The number of time steps used in the non-adaptive scheme;
  \item $T$: The time horizon;
  \item $n_{\sdc}$: The number of SDC corrections;
  \item $p$: The number of Gauss-Lobatto quadrature points;
  \item $\beta_{\up}$, $\beta_{\down}$, $\alpha$: Safety factors for
  adaptive time stepping as described in Section~\ref{s:adaptive}.
\end{itemize}
The simulations are performed in Matlab on a six-core 2.67GHz Intel
Xeon processor with 24GB of memory.  When appropriate, we also report
the number of fast multiple method (FMM) calls which is the most
expensive part of the simulations.  We use four numerical examples of
increasing complexity which we now summarize.

\begin{itemize}
\item{{\bf Relaxation}
(Tables~\ref{t:noSDCrelaxation}--\ref{t:SDC34relaxation}):} We consider
stiffness due to self-interactions by simulating a single vesicle in a
Stokes fluid with no background velocity.  Motivated by
Theorem~\ref{thm:convergence}, we would like each additional SDC
iteration to result in an additional order of accuracy.

\item{{\bf Extensional}
(Tables~\ref{t:noSDCextensional}--\ref{t:adaptiveThirdOrderExtensional}
and Figures~\ref{f:extensionalGeom}--\ref{f:extensionalSummary}):} We
consider stiffness due to vesicle-vesicle interactions by simulating two
vesicles approaching each other in an extensional flow.  We check the
convergence of different time integrators and compare constant time step
sizes and adaptive time step sizes.

\item{{\bf Stenosis}
(Tables~\ref{t:SDC01stenosis}--\ref{t:adaptiveThirdOrderStenosis} and
Figures~\ref{f:stenosisGeom}--\ref{f:stenosisSummary}):} We consider
stiffness due to vesicle-wall interactions by simulating a single
vesicle in a confined tube with a parabolic-profile flow at the intake
and outtake.  The vesicle passes through a narrow region where a smaller
time step should be taken.

\item{{\bf Couette}
(Tables~\ref{t:SDC0couette}--\ref{t:adaptiveSecondOrderCouette} and
Figures~\ref{f:couetteGeom}--\ref{f:couetteSummary}):} We consider
eight vesicles in a Couette apparatus with a long time horizon.  The
solid walls are aligned so that there is a narrow region that the
vesicles pass through.  This examples involves all the stiff terms in
our model.
\end{itemize}

Other details include:
\begin{itemize}
\item When using adaptive time stepping, no trial and error procedure
to find a time step size that achieves the desired error is required.
This is especially beneficial for simulations that require large
amounts of CPU time due to long time horizons or concentrated
suspensions.

\item The area and length of the vesicles should be unchanged by the
incompressibility and inextensibility constraints.  Therefore, we use
the error in area, $e_{A}$, and the error in length, $e_{L}$, to measure
the error of the simulation, and to estimate the local truncation error
required for adaptive time stepping.

\item For the simulations we present, the interactions between the
vesicles and solid walls are the most expensive part of the calculation.
Therefore, for all examples except the first (it has no vesicle-vesicle
or vesicle-boundary interactions), we report the total number of FMM
calls (\# fmm).  

\item For more concentrated suspensions, forming and factorizing the
block-diagonal preconditioner is the most expensive part of the
caculation.  However, it is computed only once per time step in order to
gurantee that the preconditioner is formed exactly once per time step
for all the time integrators.   When SDC is used, the preconditioner
formed at the first Gauss-Lobatto point is used for all subsequent
Gauss-Lobatto points and SDC corrections.

\item For problems that use $n_{\sdc}=0$ or BDF, we do not use
Gauss-Lobatto substeps.  However, when $n_{\sdc}>0$, there are $m(p-1)$
time steps since the solution is required at the quadrature points to
estimate $\rr$.

\item We take larger values of $N$ than required in order to fully
resolve the geometry, differential operators, and the layer potentials.
This way, the temporal error, which is the main focus of this work,
dominates the spatial error.

\item All the linear systems are solved with GMRES without restarts (no
GMRES solve requires more than 50 iterations) and with a tolerance of
$1$E$-10$.  We have experimented with tighter tolerances, but this has
little or no effect on the reported errors and only increases the number
of GMRES iterations.  The reason that the errors do not decrease is that
the condition number of the linear system, which for the report problems
is on the order of $1$E$5$, restricts the minimum achievable residual.

\item For the adaptive time stepping results, the time step size is
never upscaled by more than $\beta_{\up} = 1.5$, downscaled by more
than $\beta_{\down} = 0.6$, and the safety scaling factor is $\alpha =
\sqrt{0.9}$.  In addition, the time step size is never increased
immediately after a rejected time step as is recommended
in~\cite{hai:nor:wan1993}.

\end{itemize}

\subsection{Relaxation}

We consider a single vesicle, initialized as a three-to-one ellipse, in
a Stokes fluid with no background velocity.  We discretize the vesicle
with $N=96$ points and the time horizon is $T=2$ which is large enough
that the vesicle comes within 10\% of its steady state solution.  We use
$p=5$ Gauss-Lobatto quadrature points so that the quadrature's order of
accuracy (seven in this case) is a few orders larger than the order of
the time integrator (up to five).  For this example, there are no calls
to the FMM because the self-interactions are evaluated directly.  We
report the errors and CPU timings in
Tables~\ref{t:noSDCrelaxation}--\ref{t:SDC34relaxation}.  We observe
that:
\begin{itemize}
\item Each SDC iteration significantly improves the accuracy of the
solution.  However, we are unable to achieve third- and higher-order
results.  We expect that these convergence rates will be achieved for
smaller values of $\Delta t$, but, at these required values for $\Delta
t$, other sources of error, such as the GMRES tolerance or machine
precision, will dominate.  This behaviour is observed by Minion
in~\cite{min2003}.

\item Comparing the left entires of Tables~\ref{t:noSDCrelaxation}
and~\ref{t:SDC12relaxation}, we see that the CPU time increases by more
than four-fold.  This is a result of SDC having to form the solution at
the intermediate Gauss-Lobatto points.

\item Comparing the two second-order solvers (BDF and $n_{\sdc}=1$), we
can not conclusively pick the faster method.  However, simulations using
SDC corrections are compatible with adaptive time stepping while with
BDF, they are not.

\end{itemize}

\begin{table}[htps]
\begin{center}
\begin{tabular}{c|ccc|ccc}
 & \multicolumn{3}{c|}{$\boldnsdc{0}$} & 
   \multicolumn{3}{c}{{\bf BDF}} \\
$m$ & $e_{A}$ & $e_{L}$ & CPU & $e_{A}$ & $e_{L}$ & CPU \\
\hline
$125$  & $3.44$E$-6$  & $3.40$E$-5$ & $1.0$ 
       & $6.13$E$-8$  & $1.56$E$-6$ & $1.1$ \\
$250$  & $1.74$E$-6$  & $1.74$E$-5$ & $2.1$ 
       & $1.40$E$-8$  & $3.57$E$-7$ & $2.1$ \\
$500$  & $8.76$E$-7$  & $8.83$E$-6$ & $4.2$ 
       & $3.10$E$-9$  & $7.84$E$-8$ & $4.2$ \\
$1000$ & $4.40$E$-7$  & $4.44$E$-6$ & $8.4$ 
       & $2.41$E$-10$ & $1.75$E$-8$ & $8.3$ 
\end{tabular}
\mcaption{The errors in area and length and the CPU time for a single
vesicle in a {\bf relaxation} flow with a {\bf constant} time step size
using $\boldnsdc{0}$ (left) and {\bf BDF} (right).  The CPU times for
Tables~\ref{t:noSDCrelaxation}--\ref{t:SDC34relaxation} are relative to
the cheapest simulation ($m=125$ and $n_{\sdc}=0$) which took
approximately 38 seconds.  Both methods converge with their expected
rate of convergence.}{t:noSDCrelaxation}

\begin{tabular}{c|ccc|ccc}
 & \multicolumn{3}{c|}{$\boldnsdc{1}$} & 
   \multicolumn{3}{c}{$\boldnsdc{2}$} \\
$m$ & $e_{A}$ & $e_{L}$ & CPU & $e_{A}$ & $e_{L}$ & CPU \\
\hline
$125$  & $2.75$E$-8$  & $1.30$E$-8$  & $4.6$ 
       & $3.99$E$-10$ & $5.29$E$-11$ & $7.2$ \\
$250$  & $7.28$E$-9$  & $2.58$E$-9$  & $9.0$ 
       & $1.98$E$-11$ & $6.54$E$-12$ & $16$ \\
$500$  & $1.88$E$-9$  & $4.45$E$-10$ & $17$ 
       & $1.04$E$-11$ & $6.37$E$-13$ & $29$ \\
$1000$ & $4.78$E$-10$ & $6.88$E$-11$ & $35$ 
       & $5.37$E$-12$ & $4.77$E$-14$ & $57$
\end{tabular}
\mcaption{The errors in area and length and the CPU time for a single
vesicle in a {\bf relaxation} flow with a {\bf constant} time step size
using $\boldnsdc{1}$ (left) and $\boldnsdc{2}$ (right).  We achieve
second-order convergence with one SDC correction, but third-order
convergence is only observed for the error in length with two SDC
corrections.  However, the error in area has plateaued which indicates
that the GMRES tolerance has limited the error in
area.}{t:SDC12relaxation}

\begin{tabular}{c|ccc|ccc}
 & \multicolumn{3}{c|}{$\boldnsdc{3}$} & 
   \multicolumn{3}{c}{$\boldnsdc{4}$} \\
$m$ & $e_{A}$ & $e_{L}$ & CPU & $e_{A}$ & $e_{L}$ & CPU \\
\hline
$125$  & $9.77$E$-11$ & $3.94$E$-13$ & $10$ 
       & $1.10$E$-10$ & $1.07$E$-14$ & $13$ \\
$250$  & $2.03$E$-11$ & $3.13$E$-14$ & $20$ 
       & $2.74$E$-11$ & $3.11$E$-15$ & $26$ \\
$500$  & $4.30$E$-12$ & $2.22$E$-15$ & $41$ 
       & $6.60$E$-12$ & $1.11$E$-16$ & $49$ \\
$1000$ & $7.11$E$-13$ & $1.33$E$-16$ & $80$ 
       & $1.33$E$-12$ & $1.44$E$-15$ & $101$
\end{tabular}
\mcaption{The errors in area and length and the CPU time for a single
vesicle in a {\bf relaxation} flow with a {\bf constant} time step size
using $\boldnsdc{3}$ (left) and $\boldnsdc{4}$ (right).  The error in
length achieves fourth-order convergence with three SDC corrections, but
again, the error in area has plateaued.  The error in area continues to
plateau with four SDC corrections, and the error in length has reached
machine precision.}{t:SDC34relaxation}

\end{center}
\end{table}

\subsection{Extensional}

We consider two vesicles placed symmetrically around the origin with
the background velocity $\vv_{\infty} = (-x,y)$
(Figure~\ref{f:extensionalGeom}).  We discretize both vesicles with
$N=96$ points and the time horizon is $T=24$ which is long enough that
the distance between the vesicles at the time horizon is
$0.4\sqrt{\Delta s}$, where $\Delta s$ is the arclength spacing.  We
use $p=4$ Gauss-Lobatto quadrature points so that $\rr$ has fifth-order
accuracy which is at least two orders more accurate than all the
reported time integrators.  We report results using zero, one, and two
SDC corrections, and BDF in Tables~\ref{t:noSDCextensional}
and~\ref{t:SDC12extensional}.  Again, for some of the results, the
expected convergence rates are not observed.  Other groups (see, for
example,~\cite{min2003,wei2013}) have observed that for stiff systems,
very small time steps must be taken before the asymptotic convergence
rates are achieved.  As before, we expect that other sources of error
will dominate once the temporal asymptotic regime is achieved.

\begin{figure}[htps]
\begin{center}
  \begin{tabular}{ccccc}
  \ifTikz
  \begin{tikzpicture}[scale=0.40]

\begin{axis}[
  xmin = -1.8,
  xmax = 1.8,
  ymin = -1.8,
  ymax = 1.8,
  axis equal = true,
  hide axis,
  title = {\Huge$t=0$}
  ]

\addplot [mark=none,red,line width=1.5] table{
-1.2000e+00 1.7412e+00
-1.2376e+00 1.7375e+00
-1.2750e+00 1.7263e+00
-1.3120e+00 1.7078e+00
-1.3486e+00 1.6819e+00
-1.3846e+00 1.6488e+00
-1.4198e+00 1.6087e+00
-1.4540e+00 1.5617e+00
-1.4872e+00 1.5080e+00
-1.5191e+00 1.4478e+00
-1.5496e+00 1.3814e+00
-1.5787e+00 1.3091e+00
-1.6061e+00 1.2312e+00
-1.6318e+00 1.1481e+00
-1.6556e+00 1.0600e+00
-1.6775e+00 9.6738e-01
-1.6974e+00 8.7062e-01
-1.7151e+00 7.7013e-01
-1.7306e+00 6.6635e-01
-1.7438e+00 5.5970e-01
-1.7547e+00 4.5067e-01
-1.7633e+00 3.3970e-01
-1.7694e+00 2.2728e-01
-1.7731e+00 1.1388e-01
-1.7743e+00 1.4179e-16
-1.7731e+00 -1.1388e-01
-1.7694e+00 -2.2728e-01
-1.7633e+00 -3.3970e-01
-1.7547e+00 -4.5067e-01
-1.7438e+00 -5.5970e-01
-1.7306e+00 -6.6635e-01
-1.7151e+00 -7.7013e-01
-1.6974e+00 -8.7062e-01
-1.6775e+00 -9.6738e-01
-1.6556e+00 -1.0600e+00
-1.6318e+00 -1.1481e+00
-1.6061e+00 -1.2312e+00
-1.5787e+00 -1.3091e+00
-1.5496e+00 -1.3814e+00
-1.5191e+00 -1.4478e+00
-1.4872e+00 -1.5080e+00
-1.4540e+00 -1.5617e+00
-1.4198e+00 -1.6087e+00
-1.3846e+00 -1.6488e+00
-1.3486e+00 -1.6819e+00
-1.3120e+00 -1.7078e+00
-1.2750e+00 -1.7263e+00
-1.2376e+00 -1.7375e+00
-1.2000e+00 -1.7412e+00
-1.1624e+00 -1.7375e+00
-1.1250e+00 -1.7263e+00
-1.0880e+00 -1.7078e+00
-1.0514e+00 -1.6819e+00
-1.0154e+00 -1.6488e+00
-9.8022e-01 -1.6087e+00
-9.4599e-01 -1.5617e+00
-9.1285e-01 -1.5080e+00
-8.8093e-01 -1.4478e+00
-8.5039e-01 -1.3814e+00
-8.2134e-01 -1.3091e+00
-7.9391e-01 -1.2312e+00
-7.6822e-01 -1.1481e+00
-7.4438e-01 -1.0600e+00
-7.2249e-01 -9.6738e-01
-7.0264e-01 -8.7062e-01
-6.8492e-01 -7.7013e-01
-6.6941e-01 -6.6635e-01
-6.5618e-01 -5.5970e-01
-6.4527e-01 -4.5067e-01
-6.3673e-01 -3.3970e-01
-6.3061e-01 -2.2728e-01
-6.2693e-01 -1.1388e-01
-6.2570e-01 -3.5503e-16
-6.2693e-01 1.1388e-01
-6.3061e-01 2.2728e-01
-6.3673e-01 3.3970e-01
-6.4527e-01 4.5067e-01
-6.5618e-01 5.5970e-01
-6.6941e-01 6.6635e-01
-6.8492e-01 7.7013e-01
-7.0264e-01 8.7062e-01
-7.2249e-01 9.6738e-01
-7.4438e-01 1.0600e+00
-7.6822e-01 1.1481e+00
-7.9391e-01 1.2312e+00
-8.2134e-01 1.3091e+00
-8.5039e-01 1.3814e+00
-8.8093e-01 1.4478e+00
-9.1285e-01 1.5080e+00
-9.4599e-01 1.5617e+00
-9.8022e-01 1.6087e+00
-1.0154e+00 1.6488e+00
-1.0514e+00 1.6819e+00
-1.0880e+00 1.7078e+00
-1.1250e+00 1.7263e+00
-1.1624e+00 1.7375e+00
-1.2000e+00 1.7412e+00
};

\addplot [mark=none,red,line width=1.5] table{
1.2000e+00 1.7412e+00
1.1624e+00 1.7375e+00
1.1250e+00 1.7263e+00
1.0880e+00 1.7078e+00
1.0514e+00 1.6819e+00
1.0154e+00 1.6488e+00
9.8022e-01 1.6087e+00
9.4599e-01 1.5617e+00
9.1285e-01 1.5080e+00
8.8093e-01 1.4478e+00
8.5039e-01 1.3814e+00
8.2134e-01 1.3091e+00
7.9391e-01 1.2312e+00
7.6822e-01 1.1481e+00
7.4438e-01 1.0600e+00
7.2249e-01 9.6738e-01
7.0264e-01 8.7062e-01
6.8492e-01 7.7013e-01
6.6941e-01 6.6635e-01
6.5618e-01 5.5970e-01
6.4527e-01 4.5067e-01
6.3673e-01 3.3970e-01
6.3061e-01 2.2728e-01
6.2693e-01 1.1388e-01
6.2570e-01 1.4179e-16
6.2693e-01 -1.1388e-01
6.3061e-01 -2.2728e-01
6.3673e-01 -3.3970e-01
6.4527e-01 -4.5067e-01
6.5618e-01 -5.5970e-01
6.6941e-01 -6.6635e-01
6.8492e-01 -7.7013e-01
7.0264e-01 -8.7062e-01
7.2249e-01 -9.6738e-01
7.4438e-01 -1.0600e+00
7.6822e-01 -1.1481e+00
7.9391e-01 -1.2312e+00
8.2134e-01 -1.3091e+00
8.5039e-01 -1.3814e+00
8.8093e-01 -1.4478e+00
9.1285e-01 -1.5080e+00
9.4599e-01 -1.5617e+00
9.8022e-01 -1.6087e+00
1.0154e+00 -1.6488e+00
1.0514e+00 -1.6819e+00
1.0880e+00 -1.7078e+00
1.1250e+00 -1.7263e+00
1.1624e+00 -1.7375e+00
1.2000e+00 -1.7412e+00
1.2376e+00 -1.7375e+00
1.2750e+00 -1.7263e+00
1.3120e+00 -1.7078e+00
1.3486e+00 -1.6819e+00
1.3846e+00 -1.6488e+00
1.4198e+00 -1.6087e+00
1.4540e+00 -1.5617e+00
1.4872e+00 -1.5080e+00
1.5191e+00 -1.4478e+00
1.5496e+00 -1.3814e+00
1.5787e+00 -1.3091e+00
1.6061e+00 -1.2312e+00
1.6318e+00 -1.1481e+00
1.6556e+00 -1.0600e+00
1.6775e+00 -9.6738e-01
1.6974e+00 -8.7062e-01
1.7151e+00 -7.7013e-01
1.7306e+00 -6.6635e-01
1.7438e+00 -5.5970e-01
1.7547e+00 -4.5067e-01
1.7633e+00 -3.3970e-01
1.7694e+00 -2.2728e-01
1.7731e+00 -1.1388e-01
1.7743e+00 -3.5503e-16
1.7731e+00 1.1388e-01
1.7694e+00 2.2728e-01
1.7633e+00 3.3970e-01
1.7547e+00 4.5067e-01
1.7438e+00 5.5970e-01
1.7306e+00 6.6635e-01
1.7151e+00 7.7013e-01
1.6974e+00 8.7062e-01
1.6775e+00 9.6738e-01
1.6556e+00 1.0600e+00
1.6318e+00 1.1481e+00
1.6061e+00 1.2312e+00
1.5787e+00 1.3091e+00
1.5496e+00 1.3814e+00
1.5191e+00 1.4478e+00
1.4872e+00 1.5080e+00
1.4540e+00 1.5617e+00
1.4198e+00 1.6087e+00
1.3846e+00 1.6488e+00
1.3486e+00 1.6819e+00
1.3120e+00 1.7078e+00
1.2750e+00 1.7263e+00
1.2376e+00 1.7375e+00
1.2000e+00 1.7412e+00
};

\end{axis}

\end{tikzpicture} &
  \begin{tikzpicture}[scale=0.40]

\begin{axis}[
  xmin = -1.8,
  xmax = 1.8,
  ymin = -1.8,
  ymax = 1.8,
  axis equal = true,
  hide axis,
  title = {\Huge$t=1$}
  ]

\addplot [mark=none,red,line width=1.5] table{
-8.4626e-01 1.6919e+00
-8.8341e-01 1.6850e+00
-9.2082e-01 1.6737e+00
-9.5905e-01 1.6576e+00
-9.9832e-01 1.6360e+00
-1.0385e+00 1.6081e+00
-1.0791e+00 1.5736e+00
-1.1198e+00 1.5319e+00
-1.1598e+00 1.4832e+00
-1.1989e+00 1.4274e+00
-1.2365e+00 1.3648e+00
-1.2726e+00 1.2957e+00
-1.3069e+00 1.2206e+00
-1.3392e+00 1.1398e+00
-1.3694e+00 1.0537e+00
-1.3973e+00 9.6268e-01
-1.4229e+00 8.6727e-01
-1.4459e+00 7.6786e-01
-1.4661e+00 6.6489e-01
-1.4835e+00 5.5885e-01
-1.4980e+00 4.5023e-01
-1.5093e+00 3.3951e-01
-1.5174e+00 2.2722e-01
-1.5224e+00 1.1388e-01
-1.5240e+00 -2.1459e-15
-1.5224e+00 -1.1388e-01
-1.5174e+00 -2.2722e-01
-1.5093e+00 -3.3951e-01
-1.4980e+00 -4.5023e-01
-1.4835e+00 -5.5885e-01
-1.4661e+00 -6.6489e-01
-1.4459e+00 -7.6786e-01
-1.4229e+00 -8.6727e-01
-1.3973e+00 -9.6268e-01
-1.3694e+00 -1.0537e+00
-1.3392e+00 -1.1398e+00
-1.3069e+00 -1.2206e+00
-1.2726e+00 -1.2957e+00
-1.2365e+00 -1.3648e+00
-1.1989e+00 -1.4274e+00
-1.1598e+00 -1.4832e+00
-1.1198e+00 -1.5319e+00
-1.0791e+00 -1.5736e+00
-1.0385e+00 -1.6081e+00
-9.9832e-01 -1.6360e+00
-9.5905e-01 -1.6576e+00
-9.2082e-01 -1.6737e+00
-8.8341e-01 -1.6850e+00
-8.4626e-01 -1.6919e+00
-8.0857e-01 -1.6946e+00
-7.6955e-01 -1.6928e+00
-7.2869e-01 -1.6857e+00
-6.8600e-01 -1.6720e+00
-6.4213e-01 -1.6506e+00
-5.9836e-01 -1.6201e+00
-5.5638e-01 -1.5799e+00
-5.1798e-01 -1.5299e+00
-4.8468e-01 -1.4705e+00
-4.5747e-01 -1.4027e+00
-4.3668e-01 -1.3277e+00
-4.2205e-01 -1.2464e+00
-4.1290e-01 -1.1599e+00
-4.0831e-01 -1.0687e+00
-4.0730e-01 -9.7357e-01
-4.0889e-01 -8.7481e-01
-4.1221e-01 -7.7282e-01
-4.1651e-01 -6.6797e-01
-4.2115e-01 -5.6060e-01
-4.2562e-01 -4.5111e-01
-4.2949e-01 -3.3988e-01
-4.3248e-01 -2.2733e-01
-4.3435e-01 -1.1389e-01
-4.3499e-01 -6.3033e-15
-4.3435e-01 1.1389e-01
-4.3248e-01 2.2733e-01
-4.2949e-01 3.3988e-01
-4.2562e-01 4.5111e-01
-4.2115e-01 5.6060e-01
-4.1651e-01 6.6797e-01
-4.1221e-01 7.7282e-01
-4.0889e-01 8.7481e-01
-4.0730e-01 9.7357e-01
-4.0831e-01 1.0687e+00
-4.1290e-01 1.1599e+00
-4.2205e-01 1.2464e+00
-4.3668e-01 1.3277e+00
-4.5747e-01 1.4027e+00
-4.8468e-01 1.4705e+00
-5.1798e-01 1.5299e+00
-5.5638e-01 1.5799e+00
-5.9836e-01 1.6201e+00
-6.4213e-01 1.6506e+00
-6.8600e-01 1.6720e+00
-7.2869e-01 1.6857e+00
-7.6955e-01 1.6928e+00
-8.0857e-01 1.6946e+00
-8.4626e-01 1.6919e+00
};

\addplot [mark=none,red,line width=1.5] table{
8.4626e-01 1.6919e+00
8.0857e-01 1.6946e+00
7.6955e-01 1.6928e+00
7.2869e-01 1.6857e+00
6.8600e-01 1.6720e+00
6.4213e-01 1.6506e+00
5.9836e-01 1.6201e+00
5.5638e-01 1.5799e+00
5.1798e-01 1.5299e+00
4.8468e-01 1.4705e+00
4.5747e-01 1.4027e+00
4.3668e-01 1.3277e+00
4.2205e-01 1.2464e+00
4.1290e-01 1.1599e+00
4.0831e-01 1.0687e+00
4.0730e-01 9.7357e-01
4.0889e-01 8.7481e-01
4.1221e-01 7.7282e-01
4.1651e-01 6.6797e-01
4.2115e-01 5.6060e-01
4.2562e-01 4.5111e-01
4.2949e-01 3.3988e-01
4.3248e-01 2.2733e-01
4.3435e-01 1.1389e-01
4.3499e-01 -5.9187e-15
4.3435e-01 -1.1389e-01
4.3248e-01 -2.2733e-01
4.2949e-01 -3.3988e-01
4.2562e-01 -4.5111e-01
4.2115e-01 -5.6060e-01
4.1651e-01 -6.6797e-01
4.1221e-01 -7.7282e-01
4.0889e-01 -8.7481e-01
4.0730e-01 -9.7357e-01
4.0831e-01 -1.0687e+00
4.1290e-01 -1.1599e+00
4.2205e-01 -1.2464e+00
4.3668e-01 -1.3277e+00
4.5747e-01 -1.4027e+00
4.8468e-01 -1.4705e+00
5.1798e-01 -1.5299e+00
5.5638e-01 -1.5799e+00
5.9836e-01 -1.6201e+00
6.4213e-01 -1.6506e+00
6.8600e-01 -1.6720e+00
7.2869e-01 -1.6857e+00
7.6955e-01 -1.6928e+00
8.0857e-01 -1.6946e+00
8.4626e-01 -1.6919e+00
8.8341e-01 -1.6850e+00
9.2082e-01 -1.6737e+00
9.5905e-01 -1.6576e+00
9.9832e-01 -1.6360e+00
1.0385e+00 -1.6081e+00
1.0791e+00 -1.5736e+00
1.1198e+00 -1.5319e+00
1.1598e+00 -1.4832e+00
1.1989e+00 -1.4274e+00
1.2365e+00 -1.3648e+00
1.2726e+00 -1.2957e+00
1.3069e+00 -1.2206e+00
1.3392e+00 -1.1398e+00
1.3694e+00 -1.0537e+00
1.3973e+00 -9.6268e-01
1.4229e+00 -8.6727e-01
1.4459e+00 -7.6786e-01
1.4661e+00 -6.6489e-01
1.4835e+00 -5.5885e-01
1.4980e+00 -4.5023e-01
1.5093e+00 -3.3951e-01
1.5174e+00 -2.2722e-01
1.5224e+00 -1.1388e-01
1.5240e+00 -1.2890e-15
1.5224e+00 1.1388e-01
1.5174e+00 2.2722e-01
1.5093e+00 3.3951e-01
1.4980e+00 4.5023e-01
1.4835e+00 5.5885e-01
1.4661e+00 6.6489e-01
1.4459e+00 7.6786e-01
1.4229e+00 8.6727e-01
1.3973e+00 9.6268e-01
1.3694e+00 1.0537e+00
1.3392e+00 1.1398e+00
1.3069e+00 1.2206e+00
1.2726e+00 1.2957e+00
1.2365e+00 1.3648e+00
1.1989e+00 1.4274e+00
1.1598e+00 1.4832e+00
1.1198e+00 1.5319e+00
1.0791e+00 1.5736e+00
1.0385e+00 1.6081e+00
9.9832e-01 1.6360e+00
9.5905e-01 1.6576e+00
9.2082e-01 1.6737e+00
8.8341e-01 1.6850e+00
8.4626e-01 1.6919e+00
};

\end{axis}

\end{tikzpicture} &
  \begin{tikzpicture}[scale=0.40]

\begin{axis}[
  xmin = -1.8,
  xmax = 1.8,
  ymin = -1.8,
  ymax = 1.8,
  axis equal = true,
  hide axis,
  title = {\Huge$t=5$}
  ]

\addplot [mark=none,red,line width=1.5] table{
-6.0188e-01 1.6622e+00
-6.3899e-01 1.6550e+00
-6.7650e-01 1.6441e+00
-7.1509e-01 1.6288e+00
-7.5512e-01 1.6086e+00
-7.9656e-01 1.5827e+00
-8.3911e-01 1.5505e+00
-8.8228e-01 1.5115e+00
-9.2551e-01 1.4655e+00
-9.6825e-01 1.4125e+00
-1.0100e+00 1.3526e+00
-1.0503e+00 1.2859e+00
-1.0889e+00 1.2129e+00
-1.1254e+00 1.1339e+00
-1.1596e+00 1.0493e+00
-1.1913e+00 9.5956e-01
-1.2202e+00 8.6511e-01
-1.2462e+00 7.6643e-01
-1.2690e+00 6.6400e-01
-1.2886e+00 5.5835e-01
-1.3048e+00 4.4997e-01
-1.3175e+00 3.3941e-01
-1.3266e+00 2.2719e-01
-1.3321e+00 1.1387e-01
-1.3339e+00 -5.5197e-13
-1.3321e+00 -1.1387e-01
-1.3266e+00 -2.2719e-01
-1.3175e+00 -3.3941e-01
-1.3048e+00 -4.4997e-01
-1.2886e+00 -5.5835e-01
-1.2690e+00 -6.6400e-01
-1.2462e+00 -7.6643e-01
-1.2202e+00 -8.6511e-01
-1.1913e+00 -9.5956e-01
-1.1596e+00 -1.0493e+00
-1.1254e+00 -1.1339e+00
-1.0889e+00 -1.2129e+00
-1.0503e+00 -1.2859e+00
-1.0100e+00 -1.3526e+00
-9.6825e-01 -1.4125e+00
-9.2551e-01 -1.4655e+00
-8.8228e-01 -1.5115e+00
-8.3911e-01 -1.5505e+00
-7.9656e-01 -1.5827e+00
-7.5512e-01 -1.6086e+00
-7.1509e-01 -1.6288e+00
-6.7650e-01 -1.6441e+00
-6.3899e-01 -1.6550e+00
-6.0188e-01 -1.6622e+00
-5.6425e-01 -1.6657e+00
-5.2518e-01 -1.6653e+00
-4.8399e-01 -1.6603e+00
-4.4047e-01 -1.6496e+00
-3.9502e-01 -1.6317e+00
-3.4873e-01 -1.6052e+00
-3.0332e-01 -1.5689e+00
-2.6098e-01 -1.5222e+00
-2.2404e-01 -1.4651e+00
-1.9453e-01 -1.3983e+00
-1.7381e-01 -1.3232e+00
-1.6222e-01 -1.2415e+00
-1.5912e-01 -1.1546e+00
-1.6306e-01 -1.0634e+00
-1.7224e-01 -9.6869e-01
-1.8482e-01 -8.7072e-01
-1.9916e-01 -7.6969e-01
-2.1388e-01 -6.6578e-01
-2.2791e-01 -5.5924e-01
-2.4039e-01 -4.5037e-01
-2.5070e-01 -3.3956e-01
-2.5838e-01 -2.2723e-01
-2.6311e-01 -1.1388e-01
-2.6471e-01 -4.9031e-13
-2.6311e-01 1.1388e-01
-2.5838e-01 2.2723e-01
-2.5070e-01 3.3956e-01
-2.4039e-01 4.5037e-01
-2.2791e-01 5.5924e-01
-2.1388e-01 6.6578e-01
-1.9916e-01 7.6969e-01
-1.8482e-01 8.7072e-01
-1.7224e-01 9.6869e-01
-1.6306e-01 1.0634e+00
-1.5912e-01 1.1546e+00
-1.6222e-01 1.2415e+00
-1.7381e-01 1.3232e+00
-1.9453e-01 1.3983e+00
-2.2404e-01 1.4651e+00
-2.6098e-01 1.5222e+00
-3.0332e-01 1.5689e+00
-3.4873e-01 1.6052e+00
-3.9502e-01 1.6317e+00
-4.4047e-01 1.6496e+00
-4.8399e-01 1.6603e+00
-5.2518e-01 1.6653e+00
-5.6425e-01 1.6657e+00
-6.0188e-01 1.6622e+00
};

\addplot [mark=none,red,line width=1.5] table{
6.0188e-01 1.6622e+00
5.6425e-01 1.6657e+00
5.2518e-01 1.6653e+00
4.8399e-01 1.6603e+00
4.4047e-01 1.6496e+00
3.9502e-01 1.6317e+00
3.4873e-01 1.6052e+00
3.0332e-01 1.5689e+00
2.6098e-01 1.5222e+00
2.2404e-01 1.4651e+00
1.9453e-01 1.3983e+00
1.7381e-01 1.3232e+00
1.6222e-01 1.2415e+00
1.5912e-01 1.1546e+00
1.6306e-01 1.0634e+00
1.7224e-01 9.6869e-01
1.8482e-01 8.7072e-01
1.9916e-01 7.6969e-01
2.1388e-01 6.6578e-01
2.2791e-01 5.5924e-01
2.4039e-01 4.5037e-01
2.5070e-01 3.3956e-01
2.5838e-01 2.2723e-01
2.6311e-01 1.1388e-01
2.6471e-01 -3.0230e-13
2.6311e-01 -1.1388e-01
2.5838e-01 -2.2723e-01
2.5070e-01 -3.3956e-01
2.4039e-01 -4.5037e-01
2.2791e-01 -5.5924e-01
2.1388e-01 -6.6578e-01
1.9916e-01 -7.6969e-01
1.8482e-01 -8.7072e-01
1.7224e-01 -9.6869e-01
1.6306e-01 -1.0634e+00
1.5912e-01 -1.1546e+00
1.6222e-01 -1.2415e+00
1.7381e-01 -1.3232e+00
1.9453e-01 -1.3983e+00
2.2404e-01 -1.4651e+00
2.6098e-01 -1.5222e+00
3.0332e-01 -1.5689e+00
3.4873e-01 -1.6052e+00
3.9502e-01 -1.6317e+00
4.4047e-01 -1.6496e+00
4.8399e-01 -1.6603e+00
5.2518e-01 -1.6653e+00
5.6425e-01 -1.6657e+00
6.0188e-01 -1.6622e+00
6.3899e-01 -1.6550e+00
6.7650e-01 -1.6441e+00
7.1509e-01 -1.6288e+00
7.5512e-01 -1.6086e+00
7.9656e-01 -1.5827e+00
8.3911e-01 -1.5505e+00
8.8228e-01 -1.5115e+00
9.2551e-01 -1.4655e+00
9.6825e-01 -1.4125e+00
1.0100e+00 -1.3526e+00
1.0503e+00 -1.2859e+00
1.0889e+00 -1.2129e+00
1.1254e+00 -1.1339e+00
1.1596e+00 -1.0493e+00
1.1913e+00 -9.5956e-01
1.2202e+00 -8.6511e-01
1.2462e+00 -7.6643e-01
1.2690e+00 -6.6400e-01
1.2886e+00 -5.5835e-01
1.3048e+00 -4.4997e-01
1.3175e+00 -3.3941e-01
1.3266e+00 -2.2719e-01
1.3321e+00 -1.1387e-01
1.3339e+00 -1.9961e-13
1.3321e+00 1.1387e-01
1.3266e+00 2.2719e-01
1.3175e+00 3.3941e-01
1.3048e+00 4.4997e-01
1.2886e+00 5.5835e-01
1.2690e+00 6.6400e-01
1.2462e+00 7.6643e-01
1.2202e+00 8.6511e-01
1.1913e+00 9.5956e-01
1.1596e+00 1.0493e+00
1.1254e+00 1.1339e+00
1.0889e+00 1.2129e+00
1.0503e+00 1.2859e+00
1.0100e+00 1.3526e+00
9.6825e-01 1.4125e+00
9.2551e-01 1.4655e+00
8.8228e-01 1.5115e+00
8.3911e-01 1.5505e+00
7.9656e-01 1.5827e+00
7.5512e-01 1.6086e+00
7.1509e-01 1.6288e+00
6.7650e-01 1.6441e+00
6.3899e-01 1.6550e+00
6.0188e-01 1.6622e+00
};

\end{axis}

\end{tikzpicture} &
  \begin{tikzpicture}[scale=0.40]

\begin{axis}[
  xmin = -1.8,
  xmax = 1.8,
  ymin = -1.8,
  ymax = 1.8,
  axis equal = true,
  hide axis,
  title = {\Huge$t=12$}
  ]

\addplot [mark=none,red,line width=1.5] table{
-5.3331e-01 1.6539e+00
-5.7044e-01 1.6469e+00
-6.0802e-01 1.6362e+00
-6.4676e-01 1.6213e+00
-6.8703e-01 1.6016e+00
-7.2883e-01 1.5763e+00
-7.7188e-01 1.5447e+00
-8.1568e-01 1.5064e+00
-8.5968e-01 1.4612e+00
-9.0328e-01 1.4089e+00
-9.4596e-01 1.3496e+00
-9.8727e-01 1.2835e+00
-1.0268e+00 1.2110e+00
-1.0643e+00 1.1325e+00
-1.0994e+00 1.0482e+00
-1.1318e+00 9.5879e-01
-1.1615e+00 8.6458e-01
-1.1881e+00 7.6607e-01
-1.2115e+00 6.6379e-01
-1.2316e+00 5.5822e-01
-1.2482e+00 4.4991e-01
-1.2612e+00 3.3938e-01
-1.2706e+00 2.2718e-01
-1.2762e+00 1.1387e-01
-1.2781e+00 -5.1688e-10
-1.2762e+00 -1.1387e-01
-1.2706e+00 -2.2718e-01
-1.2612e+00 -3.3938e-01
-1.2482e+00 -4.4991e-01
-1.2316e+00 -5.5822e-01
-1.2115e+00 -6.6379e-01
-1.1881e+00 -7.6607e-01
-1.1615e+00 -8.6458e-01
-1.1318e+00 -9.5879e-01
-1.0994e+00 -1.0482e+00
-1.0643e+00 -1.1325e+00
-1.0268e+00 -1.2110e+00
-9.8727e-01 -1.2835e+00
-9.4596e-01 -1.3496e+00
-9.0328e-01 -1.4089e+00
-8.5968e-01 -1.4612e+00
-8.1568e-01 -1.5064e+00
-7.7188e-01 -1.5447e+00
-7.2883e-01 -1.5763e+00
-6.8703e-01 -1.6016e+00
-6.4676e-01 -1.6213e+00
-6.0802e-01 -1.6362e+00
-5.7044e-01 -1.6469e+00
-5.3331e-01 -1.6539e+00
-4.9568e-01 -1.6574e+00
-4.5661e-01 -1.6572e+00
-4.1538e-01 -1.6526e+00
-3.7172e-01 -1.6424e+00
-3.2596e-01 -1.6253e+00
-2.7912e-01 -1.5998e+00
-2.3286e-01 -1.5646e+00
-1.8941e-01 -1.5189e+00
-1.5128e-01 -1.4625e+00
-1.2086e-01 -1.3962e+00
-9.9915e-02 -1.3212e+00
-8.9108e-02 -1.2394e+00
-8.7789e-02 -1.1524e+00
-9.4241e-02 -1.0614e+00
-1.0625e-01 -9.6700e-01
-1.2161e-01 -8.6942e-01
-1.3847e-01 -7.6878e-01
-1.5533e-01 -6.6521e-01
-1.7110e-01 -5.5891e-01
-1.8494e-01 -4.5020e-01
-1.9626e-01 -3.3949e-01
-2.0463e-01 -2.2721e-01
-2.0976e-01 -1.1387e-01
-2.1149e-01 -4.5199e-10
-2.0976e-01 1.1387e-01
-2.0463e-01 2.2721e-01
-1.9626e-01 3.3949e-01
-1.8494e-01 4.5020e-01
-1.7110e-01 5.5891e-01
-1.5533e-01 6.6521e-01
-1.3847e-01 7.6878e-01
-1.2161e-01 8.6942e-01
-1.0625e-01 9.6700e-01
-9.4241e-02 1.0614e+00
-8.7789e-02 1.1524e+00
-8.9108e-02 1.2394e+00
-9.9915e-02 1.3212e+00
-1.2086e-01 1.3962e+00
-1.5128e-01 1.4625e+00
-1.8941e-01 1.5189e+00
-2.3286e-01 1.5646e+00
-2.7912e-01 1.5998e+00
-3.2596e-01 1.6253e+00
-3.7172e-01 1.6424e+00
-4.1538e-01 1.6526e+00
-4.5661e-01 1.6572e+00
-4.9568e-01 1.6574e+00
-5.3331e-01 1.6539e+00
};

\addplot [mark=none,red,line width=1.5] table{
5.3331e-01 1.6539e+00
4.9568e-01 1.6574e+00
4.5661e-01 1.6572e+00
4.1538e-01 1.6526e+00
3.7172e-01 1.6424e+00
3.2596e-01 1.6253e+00
2.7912e-01 1.5998e+00
2.3286e-01 1.5646e+00
1.8941e-01 1.5189e+00
1.5128e-01 1.4625e+00
1.2086e-01 1.3962e+00
9.9915e-02 1.3212e+00
8.9108e-02 1.2394e+00
8.7789e-02 1.1524e+00
9.4241e-02 1.0614e+00
1.0625e-01 9.6700e-01
1.2161e-01 8.6942e-01
1.3847e-01 7.6878e-01
1.5533e-01 6.6521e-01
1.7110e-01 5.5891e-01
1.8494e-01 4.5020e-01
1.9626e-01 3.3949e-01
2.0463e-01 2.2721e-01
2.0976e-01 1.1387e-01
2.1149e-01 -3.7564e-10
2.0976e-01 -1.1387e-01
2.0463e-01 -2.2721e-01
1.9626e-01 -3.3949e-01
1.8494e-01 -4.5020e-01
1.7110e-01 -5.5891e-01
1.5533e-01 -6.6521e-01
1.3847e-01 -7.6878e-01
1.2161e-01 -8.6942e-01
1.0625e-01 -9.6700e-01
9.4241e-02 -1.0614e+00
8.7789e-02 -1.1524e+00
8.9108e-02 -1.2394e+00
9.9915e-02 -1.3212e+00
1.2086e-01 -1.3962e+00
1.5128e-01 -1.4625e+00
1.8941e-01 -1.5189e+00
2.3286e-01 -1.5646e+00
2.7912e-01 -1.5998e+00
3.2596e-01 -1.6253e+00
3.7172e-01 -1.6424e+00
4.1538e-01 -1.6526e+00
4.5661e-01 -1.6572e+00
4.9568e-01 -1.6574e+00
5.3331e-01 -1.6539e+00
5.7044e-01 -1.6469e+00
6.0802e-01 -1.6362e+00
6.4676e-01 -1.6213e+00
6.8703e-01 -1.6016e+00
7.2883e-01 -1.5763e+00
7.7188e-01 -1.5447e+00
8.1568e-01 -1.5064e+00
8.5968e-01 -1.4612e+00
9.0328e-01 -1.4089e+00
9.4596e-01 -1.3496e+00
9.8727e-01 -1.2835e+00
1.0268e+00 -1.2110e+00
1.0643e+00 -1.1325e+00
1.0994e+00 -1.0482e+00
1.1318e+00 -9.5879e-01
1.1615e+00 -8.6458e-01
1.1881e+00 -7.6607e-01
1.2115e+00 -6.6379e-01
1.2316e+00 -5.5822e-01
1.2482e+00 -4.4991e-01
1.2612e+00 -3.3938e-01
1.2706e+00 -2.2718e-01
1.2762e+00 -1.1387e-01
1.2781e+00 -3.1094e-10
1.2762e+00 1.1387e-01
1.2706e+00 2.2718e-01
1.2612e+00 3.3938e-01
1.2482e+00 4.4991e-01
1.2316e+00 5.5822e-01
1.2115e+00 6.6379e-01
1.1881e+00 7.6607e-01
1.1615e+00 8.6458e-01
1.1318e+00 9.5879e-01
1.0994e+00 1.0482e+00
1.0643e+00 1.1325e+00
1.0268e+00 1.2110e+00
9.8727e-01 1.2835e+00
9.4596e-01 1.3496e+00
9.0328e-01 1.4089e+00
8.5968e-01 1.4612e+00
8.1568e-01 1.5064e+00
7.7188e-01 1.5447e+00
7.2883e-01 1.5763e+00
6.8703e-01 1.6016e+00
6.4676e-01 1.6213e+00
6.0802e-01 1.6362e+00
5.3331e-01 1.6539e+00
5.7044e-01 1.6469e+00
};

\end{axis}

\end{tikzpicture} &
  \begin{tikzpicture}[scale=0.40]

\begin{axis}[
  xmin = -1.8,
  xmax = 1.8,
  ymin = -1.8,
  ymax = 1.8,
  axis equal = true,
  hide axis,
  title = {\Huge$t=24$}
  ]

\addplot [mark=none,red,line width=1.5] table{
-5.0194e-01 1.6512e+00
-5.3909e-01 1.6442e+00
-5.7669e-01 1.6336e+00
-6.1548e-01 1.6188e+00
-6.5582e-01 1.5993e+00
-6.9775e-01 1.5741e+00
-7.4095e-01 1.5428e+00
-7.8496e-01 1.5048e+00
-8.2920e-01 1.4598e+00
-8.7307e-01 1.4077e+00
-9.1604e-01 1.3486e+00
-9.5764e-01 1.2827e+00
-9.9747e-01 1.2104e+00
-1.0352e+00 1.1320e+00
-1.0706e+00 1.0479e+00
-1.1033e+00 9.5849e-01
-1.1332e+00 8.6435e-01
-1.1600e+00 7.6590e-01
-1.1836e+00 6.6365e-01
-1.2038e+00 5.5811e-01
-1.2205e+00 4.4982e-01
-1.2336e+00 3.3930e-01
-1.2430e+00 2.2711e-01
-1.2487e+00 1.1380e-01
-1.2506e+00 -6.9370e-05
-1.2487e+00 -1.1394e-01
-1.2430e+00 -2.2725e-01
-1.2336e+00 -3.3944e-01
-1.2205e+00 -4.4996e-01
-1.2038e+00 -5.5825e-01
-1.1836e+00 -6.6379e-01
-1.1600e+00 -7.6604e-01
-1.1332e+00 -8.6449e-01
-1.1033e+00 -9.5863e-01
-1.0706e+00 -1.0480e+00
-1.0352e+00 -1.1321e+00
-9.9747e-01 -1.2105e+00
-9.5764e-01 -1.2829e+00
-9.1604e-01 -1.3487e+00
-8.7307e-01 -1.4078e+00
-8.2919e-01 -1.4599e+00
-7.8496e-01 -1.5049e+00
-7.4095e-01 -1.5429e+00
-6.9775e-01 -1.5743e+00
-6.5582e-01 -1.5994e+00
-6.1548e-01 -1.6190e+00
-5.7669e-01 -1.6337e+00
-5.3908e-01 -1.6443e+00
-5.0194e-01 -1.6513e+00
-4.6431e-01 -1.6548e+00
-4.2524e-01 -1.6547e+00
-3.8400e-01 -1.6502e+00
-3.4029e-01 -1.6402e+00
-2.9443e-01 -1.6233e+00
-2.4739e-01 -1.5982e+00
-2.0083e-01 -1.5634e+00
-1.5697e-01 -1.5181e+00
-1.1841e-01 -1.4620e+00
-8.7681e-02 -1.3958e+00
-6.6776e-02 -1.3209e+00
-5.6484e-02 -1.2390e+00
-5.6128e-02 -1.1520e+00
-6.3728e-02 -1.0611e+00
-7.6812e-02 -9.6681e-01
-9.3003e-02 -8.6936e-01
-1.1040e-01 -7.6882e-01
-1.2755e-01 -6.6528e-01
-1.4341e-01 -5.5900e-01
-1.5722e-01 -4.5029e-01
-1.6845e-01 -3.3957e-01
-1.7671e-01 -2.2728e-01
-1.8176e-01 -1.1394e-01
-1.8346e-01 -6.7945e-05
-1.8176e-01 1.1381e-01
-1.7671e-01 2.2715e-01
-1.6845e-01 3.3943e-01
-1.5722e-01 4.5016e-01
-1.4341e-01 5.5887e-01
-1.2755e-01 6.6515e-01
-1.1040e-01 7.6868e-01
-9.3004e-02 8.6923e-01
-7.6814e-02 9.6667e-01
-6.3730e-02 1.0609e+00
-5.6129e-02 1.1519e+00
-5.6486e-02 1.2389e+00
-6.6778e-02 1.3207e+00
-8.7683e-02 1.3957e+00
-1.1841e-01 1.4619e+00
-1.5697e-01 1.5180e+00
-2.0083e-01 1.5633e+00
-2.4739e-01 1.5981e+00
-2.9444e-01 1.6232e+00
-3.4029e-01 1.6400e+00
-3.8400e-01 1.6500e+00
-4.2524e-01 1.6546e+00
-4.6432e-01 1.6547e+00
-5.0194e-01 1.6512e+00
};

\addplot [mark=none,red,line width=1.5] table{
5.0194e-01 1.6512e+00
4.6431e-01 1.6547e+00
4.2524e-01 1.6546e+00
3.8400e-01 1.6500e+00
3.4029e-01 1.6400e+00
2.9443e-01 1.6232e+00
2.4739e-01 1.5981e+00
2.0083e-01 1.5633e+00
1.5697e-01 1.5180e+00
1.1841e-01 1.4619e+00
8.7681e-02 1.3957e+00
6.6776e-02 1.3207e+00
5.6484e-02 1.2389e+00
5.6128e-02 1.1519e+00
6.3728e-02 1.0609e+00
7.6812e-02 9.6667e-01
9.3003e-02 8.6923e-01
1.1040e-01 7.6868e-01
1.2755e-01 6.6515e-01
1.4341e-01 5.5887e-01
1.5722e-01 4.5016e-01
1.6845e-01 3.3943e-01
1.7671e-01 2.2715e-01
1.8176e-01 1.1381e-01
1.8346e-01 -6.6706e-05
1.8176e-01 -1.1394e-01
1.7671e-01 -2.2728e-01
1.6845e-01 -3.3957e-01
1.5722e-01 -4.5029e-01
1.4341e-01 -5.5900e-01
1.2755e-01 -6.6528e-01
1.1040e-01 -7.6881e-01
9.3004e-02 -8.6936e-01
7.6814e-02 -9.6680e-01
6.3730e-02 -1.0611e+00
5.6129e-02 -1.1520e+00
5.6486e-02 -1.2390e+00
6.6778e-02 -1.3209e+00
8.7683e-02 -1.3958e+00
1.1841e-01 -1.4620e+00
1.5697e-01 -1.5181e+00
2.0083e-01 -1.5634e+00
2.4739e-01 -1.5982e+00
2.9444e-01 -1.6233e+00
3.4029e-01 -1.6402e+00
3.8400e-01 -1.6502e+00
4.2524e-01 -1.6547e+00
4.6432e-01 -1.6548e+00
5.0194e-01 -1.6513e+00
5.3909e-01 -1.6443e+00
5.7669e-01 -1.6337e+00
6.1548e-01 -1.6190e+00
6.5582e-01 -1.5994e+00
6.9775e-01 -1.5743e+00
7.4095e-01 -1.5429e+00
7.8496e-01 -1.5049e+00
8.2920e-01 -1.4599e+00
8.7307e-01 -1.4078e+00
9.1604e-01 -1.3487e+00
9.5764e-01 -1.2829e+00
9.9747e-01 -1.2105e+00
1.0352e+00 -1.1321e+00
1.0706e+00 -1.0480e+00
1.1033e+00 -9.5863e-01
1.1332e+00 -8.6448e-01
1.1600e+00 -7.6603e-01
1.1836e+00 -6.6379e-01
1.2038e+00 -5.5825e-01
1.2205e+00 -4.4995e-01
1.2336e+00 -3.3944e-01
1.2430e+00 -2.2725e-01
1.2487e+00 -1.1394e-01
1.2506e+00 -6.5281e-05
1.2487e+00 1.1381e-01
1.2430e+00 2.2712e-01
1.2336e+00 3.3931e-01
1.2205e+00 4.4982e-01
1.2038e+00 5.5812e-01
1.1836e+00 6.6365e-01
1.1600e+00 7.6590e-01
1.1332e+00 8.6435e-01
1.1033e+00 9.5850e-01
1.0706e+00 1.0479e+00
1.0352e+00 1.1320e+00
9.9747e-01 1.2104e+00
9.5764e-01 1.2827e+00
9.1604e-01 1.3486e+00
8.7307e-01 1.4077e+00
8.2919e-01 1.4598e+00
7.8496e-01 1.5048e+00
7.4095e-01 1.5428e+00
6.9775e-01 1.5741e+00
6.5582e-01 1.5993e+00
6.1548e-01 1.6189e+00
5.7669e-01 1.6336e+00
5.3908e-01 1.6442e+00
5.0194e-01 1.6512e+00
};

\end{axis}

\end{tikzpicture} 
  \fi
  \end{tabular}
\end{center}
\mcaption{Two vesicles discretized with $N=96$ points.  The vesicles are
initially placed symmetrically around the origin and the background
velocity is $\vv_{\infty} = (-x,y)$.}{f:extensionalGeom} 
\end{figure}

As the vesicles approach one another, their shape is nearly static and
their velocities decrease.  Therefore, we test our adaptive time
stepping strategy using zero, one, and two SDC corrections.  We expect
that larger time steps can be taken as the vesicles come closer
together.  The errors in area and length, the number of accepted and
rejected time steps, the number of FMM calls, and the CPU times are
reported in Table~\ref{t:adaptiveFirstOrderExtensional} ($n_{\sdc}=0$),
Table~\ref{t:adaptiveSecondOrderExtensional} ($n_{\sdc}=1$), and
Table~\ref{t:adaptiveThirdOrderExtensional} ($n_{\sdc}=2$).  The
adaptive time stepping strategy does a good job of attaining the desired
tolerance while not having too many rejected time steps.  In
Figure~\ref{f:extensionalSummary}, we plot the time step size, the
location of the rejected time steps, and the errors in area and length
for a constant time step size and for an adaptive time step size.  While
the constant time step size commits a negligible amount of error shortly
after the initial condition, too much error is accumulated at the
beginning of the simulation.  This indicates that a smaller time step
should be taken near the start of the simulation, and then it can be
increased later in the simulation.  The left plot of
Figure~\ref{f:extensionalSummary} exactly demonstrates this behaviour.

Comparing
Tables~\ref{t:noSDCextensional}--\ref{t:adaptiveThirdOrderExtensional},
we observe the following behaviors:
\begin{itemize}
  \item Unsurprisingly, the errors resulting from the first-order time
  integrator are much larger than those resulting from higher-order
  time integrators.  These larger errors have a very adverse effect
  when using adaptive time stepping since very small time steps must be
  taken to maintain the requested local truncation error.

  \item As we saw in the {\em relaxation} example, the two second-order
  methods (BDF and $n_{\sdc}=1$) achieve similar errors with respect to
  CPU time.  Again, only the integrator $n_{\sdc}=1$ is compatible with
  adaptive time stepping, so we no longer report results using BDF.

  \item Considering adaptive time steps, to achieve a four digits of
  accuracy, our second-order adaptive time integrator is 11 times
  faster than our first-order adaptive time integrator.  Furthermore,
  if seven digits of accuracy is requested, our third-order adaptive
  time integrator is 23\% faster.


  \item In general, when smaller tolerances are requested, additional
  SDC iterations, rather than smaller time step sizes, should be used
  to increase the accuracy.  For instance, to achieve seven digits of
  accuracy, using two rather than one SDC correction results in a 25\%
  savings in CPU time.  This behaviour was also observed by
  Minion~\cite{min2003} for the Van der Pol equation: ``higher-order
  methods are again more efficient when higher precision is required."

\end{itemize}

\begin{table}[htps]
\begin{center}
\begin{tabular}{c|cccc|cccc}
 & \multicolumn{4}{c|}{$\boldnsdc{0}$} & 
   \multicolumn{4}{c}{{\bf BDF}} \\
$m$ & $e_{A}$ & $e_{L}$ & \# fmm & CPU & 
      $e_{A}$ & $e_{L}$ & \# fmm & CPU \\
\hline
$300$  & $2.46$E$-4$ & $1.27$E$-3$ & $3.99$E$3$ & $1.0$
       & $5.76$E$-5$ & $5.67$E$-4$ & $3.88$E$3$ & $0.9$ \\
$600$  & $1.24$E$-4$ & $6.64$E$-4$ & $7.53$E$3$ & $1.8$
       & $1.55$E$-5$ & $2.00$E$-4$ & $7.34$E$3$ & $1.8$ \\ 
$1200$ & $6.24$E$-5$ & $3.46$E$-4$ & $1.44$E$4$ & $3.6$
       & $4.12$E$-6$ & $6.54$E$-5$ & $1.40$E$4$ & $3.6$ \\ 
$2400$ & $3.15$E$-5$ & $1.79$E$-4$ & $2.75$E$4$ & $6.8$
       & $1.06$E$-6$ & $1.92$E$-5$ & $2.71$E$4$ & $6.7$
\end{tabular}
\mcaption{The errors in area and length and the CPU time for two
vesicles in an {\bf extensional} flow with a {\bf constant} time step
size using $\boldnsdc{0}$ (left) and ${\bf BDF}$ (right).  The CPU
times for
Tables~\ref{t:noSDCextensional}--\ref{t:adaptiveThirdOrderExtensional}
are relative to the cheapest simulation ($m=300$ and $n_{\sdc}=0$)
which took approximately 613 seconds.  We achieve the desired
first-order results, but second-order results are not yet achieved.
With additional time steps, second-order convergence is achieved (see
Table $7$ in~\cite{qua:bir2013b}).}{t:noSDCextensional}

\begin{tabular}{c|cccc|cccc}
 & \multicolumn{4}{c|}{$\boldnsdc{1}$} & 
   \multicolumn{4}{c}{$\boldnsdc{2}$} \\ 
$m$ & $e_{A}$ & $e_{L}$ & \# fmm & CPU & 
      $e_{A}$ & $e_{L}$ & \# fmm & CPU \\
\hline
$300$  & $1.37$E$-6$ & $3.08$E$-6$  & $2.38$E$4$ & $4.6$
       & $7.39$E$-7$ & $7.16$E$-8$  & $3.61$E$4$ & $7.1$ \\
$600$  & $1.03$E$-6$ & $1.32$E$-6$  & $4.59$E$4$ & $8.9$
       & $1.56$E$-7$ & $8.78$E$-9$  & $6.95$E$4$ & $14$ \\ 
$1200$ & $4.51$E$-7$ & $5.24$E$-7$  & $8.91$E$4$ & $18$
       & $3.26$E$-8$ & $2.23$E$-9$  & $1.35$E$5$ & $26$ \\ 
$2400$ & $1.60$E$-7$ & $1.76$E$-7$  & $1.75$E$5$ & $34$
       & $8.39$E$-9$ & $7.25$E$-10$ & $2.65$E$5$ & $52$ \\
$4800$ & $4.95$E$-8$ & $4.92$E$-8$  & $3.46$E$5$ & $67$
       & $2.90$E$-9$ & $1.94$E$-10$ & $5.28$E$5$ & $104$
\end{tabular}
\mcaption{The errors in area and length and the CPU time for two
vesicles in an {\bf extensional} flow with a {\bf constant} time step
size using $\boldnsdc{1}$ (left) and $\boldnsdc{2}$ (right).  While
each SDC correction reduces the error, the desired asymptotic rates of
convergence are not achieved.  As the ratios of successive errors are
approaching the expected values, we expect that we have not taken
enough time steps to observe the asymptotic rate of convergence.  Also,
we observe that one SDC correction and BDF have comparable errors with
respect to computational work.}{t:SDC12extensional}

\begin{tabular}{ccccccc}
Tolerance & $e_{A}$ & $e_{L}$ & Accepts & Rejects & \# fmm & CPU \\
\hline
$1$E$-2$ & $2.10$E$-4$ & $1.34$E$-3$ 
         & $76$   & $28$ & $1.16$E$3$ & $0.3$ \\
$1$E$-3$ & $3.55$E$-5$ & $2.22$E$-4$ 
         & $510$  & $30$ & $4.69$E$3$ & $1.4$ \\
$1$E$-4$ & $6.73$E$-6$ & $4.11$E$-5$ 
         & $4803$ & $34$ & $3.98$E$4$ & $12$ 
\end{tabular}
\mcaption{The errors in area and length, the CPU time, and the number
of accepted and rejected time steps for two vesicles in an {\bf
extensional} flow with an {\bf adaptive} time step size using
$\boldnsdc{0}$.  For the larger tolerances, the desired tolerance is
achieved in an acceptable amount of CPU time.  However, first-order
methods require too small of time steps to achieve four digits of
accuracy.  With these smaller tolerances, higher-order methods should
be used.}{t:adaptiveFirstOrderExtensional}

\begin{tabular}{ccccccc}
Tolerance & $e_{A}$ & $e_{L}$ & Accepts & Rejects & \# fmm & CPU \\
\hline $1$E$-2$ & $4.24$E$-3$ & $2.93$E$-4$ & $26$  & $18$ & $5.12$E$3$
& $0.9$ \\ $1$E$-3$ & $2.90$E$-4$ & $8.42$E$-5$ & $28$  & $16$ &
$4.32$E$3$ & $0.8$ \\ $1$E$-4$ & $1.01$E$-6$ & $7.10$E$-6$ & $45$  &
$24$ & $5.90$E$3$ & $1.1$ \\ $1$E$-5$ & $7.08$E$-6$ & $8.97$E$-7$ 
         & $97$  & $16$ & $8.14$E$3$ & $1.5$ \\
$1$E$-6$ & $8.79$E$-7$ & $4.40$E$-7$ 
         & $289$ & $21$ & $2.13$E$4$ & $4.1$ \\
$1$E$-7$ & $9.18$E$-8$ & $3.04$E$-9$ 
         & $891$ & $24$ & $6.07$E$4$ & $12$
\end{tabular}
\mcaption{The errors in area and length, the CPU time, and the number of
accepted and rejected time steps for two vesicles in an {\bf
extensional} flow with an {\bf adaptive} time step size and
$\boldnsdc{1}$.  This second-order method is able to achieve much
smaller tolerances than $n_{\sdc}=0$.  However, when seven digits
of accuracy is requested, the number of required time steps becomes
unacceptable.  In this case, a third-order method should be
used.}{t:adaptiveSecondOrderExtensional}

\begin{tabular}{ccccccc}
Tolerance & $e_{A}$ & $e_{L}$ & Accepts & Rejects & \# fmm & CPU \\
\hline
$1$E$-5$ & $3.07$E$-6$ & $1.53$E$-6$ 
         & $59$  & $23$ & $1.06$E$4$ & $2.0$ \\
$1$E$-6$ & $6.48$E$-7$ & $1.11$E$-7$ 
         & $143$ & $33$ & $2.06$E$4$ & $4.0$ \\
$1$E$-7$ & $8.90$E$-8$ & $1.68$E$-9$ 
         & $430$ & $22$ & $4.75$E$4$ & $9.2$ 
\end{tabular}
\mcaption{The errors in area and length, the CPU time, and the number of
accepted and rejected time steps for two vesicles in an {\bf
extensional} flow with an {\bf adaptive} time step size using
$\boldnsdc{2}$.  We see that if smaller tolerances are desired, it is
advantageous to use additional SDC corrections to allow for larger time
steps.}{t:adaptiveThirdOrderExtensional}

\end{center}
\end{table}

\begin{figure}[htps]
\begin{tabular}{cc}
\ifTikz
\input{extensionalAdaptiveDT.tikz} &
\input{extensionalAdaptiveErrors.tikz}
\fi
\end{tabular}
\mcaption{Results for the {\bf extensional} flow using $\boldnsdc{1}$.
Left: The time step size using 2400 {\bf constant} time step sizes
(dashed) and an {\bf adaptive} time step size (solid).  The open
circles indicate the 24 times when the time step size is rejected.  As
expected, the time step size increases as the vesicles come closer to a
steady state.  Right: The errors in area and length using {\bf
constant} (dashed) and {\bf adaptive} (solid) time steps, and the
desired error (black) of the adaptive time step.  When using adaptive
time stepping, the error in area nearly achieves the desired error
indicating that we are almost selecting the optimal time step size.
However, when using a constant time step size, a large amount of error
is committed at the start of the simulation, and then very little error
is committed for the duration of the simulation.  The result is a
nearly three-fold increase in the CPU time.  The CPU time is further
reduced by using adaptive time stepping with
$n_{\sdc}=2$.}{f:extensionalSummary}
\end{figure}

\subsection{Stenosis}

We consider a single vesicle discretized with $N=128$ points in a
constricted tube discretized with $N_{\mathrm{wall}}=256$ points
(Figure~\ref{f:stenosisGeom}).  At this resolution, our FMM
implementation of the double-layer potential is slower than a direct
evaluation.  Therefore, the FMM is only used for the single-layer
potentials.  The time horizon is $T=15$ which is sufficiently long that
the vesicle passes through the constriction.  We again use $p=4$
Gauss-Lobatto quadrature points.  We check the rates of convergence for
a varying number of SDC corrections in Tables~\ref{t:SDC01stenosis}
and~\ref{t:SDC23stenosis}.  We see that first- and second-order
convergence is achieved in Table~\ref{t:SDC01stenosis}.  While the
error continues to decrease with each SDC iteration, as before,
additional orders of convergence are not achieved for the presented
values of $\Delta t$.

\begin{figure}[htps]
\centering
\begin{tabular}{ccc}
\ifTikz
\begin{tikzpicture}[scale=0.5]

\begin{axis}[
  xmin = -11,
  xmax = 11,
  ymin = -3.2,
  ymax = 3.2,
  scale only axis,
  axis equal image,
  hide axis,
  title = {\Huge$t=5$}
  ]

\addplot [mark=none,black,line width=1.5] table{
1.0000e+01 0.0000e+00
1.0000e+01 1.2574e-01
1.0000e+01 2.5192e-01
1.0000e+01 3.7899e-01
1.0000e+01 5.0741e-01
1.0000e+01 6.3767e-01
1.0000e+01 7.7027e-01
9.9999e+00 9.0575e-01
9.9997e+00 1.0447e+00
9.9992e+00 1.1877e+00
9.9981e+00 1.3354e+00
9.9954e+00 1.4885e+00
9.9896e+00 1.6476e+00
9.9776e+00 1.8128e+00
9.9535e+00 1.9839e+00
9.9070e+00 2.1594e+00
9.8736e+00 2.2404e+00
9.8298e+00 2.3209e+00
9.7729e+00 2.4002e+00
9.7001e+00 2.4773e+00
9.6085e+00 2.5512e+00
9.4956e+00 2.6206e+00
9.3594e+00 2.6846e+00
9.1991e+00 2.7421e+00
9.0150e+00 2.7926e+00
8.8086e+00 2.8359e+00
8.5829e+00 2.8720e+00
8.3412e+00 2.9015e+00
8.0874e+00 2.9252e+00
7.8251e+00 2.9437e+00
7.5577e+00 2.9581e+00
7.2879e+00 2.9691e+00
7.0179e+00 2.9773e+00
6.7494e+00 2.9835e+00
6.4835e+00 2.9881e+00
6.2211e+00 2.9915e+00
5.9627e+00 2.9940e+00
5.7086e+00 2.9957e+00
5.4588e+00 2.9970e+00
5.2134e+00 2.9979e+00
4.9724e+00 2.9986e+00
4.7355e+00 2.9991e+00
4.5026e+00 2.9994e+00
4.2736e+00 2.9996e+00
4.0482e+00 2.9997e+00
3.8262e+00 2.9998e+00
3.6073e+00 2.9999e+00
3.3915e+00 2.9999e+00
3.1784e+00 3.0000e+00
2.9679e+00 2.9831e+00
2.7598e+00 2.9190e+00
2.5539e+00 2.8112e+00
2.3500e+00 2.6655e+00
2.1479e+00 2.4888e+00
1.9475e+00 2.2888e+00
1.7486e+00 2.0740e+00
1.5510e+00 1.8527e+00
1.3546e+00 1.6337e+00
1.1592e+00 1.4249e+00
9.6469e-01 1.2341e+00
7.7089e-01 1.0680e+00
5.7767e-01 9.3254e-01
3.8487e-01 8.3230e-01
1.9237e-01 7.7075e-01
6.1232e-16 7.5000e-01
-1.9237e-01 7.7075e-01
-3.8487e-01 8.3230e-01
-5.7767e-01 9.3254e-01
-7.7089e-01 1.0680e+00
-9.6469e-01 1.2341e+00
-1.1592e+00 1.4249e+00
-1.3546e+00 1.6337e+00
-1.5510e+00 1.8527e+00
-1.7486e+00 2.0740e+00
-1.9475e+00 2.2888e+00
-2.1479e+00 2.4888e+00
-2.3500e+00 2.6655e+00
-2.5539e+00 2.8112e+00
-2.7598e+00 2.9190e+00
-2.9679e+00 2.9831e+00
-3.1784e+00 3.0000e+00
-3.3915e+00 2.9999e+00
-3.6073e+00 2.9999e+00
-3.8262e+00 2.9998e+00
-4.0482e+00 2.9997e+00
-4.2736e+00 2.9996e+00
-4.5026e+00 2.9994e+00
-4.7355e+00 2.9991e+00
-4.9724e+00 2.9986e+00
-5.2134e+00 2.9979e+00
-5.4588e+00 2.9970e+00
-5.7086e+00 2.9957e+00
-5.9627e+00 2.9940e+00
-6.2211e+00 2.9915e+00
-6.4835e+00 2.9881e+00
-6.7494e+00 2.9835e+00
-7.0179e+00 2.9773e+00
-7.2879e+00 2.9691e+00
-7.5577e+00 2.9581e+00
-7.8251e+00 2.9437e+00
-8.0874e+00 2.9252e+00
-8.3412e+00 2.9015e+00
-8.5829e+00 2.8720e+00
-8.8086e+00 2.8359e+00
-9.0150e+00 2.7926e+00
-9.1991e+00 2.7421e+00
-9.3594e+00 2.6846e+00
-9.4956e+00 2.6206e+00
-9.6085e+00 2.5512e+00
-9.7001e+00 2.4773e+00
-9.7729e+00 2.4002e+00
-9.8298e+00 2.3209e+00
-9.8736e+00 2.2404e+00
-9.9070e+00 2.1594e+00
-9.9535e+00 1.9839e+00
-9.9776e+00 1.8128e+00
-9.9896e+00 1.6476e+00
-9.9954e+00 1.4885e+00
-9.9981e+00 1.3354e+00
-9.9992e+00 1.1877e+00
-9.9997e+00 1.0447e+00
-9.9999e+00 9.0575e-01
-1.0000e+01 7.7027e-01
-1.0000e+01 6.3767e-01
-1.0000e+01 5.0741e-01
-1.0000e+01 3.7899e-01
-1.0000e+01 2.5192e-01
-1.0000e+01 1.2574e-01
-1.0000e+01 3.6739e-16
-1.0000e+01 -1.2574e-01
-1.0000e+01 -2.5192e-01
-1.0000e+01 -3.7899e-01
-1.0000e+01 -5.0741e-01
-1.0000e+01 -6.3767e-01
-1.0000e+01 -7.7027e-01
-9.9999e+00 -9.0575e-01
-9.9997e+00 -1.0447e+00
-9.9992e+00 -1.1877e+00
-9.9981e+00 -1.3354e+00
-9.9954e+00 -1.4885e+00
-9.9896e+00 -1.6476e+00
-9.9776e+00 -1.8128e+00
-9.9535e+00 -1.9839e+00
-9.9070e+00 -2.1594e+00
-9.8736e+00 -2.2404e+00
-9.8298e+00 -2.3209e+00
-9.7729e+00 -2.4002e+00
-9.7001e+00 -2.4773e+00
-9.6085e+00 -2.5512e+00
-9.4956e+00 -2.6206e+00
-9.3594e+00 -2.6846e+00
-9.1991e+00 -2.7421e+00
-9.0150e+00 -2.7926e+00
-8.8086e+00 -2.8359e+00
-8.5829e+00 -2.8720e+00
-8.3412e+00 -2.9015e+00
-8.0874e+00 -2.9252e+00
-7.8251e+00 -2.9437e+00
-7.5577e+00 -2.9581e+00
-7.2879e+00 -2.9691e+00
-7.0179e+00 -2.9773e+00
-6.7494e+00 -2.9835e+00
-6.4835e+00 -2.9881e+00
-6.2211e+00 -2.9915e+00
-5.9627e+00 -2.9940e+00
-5.7086e+00 -2.9957e+00
-5.4588e+00 -2.9970e+00
-5.2134e+00 -2.9979e+00
-4.9724e+00 -2.9986e+00
-4.7355e+00 -2.9991e+00
-4.5026e+00 -2.9994e+00
-4.2736e+00 -2.9996e+00
-4.0482e+00 -2.9997e+00
-3.8262e+00 -2.9998e+00
-3.6073e+00 -2.9999e+00
-3.3915e+00 -2.9999e+00
-3.1784e+00 -3.0000e+00
-2.9679e+00 -2.9831e+00
-2.7598e+00 -2.9190e+00
-2.5539e+00 -2.8112e+00
-2.3500e+00 -2.6655e+00
-2.1479e+00 -2.4888e+00
-1.9475e+00 -2.2888e+00
-1.7486e+00 -2.0740e+00
-1.5510e+00 -1.8527e+00
-1.3546e+00 -1.6337e+00
-1.1592e+00 -1.4249e+00
-9.6469e-01 -1.2341e+00
-7.7089e-01 -1.0680e+00
-5.7767e-01 -9.3254e-01
-3.8487e-01 -8.3230e-01
-1.9237e-01 -7.7075e-01
-1.8370e-15 -7.5000e-01
1.9237e-01 -7.7075e-01
3.8487e-01 -8.3230e-01
5.7767e-01 -9.3254e-01
7.7089e-01 -1.0680e+00
9.6469e-01 -1.2341e+00
1.1592e+00 -1.4249e+00
1.3546e+00 -1.6337e+00
1.5510e+00 -1.8527e+00
1.7486e+00 -2.0740e+00
1.9475e+00 -2.2888e+00
2.1479e+00 -2.4888e+00
2.3500e+00 -2.6655e+00
2.5539e+00 -2.8112e+00
2.7598e+00 -2.9190e+00
2.9679e+00 -2.9831e+00
3.1784e+00 -3.0000e+00
3.3915e+00 -2.9999e+00
3.6073e+00 -2.9999e+00
3.8262e+00 -2.9998e+00
4.0482e+00 -2.9997e+00
4.2736e+00 -2.9996e+00
4.5026e+00 -2.9994e+00
4.7355e+00 -2.9991e+00
4.9724e+00 -2.9986e+00
5.2134e+00 -2.9979e+00
5.4588e+00 -2.9970e+00
5.7086e+00 -2.9957e+00
5.9627e+00 -2.9940e+00
6.2211e+00 -2.9915e+00
6.4835e+00 -2.9881e+00
6.7494e+00 -2.9835e+00
7.0179e+00 -2.9773e+00
7.2879e+00 -2.9691e+00
7.5577e+00 -2.9581e+00
7.8251e+00 -2.9437e+00
8.0874e+00 -2.9252e+00
8.3412e+00 -2.9015e+00
8.5829e+00 -2.8720e+00
8.8086e+00 -2.8359e+00
9.0150e+00 -2.7926e+00
9.1991e+00 -2.7421e+00
9.3594e+00 -2.6846e+00
9.4956e+00 -2.6206e+00
9.6085e+00 -2.5512e+00
9.7001e+00 -2.4773e+00
9.7729e+00 -2.4002e+00
9.8298e+00 -2.3209e+00
9.8736e+00 -2.2404e+00
9.9070e+00 -2.1594e+00
9.9535e+00 -1.9839e+00
9.9776e+00 -1.8128e+00
9.9896e+00 -1.6476e+00
9.9954e+00 -1.4885e+00
9.9981e+00 -1.3354e+00
9.9992e+00 -1.1877e+00
9.9997e+00 -1.0447e+00
9.9999e+00 -9.0575e-01
1.0000e+01 -7.7027e-01
1.0000e+01 -6.3767e-01
1.0000e+01 -5.0741e-01
1.0000e+01 -3.7899e-01
1.0000e+01 -2.5192e-01
1.0000e+01 -1.2574e-01
1.0000e+01 0.0000e+00
};

\addplot [mark=none,red,line width=1.5] table{
-5.1031e+00 1.5091e+00
-5.1305e+00 1.5158e+00
-5.1589e+00 1.5210e+00
-5.1886e+00 1.5247e+00
-5.2200e+00 1.5263e+00
-5.2533e+00 1.5256e+00
-5.2886e+00 1.5218e+00
-5.3258e+00 1.5140e+00
-5.3643e+00 1.5016e+00
-5.4035e+00 1.4835e+00
-5.4423e+00 1.4588e+00
-5.4792e+00 1.4272e+00
-5.5129e+00 1.3881e+00
-5.5414e+00 1.3419e+00
-5.5633e+00 1.2892e+00
-5.5774e+00 1.2312e+00
-5.5829e+00 1.1691e+00
-5.5797e+00 1.1043e+00
-5.5683e+00 1.0380e+00
-5.5493e+00 9.7099e-01
-5.5240e+00 9.0383e-01
-5.4936e+00 8.3658e-01
-5.4594e+00 7.6906e-01
-5.4227e+00 7.0088e-01
-5.3848e+00 6.3154e-01
-5.3471e+00 5.6050e-01
-5.3109e+00 4.8733e-01
-5.2774e+00 4.1169e-01
-5.2479e+00 3.3345e-01
-5.2235e+00 2.5271e-01
-5.2053e+00 1.6981e-01
-5.1940e+00 8.5329e-02
-5.1902e+00 4.1535e-14
-5.1940e+00 -8.5329e-02
-5.2053e+00 -1.6981e-01
-5.2235e+00 -2.5271e-01
-5.2479e+00 -3.3345e-01
-5.2774e+00 -4.1169e-01
-5.3109e+00 -4.8733e-01
-5.3471e+00 -5.6050e-01
-5.3848e+00 -6.3154e-01
-5.4227e+00 -7.0088e-01
-5.4594e+00 -7.6906e-01
-5.4936e+00 -8.3658e-01
-5.5240e+00 -9.0383e-01
-5.5493e+00 -9.7099e-01
-5.5683e+00 -1.0380e+00
-5.5797e+00 -1.1043e+00
-5.5829e+00 -1.1691e+00
-5.5774e+00 -1.2312e+00
-5.5633e+00 -1.2892e+00
-5.5414e+00 -1.3419e+00
-5.5129e+00 -1.3881e+00
-5.4792e+00 -1.4272e+00
-5.4423e+00 -1.4588e+00
-5.4035e+00 -1.4835e+00
-5.3643e+00 -1.5016e+00
-5.3258e+00 -1.5140e+00
-5.2886e+00 -1.5218e+00
-5.2533e+00 -1.5256e+00
-5.2200e+00 -1.5263e+00
-5.1886e+00 -1.5247e+00
-5.1589e+00 -1.5210e+00
-5.1305e+00 -1.5158e+00
-5.1031e+00 -1.5091e+00
-5.0760e+00 -1.5009e+00
-5.0488e+00 -1.4913e+00
-5.0211e+00 -1.4800e+00
-4.9926e+00 -1.4668e+00
-4.9630e+00 -1.4515e+00
-4.9322e+00 -1.4338e+00
-4.9002e+00 -1.4136e+00
-4.8668e+00 -1.3907e+00
-4.8322e+00 -1.3649e+00
-4.7963e+00 -1.3362e+00
-4.7593e+00 -1.3044e+00
-4.7212e+00 -1.2696e+00
-4.6821e+00 -1.2319e+00
-4.6421e+00 -1.1912e+00
-4.6012e+00 -1.1476e+00
-4.5595e+00 -1.1012e+00
-4.5171e+00 -1.0520e+00
-4.4742e+00 -1.0002e+00
-4.4310e+00 -9.4566e-01
-4.3876e+00 -8.8852e-01
-4.3443e+00 -8.2876e-01
-4.3015e+00 -7.6635e-01
-4.2597e+00 -7.0125e-01
-4.2192e+00 -6.3341e-01
-4.1808e+00 -5.6278e-01
-4.1452e+00 -4.8935e-01
-4.1130e+00 -4.1316e-01
-4.0851e+00 -3.3433e-01
-4.0624e+00 -2.5313e-01
-4.0455e+00 -1.6995e-01
-4.0351e+00 -8.5346e-02
-4.0316e+00 6.3612e-13
-4.0351e+00 8.5346e-02
-4.0455e+00 1.6995e-01
-4.0624e+00 2.5313e-01
-4.0851e+00 3.3433e-01
-4.1130e+00 4.1316e-01
-4.1452e+00 4.8935e-01
-4.1808e+00 5.6278e-01
-4.2192e+00 6.3341e-01
-4.2597e+00 7.0125e-01
-4.3015e+00 7.6635e-01
-4.3443e+00 8.2876e-01
-4.3876e+00 8.8852e-01
-4.4310e+00 9.4566e-01
-4.4742e+00 1.0002e+00
-4.5171e+00 1.0520e+00
-4.5595e+00 1.1012e+00
-4.6012e+00 1.1476e+00
-4.6421e+00 1.1912e+00
-4.6821e+00 1.2319e+00
-4.7212e+00 1.2696e+00
-4.7593e+00 1.3044e+00
-4.7963e+00 1.3362e+00
-4.8322e+00 1.3649e+00
-4.8668e+00 1.3907e+00
-4.9002e+00 1.4136e+00
-4.9322e+00 1.4338e+00
-4.9630e+00 1.4515e+00
-4.9926e+00 1.4668e+00
-5.0211e+00 1.4800e+00
-5.0488e+00 1.4913e+00
-5.0760e+00 1.5009e+00
-5.1031e+00 1.5091e+00
};

\end{axis}

\end{tikzpicture} &
\begin{tikzpicture}[scale=0.5]

\begin{axis}[
  xmin = -11,
  xmax = 11,
  ymin = -3.2,
  ymax = 3.2,
  scale only axis,
  axis equal image,
  hide axis,
  title = {\Huge$t=6$}
  ]

\addplot [mark=none,black,line width=1.5] table{
1.0000e+01 0.0000e+00
1.0000e+01 1.2574e-01
1.0000e+01 2.5192e-01
1.0000e+01 3.7899e-01
1.0000e+01 5.0741e-01
1.0000e+01 6.3767e-01
1.0000e+01 7.7027e-01
9.9999e+00 9.0575e-01
9.9997e+00 1.0447e+00
9.9992e+00 1.1877e+00
9.9981e+00 1.3354e+00
9.9954e+00 1.4885e+00
9.9896e+00 1.6476e+00
9.9776e+00 1.8128e+00
9.9535e+00 1.9839e+00
9.9070e+00 2.1594e+00
9.8736e+00 2.2404e+00
9.8298e+00 2.3209e+00
9.7729e+00 2.4002e+00
9.7001e+00 2.4773e+00
9.6085e+00 2.5512e+00
9.4956e+00 2.6206e+00
9.3594e+00 2.6846e+00
9.1991e+00 2.7421e+00
9.0150e+00 2.7926e+00
8.8086e+00 2.8359e+00
8.5829e+00 2.8720e+00
8.3412e+00 2.9015e+00
8.0874e+00 2.9252e+00
7.8251e+00 2.9437e+00
7.5577e+00 2.9581e+00
7.2879e+00 2.9691e+00
7.0179e+00 2.9773e+00
6.7494e+00 2.9835e+00
6.4835e+00 2.9881e+00
6.2211e+00 2.9915e+00
5.9627e+00 2.9940e+00
5.7086e+00 2.9957e+00
5.4588e+00 2.9970e+00
5.2134e+00 2.9979e+00
4.9724e+00 2.9986e+00
4.7355e+00 2.9991e+00
4.5026e+00 2.9994e+00
4.2736e+00 2.9996e+00
4.0482e+00 2.9997e+00
3.8262e+00 2.9998e+00
3.6073e+00 2.9999e+00
3.3915e+00 2.9999e+00
3.1784e+00 3.0000e+00
2.9679e+00 2.9831e+00
2.7598e+00 2.9190e+00
2.5539e+00 2.8112e+00
2.3500e+00 2.6655e+00
2.1479e+00 2.4888e+00
1.9475e+00 2.2888e+00
1.7486e+00 2.0740e+00
1.5510e+00 1.8527e+00
1.3546e+00 1.6337e+00
1.1592e+00 1.4249e+00
9.6469e-01 1.2341e+00
7.7089e-01 1.0680e+00
5.7767e-01 9.3254e-01
3.8487e-01 8.3230e-01
1.9237e-01 7.7075e-01
6.1232e-16 7.5000e-01
-1.9237e-01 7.7075e-01
-3.8487e-01 8.3230e-01
-5.7767e-01 9.3254e-01
-7.7089e-01 1.0680e+00
-9.6469e-01 1.2341e+00
-1.1592e+00 1.4249e+00
-1.3546e+00 1.6337e+00
-1.5510e+00 1.8527e+00
-1.7486e+00 2.0740e+00
-1.9475e+00 2.2888e+00
-2.1479e+00 2.4888e+00
-2.3500e+00 2.6655e+00
-2.5539e+00 2.8112e+00
-2.7598e+00 2.9190e+00
-2.9679e+00 2.9831e+00
-3.1784e+00 3.0000e+00
-3.3915e+00 2.9999e+00
-3.6073e+00 2.9999e+00
-3.8262e+00 2.9998e+00
-4.0482e+00 2.9997e+00
-4.2736e+00 2.9996e+00
-4.5026e+00 2.9994e+00
-4.7355e+00 2.9991e+00
-4.9724e+00 2.9986e+00
-5.2134e+00 2.9979e+00
-5.4588e+00 2.9970e+00
-5.7086e+00 2.9957e+00
-5.9627e+00 2.9940e+00
-6.2211e+00 2.9915e+00
-6.4835e+00 2.9881e+00
-6.7494e+00 2.9835e+00
-7.0179e+00 2.9773e+00
-7.2879e+00 2.9691e+00
-7.5577e+00 2.9581e+00
-7.8251e+00 2.9437e+00
-8.0874e+00 2.9252e+00
-8.3412e+00 2.9015e+00
-8.5829e+00 2.8720e+00
-8.8086e+00 2.8359e+00
-9.0150e+00 2.7926e+00
-9.1991e+00 2.7421e+00
-9.3594e+00 2.6846e+00
-9.4956e+00 2.6206e+00
-9.6085e+00 2.5512e+00
-9.7001e+00 2.4773e+00
-9.7729e+00 2.4002e+00
-9.8298e+00 2.3209e+00
-9.8736e+00 2.2404e+00
-9.9070e+00 2.1594e+00
-9.9535e+00 1.9839e+00
-9.9776e+00 1.8128e+00
-9.9896e+00 1.6476e+00
-9.9954e+00 1.4885e+00
-9.9981e+00 1.3354e+00
-9.9992e+00 1.1877e+00
-9.9997e+00 1.0447e+00
-9.9999e+00 9.0575e-01
-1.0000e+01 7.7027e-01
-1.0000e+01 6.3767e-01
-1.0000e+01 5.0741e-01
-1.0000e+01 3.7899e-01
-1.0000e+01 2.5192e-01
-1.0000e+01 1.2574e-01
-1.0000e+01 3.6739e-16
-1.0000e+01 -1.2574e-01
-1.0000e+01 -2.5192e-01
-1.0000e+01 -3.7899e-01
-1.0000e+01 -5.0741e-01
-1.0000e+01 -6.3767e-01
-1.0000e+01 -7.7027e-01
-9.9999e+00 -9.0575e-01
-9.9997e+00 -1.0447e+00
-9.9992e+00 -1.1877e+00
-9.9981e+00 -1.3354e+00
-9.9954e+00 -1.4885e+00
-9.9896e+00 -1.6476e+00
-9.9776e+00 -1.8128e+00
-9.9535e+00 -1.9839e+00
-9.9070e+00 -2.1594e+00
-9.8736e+00 -2.2404e+00
-9.8298e+00 -2.3209e+00
-9.7729e+00 -2.4002e+00
-9.7001e+00 -2.4773e+00
-9.6085e+00 -2.5512e+00
-9.4956e+00 -2.6206e+00
-9.3594e+00 -2.6846e+00
-9.1991e+00 -2.7421e+00
-9.0150e+00 -2.7926e+00
-8.8086e+00 -2.8359e+00
-8.5829e+00 -2.8720e+00
-8.3412e+00 -2.9015e+00
-8.0874e+00 -2.9252e+00
-7.8251e+00 -2.9437e+00
-7.5577e+00 -2.9581e+00
-7.2879e+00 -2.9691e+00
-7.0179e+00 -2.9773e+00
-6.7494e+00 -2.9835e+00
-6.4835e+00 -2.9881e+00
-6.2211e+00 -2.9915e+00
-5.9627e+00 -2.9940e+00
-5.7086e+00 -2.9957e+00
-5.4588e+00 -2.9970e+00
-5.2134e+00 -2.9979e+00
-4.9724e+00 -2.9986e+00
-4.7355e+00 -2.9991e+00
-4.5026e+00 -2.9994e+00
-4.2736e+00 -2.9996e+00
-4.0482e+00 -2.9997e+00
-3.8262e+00 -2.9998e+00
-3.6073e+00 -2.9999e+00
-3.3915e+00 -2.9999e+00
-3.1784e+00 -3.0000e+00
-2.9679e+00 -2.9831e+00
-2.7598e+00 -2.9190e+00
-2.5539e+00 -2.8112e+00
-2.3500e+00 -2.6655e+00
-2.1479e+00 -2.4888e+00
-1.9475e+00 -2.2888e+00
-1.7486e+00 -2.0740e+00
-1.5510e+00 -1.8527e+00
-1.3546e+00 -1.6337e+00
-1.1592e+00 -1.4249e+00
-9.6469e-01 -1.2341e+00
-7.7089e-01 -1.0680e+00
-5.7767e-01 -9.3254e-01
-3.8487e-01 -8.3230e-01
-1.9237e-01 -7.7075e-01
-1.8370e-15 -7.5000e-01
1.9237e-01 -7.7075e-01
3.8487e-01 -8.3230e-01
5.7767e-01 -9.3254e-01
7.7089e-01 -1.0680e+00
9.6469e-01 -1.2341e+00
1.1592e+00 -1.4249e+00
1.3546e+00 -1.6337e+00
1.5510e+00 -1.8527e+00
1.7486e+00 -2.0740e+00
1.9475e+00 -2.2888e+00
2.1479e+00 -2.4888e+00
2.3500e+00 -2.6655e+00
2.5539e+00 -2.8112e+00
2.7598e+00 -2.9190e+00
2.9679e+00 -2.9831e+00
3.1784e+00 -3.0000e+00
3.3915e+00 -2.9999e+00
3.6073e+00 -2.9999e+00
3.8262e+00 -2.9998e+00
4.0482e+00 -2.9997e+00
4.2736e+00 -2.9996e+00
4.5026e+00 -2.9994e+00
4.7355e+00 -2.9991e+00
4.9724e+00 -2.9986e+00
5.2134e+00 -2.9979e+00
5.4588e+00 -2.9970e+00
5.7086e+00 -2.9957e+00
5.9627e+00 -2.9940e+00
6.2211e+00 -2.9915e+00
6.4835e+00 -2.9881e+00
6.7494e+00 -2.9835e+00
7.0179e+00 -2.9773e+00
7.2879e+00 -2.9691e+00
7.5577e+00 -2.9581e+00
7.8251e+00 -2.9437e+00
8.0874e+00 -2.9252e+00
8.3412e+00 -2.9015e+00
8.5829e+00 -2.8720e+00
8.8086e+00 -2.8359e+00
9.0150e+00 -2.7926e+00
9.1991e+00 -2.7421e+00
9.3594e+00 -2.6846e+00
9.4956e+00 -2.6206e+00
9.6085e+00 -2.5512e+00
9.7001e+00 -2.4773e+00
9.7729e+00 -2.4002e+00
9.8298e+00 -2.3209e+00
9.8736e+00 -2.2404e+00
9.9070e+00 -2.1594e+00
9.9535e+00 -1.9839e+00
9.9776e+00 -1.8128e+00
9.9896e+00 -1.6476e+00
9.9954e+00 -1.4885e+00
9.9981e+00 -1.3354e+00
9.9992e+00 -1.1877e+00
9.9997e+00 -1.0447e+00
9.9999e+00 -9.0575e-01
1.0000e+01 -7.7027e-01
1.0000e+01 -6.3767e-01
1.0000e+01 -5.0741e-01
1.0000e+01 -3.7899e-01
1.0000e+01 -2.5192e-01
1.0000e+01 -1.2574e-01
1.0000e+01 0.0000e+00
};

\addplot [mark=none,red,line width=1.5] table{
-4.2858e+00 1.4447e+00
-4.3128e+00 1.4532e+00
-4.3407e+00 1.4604e+00
-4.3701e+00 1.4662e+00
-4.4012e+00 1.4704e+00
-4.4345e+00 1.4724e+00
-4.4700e+00 1.4717e+00
-4.5077e+00 1.4674e+00
-4.5472e+00 1.4587e+00
-4.5880e+00 1.4445e+00
-4.6290e+00 1.4237e+00
-4.6688e+00 1.3957e+00
-4.7058e+00 1.3598e+00
-4.7380e+00 1.3161e+00
-4.7636e+00 1.2651e+00
-4.7810e+00 1.2080e+00
-4.7892e+00 1.1462e+00
-4.7877e+00 1.0814e+00
-4.7769e+00 1.0150e+00
-4.7576e+00 9.4812e-01
-4.7307e+00 8.8155e-01
-4.6978e+00 8.1549e-01
-4.6603e+00 7.4978e-01
-4.6195e+00 6.8397e-01
-4.5769e+00 6.1742e-01
-4.5339e+00 5.4942e-01
-4.4921e+00 4.7926e-01
-4.4530e+00 4.0638e-01
-4.4180e+00 3.3042e-01
-4.3888e+00 2.5132e-01
-4.3667e+00 1.6937e-01
-4.3529e+00 8.5272e-02
-4.3482e+00 1.3664e-13
-4.3529e+00 -8.5272e-02
-4.3667e+00 -1.6937e-01
-4.3888e+00 -2.5132e-01
-4.4180e+00 -3.3042e-01
-4.4530e+00 -4.0638e-01
-4.4921e+00 -4.7926e-01
-4.5339e+00 -5.4942e-01
-4.5769e+00 -6.1742e-01
-4.6195e+00 -6.8397e-01
-4.6603e+00 -7.4978e-01
-4.6978e+00 -8.1549e-01
-4.7307e+00 -8.8155e-01
-4.7576e+00 -9.4812e-01
-4.7769e+00 -1.0150e+00
-4.7877e+00 -1.0814e+00
-4.7892e+00 -1.1462e+00
-4.7810e+00 -1.2080e+00
-4.7636e+00 -1.2651e+00
-4.7380e+00 -1.3161e+00
-4.7058e+00 -1.3598e+00
-4.6688e+00 -1.3957e+00
-4.6290e+00 -1.4237e+00
-4.5880e+00 -1.4445e+00
-4.5472e+00 -1.4587e+00
-4.5077e+00 -1.4674e+00
-4.4700e+00 -1.4717e+00
-4.4345e+00 -1.4724e+00
-4.4012e+00 -1.4704e+00
-4.3701e+00 -1.4662e+00
-4.3407e+00 -1.4604e+00
-4.3128e+00 -1.4532e+00
-4.2858e+00 -1.4447e+00
-4.2592e+00 -1.4350e+00
-4.2326e+00 -1.4239e+00
-4.2054e+00 -1.4113e+00
-4.1775e+00 -1.3970e+00
-4.1485e+00 -1.3806e+00
-4.1182e+00 -1.3621e+00
-4.0866e+00 -1.3412e+00
-4.0535e+00 -1.3179e+00
-4.0190e+00 -1.2919e+00
-3.9831e+00 -1.2632e+00
-3.9458e+00 -1.2318e+00
-3.9070e+00 -1.1978e+00
-3.8669e+00 -1.1612e+00
-3.8253e+00 -1.1220e+00
-3.7824e+00 -1.0804e+00
-3.7381e+00 -1.0364e+00
-3.6926e+00 -9.9021e-01
-3.6458e+00 -9.4182e-01
-3.5979e+00 -8.9131e-01
-3.5491e+00 -8.3872e-01
-3.4996e+00 -7.8403e-01
-3.4497e+00 -7.2716e-01
-3.3998e+00 -6.6797e-01
-3.3505e+00 -6.0626e-01
-3.3025e+00 -5.4176e-01
-3.2567e+00 -4.7416e-01
-3.2142e+00 -4.0318e-01
-3.1764e+00 -3.2862e-01
-3.1447e+00 -2.5048e-01
-3.1206e+00 -1.6911e-01
-3.1056e+00 -8.5237e-02
-3.1004e+00 7.3356e-13
-3.1056e+00 8.5237e-02
-3.1206e+00 1.6911e-01
-3.1447e+00 2.5048e-01
-3.1764e+00 3.2862e-01
-3.2142e+00 4.0318e-01
-3.2567e+00 4.7416e-01
-3.3025e+00 5.4176e-01
-3.3505e+00 6.0626e-01
-3.3998e+00 6.6797e-01
-3.4497e+00 7.2716e-01
-3.4996e+00 7.8403e-01
-3.5491e+00 8.3872e-01
-3.5979e+00 8.9131e-01
-3.6458e+00 9.4182e-01
-3.6926e+00 9.9021e-01
-3.7381e+00 1.0364e+00
-3.7824e+00 1.0804e+00
-3.8253e+00 1.1220e+00
-3.8669e+00 1.1612e+00
-3.9070e+00 1.1978e+00
-3.9458e+00 1.2318e+00
-3.9831e+00 1.2632e+00
-4.0190e+00 1.2919e+00
-4.0535e+00 1.3179e+00
-4.0866e+00 1.3412e+00
-4.1182e+00 1.3621e+00
-4.1485e+00 1.3806e+00
-4.1775e+00 1.3970e+00
-4.2054e+00 1.4113e+00
-4.2326e+00 1.4239e+00
-4.2592e+00 1.4350e+00
-4.2858e+00 1.4447e+00
};

\end{axis}

\end{tikzpicture} &
\begin{tikzpicture}[scale=0.5]

\begin{axis}[
  xmin = -11,
  xmax = 11,
  ymin = -3.2,
  ymax = 3.2,
  scale only axis,
  axis equal image,
  hide axis,
  title = {\Huge$t=7$}
  ]

\addplot [mark=none,black,line width=1.5] table{
1.0000e+01 0.0000e+00
1.0000e+01 1.2574e-01
1.0000e+01 2.5192e-01
1.0000e+01 3.7899e-01
1.0000e+01 5.0741e-01
1.0000e+01 6.3767e-01
1.0000e+01 7.7027e-01
9.9999e+00 9.0575e-01
9.9997e+00 1.0447e+00
9.9992e+00 1.1877e+00
9.9981e+00 1.3354e+00
9.9954e+00 1.4885e+00
9.9896e+00 1.6476e+00
9.9776e+00 1.8128e+00
9.9535e+00 1.9839e+00
9.9070e+00 2.1594e+00
9.8736e+00 2.2404e+00
9.8298e+00 2.3209e+00
9.7729e+00 2.4002e+00
9.7001e+00 2.4773e+00
9.6085e+00 2.5512e+00
9.4956e+00 2.6206e+00
9.3594e+00 2.6846e+00
9.1991e+00 2.7421e+00
9.0150e+00 2.7926e+00
8.8086e+00 2.8359e+00
8.5829e+00 2.8720e+00
8.3412e+00 2.9015e+00
8.0874e+00 2.9252e+00
7.8251e+00 2.9437e+00
7.5577e+00 2.9581e+00
7.2879e+00 2.9691e+00
7.0179e+00 2.9773e+00
6.7494e+00 2.9835e+00
6.4835e+00 2.9881e+00
6.2211e+00 2.9915e+00
5.9627e+00 2.9940e+00
5.7086e+00 2.9957e+00
5.4588e+00 2.9970e+00
5.2134e+00 2.9979e+00
4.9724e+00 2.9986e+00
4.7355e+00 2.9991e+00
4.5026e+00 2.9994e+00
4.2736e+00 2.9996e+00
4.0482e+00 2.9997e+00
3.8262e+00 2.9998e+00
3.6073e+00 2.9999e+00
3.3915e+00 2.9999e+00
3.1784e+00 3.0000e+00
2.9679e+00 2.9831e+00
2.7598e+00 2.9190e+00
2.5539e+00 2.8112e+00
2.3500e+00 2.6655e+00
2.1479e+00 2.4888e+00
1.9475e+00 2.2888e+00
1.7486e+00 2.0740e+00
1.5510e+00 1.8527e+00
1.3546e+00 1.6337e+00
1.1592e+00 1.4249e+00
9.6469e-01 1.2341e+00
7.7089e-01 1.0680e+00
5.7767e-01 9.3254e-01
3.8487e-01 8.3230e-01
1.9237e-01 7.7075e-01
6.1232e-16 7.5000e-01
-1.9237e-01 7.7075e-01
-3.8487e-01 8.3230e-01
-5.7767e-01 9.3254e-01
-7.7089e-01 1.0680e+00
-9.6469e-01 1.2341e+00
-1.1592e+00 1.4249e+00
-1.3546e+00 1.6337e+00
-1.5510e+00 1.8527e+00
-1.7486e+00 2.0740e+00
-1.9475e+00 2.2888e+00
-2.1479e+00 2.4888e+00
-2.3500e+00 2.6655e+00
-2.5539e+00 2.8112e+00
-2.7598e+00 2.9190e+00
-2.9679e+00 2.9831e+00
-3.1784e+00 3.0000e+00
-3.3915e+00 2.9999e+00
-3.6073e+00 2.9999e+00
-3.8262e+00 2.9998e+00
-4.0482e+00 2.9997e+00
-4.2736e+00 2.9996e+00
-4.5026e+00 2.9994e+00
-4.7355e+00 2.9991e+00
-4.9724e+00 2.9986e+00
-5.2134e+00 2.9979e+00
-5.4588e+00 2.9970e+00
-5.7086e+00 2.9957e+00
-5.9627e+00 2.9940e+00
-6.2211e+00 2.9915e+00
-6.4835e+00 2.9881e+00
-6.7494e+00 2.9835e+00
-7.0179e+00 2.9773e+00
-7.2879e+00 2.9691e+00
-7.5577e+00 2.9581e+00
-7.8251e+00 2.9437e+00
-8.0874e+00 2.9252e+00
-8.3412e+00 2.9015e+00
-8.5829e+00 2.8720e+00
-8.8086e+00 2.8359e+00
-9.0150e+00 2.7926e+00
-9.1991e+00 2.7421e+00
-9.3594e+00 2.6846e+00
-9.4956e+00 2.6206e+00
-9.6085e+00 2.5512e+00
-9.7001e+00 2.4773e+00
-9.7729e+00 2.4002e+00
-9.8298e+00 2.3209e+00
-9.8736e+00 2.2404e+00
-9.9070e+00 2.1594e+00
-9.9535e+00 1.9839e+00
-9.9776e+00 1.8128e+00
-9.9896e+00 1.6476e+00
-9.9954e+00 1.4885e+00
-9.9981e+00 1.3354e+00
-9.9992e+00 1.1877e+00
-9.9997e+00 1.0447e+00
-9.9999e+00 9.0575e-01
-1.0000e+01 7.7027e-01
-1.0000e+01 6.3767e-01
-1.0000e+01 5.0741e-01
-1.0000e+01 3.7899e-01
-1.0000e+01 2.5192e-01
-1.0000e+01 1.2574e-01
-1.0000e+01 3.6739e-16
-1.0000e+01 -1.2574e-01
-1.0000e+01 -2.5192e-01
-1.0000e+01 -3.7899e-01
-1.0000e+01 -5.0741e-01
-1.0000e+01 -6.3767e-01
-1.0000e+01 -7.7027e-01
-9.9999e+00 -9.0575e-01
-9.9997e+00 -1.0447e+00
-9.9992e+00 -1.1877e+00
-9.9981e+00 -1.3354e+00
-9.9954e+00 -1.4885e+00
-9.9896e+00 -1.6476e+00
-9.9776e+00 -1.8128e+00
-9.9535e+00 -1.9839e+00
-9.9070e+00 -2.1594e+00
-9.8736e+00 -2.2404e+00
-9.8298e+00 -2.3209e+00
-9.7729e+00 -2.4002e+00
-9.7001e+00 -2.4773e+00
-9.6085e+00 -2.5512e+00
-9.4956e+00 -2.6206e+00
-9.3594e+00 -2.6846e+00
-9.1991e+00 -2.7421e+00
-9.0150e+00 -2.7926e+00
-8.8086e+00 -2.8359e+00
-8.5829e+00 -2.8720e+00
-8.3412e+00 -2.9015e+00
-8.0874e+00 -2.9252e+00
-7.8251e+00 -2.9437e+00
-7.5577e+00 -2.9581e+00
-7.2879e+00 -2.9691e+00
-7.0179e+00 -2.9773e+00
-6.7494e+00 -2.9835e+00
-6.4835e+00 -2.9881e+00
-6.2211e+00 -2.9915e+00
-5.9627e+00 -2.9940e+00
-5.7086e+00 -2.9957e+00
-5.4588e+00 -2.9970e+00
-5.2134e+00 -2.9979e+00
-4.9724e+00 -2.9986e+00
-4.7355e+00 -2.9991e+00
-4.5026e+00 -2.9994e+00
-4.2736e+00 -2.9996e+00
-4.0482e+00 -2.9997e+00
-3.8262e+00 -2.9998e+00
-3.6073e+00 -2.9999e+00
-3.3915e+00 -2.9999e+00
-3.1784e+00 -3.0000e+00
-2.9679e+00 -2.9831e+00
-2.7598e+00 -2.9190e+00
-2.5539e+00 -2.8112e+00
-2.3500e+00 -2.6655e+00
-2.1479e+00 -2.4888e+00
-1.9475e+00 -2.2888e+00
-1.7486e+00 -2.0740e+00
-1.5510e+00 -1.8527e+00
-1.3546e+00 -1.6337e+00
-1.1592e+00 -1.4249e+00
-9.6469e-01 -1.2341e+00
-7.7089e-01 -1.0680e+00
-5.7767e-01 -9.3254e-01
-3.8487e-01 -8.3230e-01
-1.9237e-01 -7.7075e-01
-1.8370e-15 -7.5000e-01
1.9237e-01 -7.7075e-01
3.8487e-01 -8.3230e-01
5.7767e-01 -9.3254e-01
7.7089e-01 -1.0680e+00
9.6469e-01 -1.2341e+00
1.1592e+00 -1.4249e+00
1.3546e+00 -1.6337e+00
1.5510e+00 -1.8527e+00
1.7486e+00 -2.0740e+00
1.9475e+00 -2.2888e+00
2.1479e+00 -2.4888e+00
2.3500e+00 -2.6655e+00
2.5539e+00 -2.8112e+00
2.7598e+00 -2.9190e+00
2.9679e+00 -2.9831e+00
3.1784e+00 -3.0000e+00
3.3915e+00 -2.9999e+00
3.6073e+00 -2.9999e+00
3.8262e+00 -2.9998e+00
4.0482e+00 -2.9997e+00
4.2736e+00 -2.9996e+00
4.5026e+00 -2.9994e+00
4.7355e+00 -2.9991e+00
4.9724e+00 -2.9986e+00
5.2134e+00 -2.9979e+00
5.4588e+00 -2.9970e+00
5.7086e+00 -2.9957e+00
5.9627e+00 -2.9940e+00
6.2211e+00 -2.9915e+00
6.4835e+00 -2.9881e+00
6.7494e+00 -2.9835e+00
7.0179e+00 -2.9773e+00
7.2879e+00 -2.9691e+00
7.5577e+00 -2.9581e+00
7.8251e+00 -2.9437e+00
8.0874e+00 -2.9252e+00
8.3412e+00 -2.9015e+00
8.5829e+00 -2.8720e+00
8.8086e+00 -2.8359e+00
9.0150e+00 -2.7926e+00
9.1991e+00 -2.7421e+00
9.3594e+00 -2.6846e+00
9.4956e+00 -2.6206e+00
9.6085e+00 -2.5512e+00
9.7001e+00 -2.4773e+00
9.7729e+00 -2.4002e+00
9.8298e+00 -2.3209e+00
9.8736e+00 -2.2404e+00
9.9070e+00 -2.1594e+00
9.9535e+00 -1.9839e+00
9.9776e+00 -1.8128e+00
9.9896e+00 -1.6476e+00
9.9954e+00 -1.4885e+00
9.9981e+00 -1.3354e+00
9.9992e+00 -1.1877e+00
9.9997e+00 -1.0447e+00
9.9999e+00 -9.0575e-01
1.0000e+01 -7.7027e-01
1.0000e+01 -6.3767e-01
1.0000e+01 -5.0741e-01
1.0000e+01 -3.7899e-01
1.0000e+01 -2.5192e-01
1.0000e+01 -1.2574e-01
1.0000e+01 0.0000e+00
};

\addplot [mark=none,red,line width=1.5] table{
-3.3854e+00 1.2849e+00
-3.4110e+00 1.2969e+00
-3.4375e+00 1.3083e+00
-3.4655e+00 1.3191e+00
-3.4953e+00 1.3290e+00
-3.5275e+00 1.3378e+00
-3.5623e+00 1.3450e+00
-3.5999e+00 1.3500e+00
-3.6404e+00 1.3518e+00
-3.6835e+00 1.3493e+00
-3.7287e+00 1.3413e+00
-3.7751e+00 1.3263e+00
-3.8212e+00 1.3031e+00
-3.8648e+00 1.2708e+00
-3.9037e+00 1.2290e+00
-3.9351e+00 1.1783e+00
-3.9569e+00 1.1199e+00
-3.9676e+00 1.0559e+00
-3.9665e+00 9.8868e-01
-3.9539e+00 9.2027e-01
-3.9310e+00 8.5229e-01
-3.8993e+00 7.8565e-01
-3.8606e+00 7.2061e-01
-3.8167e+00 6.5682e-01
-3.7694e+00 5.9355e-01
-3.7204e+00 5.2978e-01
-3.6715e+00 4.6437e-01
-3.6246e+00 3.9622e-01
-3.5817e+00 3.2442e-01
-3.5451e+00 2.4847e-01
-3.5168e+00 1.6846e-01
-3.4990e+00 8.5151e-02
-3.4929e+00 2.8084e-13
-3.4990e+00 -8.5151e-02
-3.5168e+00 -1.6846e-01
-3.5451e+00 -2.4847e-01
-3.5817e+00 -3.2442e-01
-3.6246e+00 -3.9622e-01
-3.6715e+00 -4.6437e-01
-3.7204e+00 -5.2978e-01
-3.7694e+00 -5.9355e-01
-3.8167e+00 -6.5682e-01
-3.8606e+00 -7.2061e-01
-3.8993e+00 -7.8565e-01
-3.9310e+00 -8.5229e-01
-3.9539e+00 -9.2027e-01
-3.9665e+00 -9.8868e-01
-3.9676e+00 -1.0559e+00
-3.9569e+00 -1.1199e+00
-3.9351e+00 -1.1783e+00
-3.9037e+00 -1.2290e+00
-3.8648e+00 -1.2708e+00
-3.8212e+00 -1.3031e+00
-3.7751e+00 -1.3263e+00
-3.7287e+00 -1.3413e+00
-3.6835e+00 -1.3493e+00
-3.6404e+00 -1.3518e+00
-3.5999e+00 -1.3500e+00
-3.5623e+00 -1.3450e+00
-3.5275e+00 -1.3378e+00
-3.4953e+00 -1.3290e+00
-3.4655e+00 -1.3191e+00
-3.4375e+00 -1.3083e+00
-3.4110e+00 -1.2969e+00
-3.3854e+00 -1.2849e+00
-3.3602e+00 -1.2720e+00
-3.3349e+00 -1.2581e+00
-3.3091e+00 -1.2431e+00
-3.2823e+00 -1.2267e+00
-3.2543e+00 -1.2086e+00
-3.2249e+00 -1.1888e+00
-3.1938e+00 -1.1671e+00
-3.1611e+00 -1.1433e+00
-3.1265e+00 -1.1174e+00
-3.0900e+00 -1.0895e+00
-3.0515e+00 -1.0596e+00
-3.0110e+00 -1.0276e+00
-2.9684e+00 -9.9382e-01
-2.9238e+00 -9.5824e-01
-2.8770e+00 -9.2103e-01
-2.8280e+00 -8.8234e-01
-2.7769e+00 -8.4231e-01
-2.7237e+00 -8.0107e-01
-2.6685e+00 -7.5874e-01
-2.6113e+00 -7.1540e-01
-2.5523e+00 -6.7107e-01
-2.4917e+00 -6.2568e-01
-2.4299e+00 -5.7907e-01
-2.3673e+00 -5.3093e-01
-2.3044e+00 -4.8072e-01
-2.2423e+00 -4.2770e-01
-2.1822e+00 -3.7084e-01
-2.1260e+00 -3.0893e-01
-2.0763e+00 -2.4080e-01
-2.0365e+00 -1.6587e-01
-2.0105e+00 -8.4798e-02
-2.0014e+00 6.6878e-13
-2.0105e+00 8.4798e-02
-2.0365e+00 1.6587e-01
-2.0763e+00 2.4080e-01
-2.1260e+00 3.0893e-01
-2.1822e+00 3.7084e-01
-2.2423e+00 4.2770e-01
-2.3044e+00 4.8072e-01
-2.3673e+00 5.3093e-01
-2.4299e+00 5.7907e-01
-2.4917e+00 6.2568e-01
-2.5523e+00 6.7107e-01
-2.6113e+00 7.1540e-01
-2.6685e+00 7.5874e-01
-2.7237e+00 8.0107e-01
-2.7769e+00 8.4231e-01
-2.8280e+00 8.8234e-01
-2.8770e+00 9.2103e-01
-2.9238e+00 9.5824e-01
-2.9684e+00 9.9382e-01
-3.0110e+00 1.0276e+00
-3.0515e+00 1.0596e+00
-3.0900e+00 1.0895e+00
-3.1265e+00 1.1174e+00
-3.1611e+00 1.1433e+00
-3.1938e+00 1.1671e+00
-3.2249e+00 1.1888e+00
-3.2543e+00 1.2086e+00
-3.2823e+00 1.2267e+00
-3.3091e+00 1.2431e+00
-3.3349e+00 1.2581e+00
-3.3602e+00 1.2720e+00
-3.3854e+00 1.2849e+00
};

\end{axis}

\end{tikzpicture} \\
\begin{tikzpicture}[scale=0.5]

\begin{axis}[
  xmin = -11,
  xmax = 11,
  ymin = -3.2,
  ymax = 3.2,
  scale only axis,
  axis equal image,
  hide axis,
  title = {\Huge$t=8$}
  ]

\addplot [mark=none,black,line width=1.5] table{
1.0000e+01 0.0000e+00
1.0000e+01 1.2574e-01
1.0000e+01 2.5192e-01
1.0000e+01 3.7899e-01
1.0000e+01 5.0741e-01
1.0000e+01 6.3767e-01
1.0000e+01 7.7027e-01
9.9999e+00 9.0575e-01
9.9997e+00 1.0447e+00
9.9992e+00 1.1877e+00
9.9981e+00 1.3354e+00
9.9954e+00 1.4885e+00
9.9896e+00 1.6476e+00
9.9776e+00 1.8128e+00
9.9535e+00 1.9839e+00
9.9070e+00 2.1594e+00
9.8736e+00 2.2404e+00
9.8298e+00 2.3209e+00
9.7729e+00 2.4002e+00
9.7001e+00 2.4773e+00
9.6085e+00 2.5512e+00
9.4956e+00 2.6206e+00
9.3594e+00 2.6846e+00
9.1991e+00 2.7421e+00
9.0150e+00 2.7926e+00
8.8086e+00 2.8359e+00
8.5829e+00 2.8720e+00
8.3412e+00 2.9015e+00
8.0874e+00 2.9252e+00
7.8251e+00 2.9437e+00
7.5577e+00 2.9581e+00
7.2879e+00 2.9691e+00
7.0179e+00 2.9773e+00
6.7494e+00 2.9835e+00
6.4835e+00 2.9881e+00
6.2211e+00 2.9915e+00
5.9627e+00 2.9940e+00
5.7086e+00 2.9957e+00
5.4588e+00 2.9970e+00
5.2134e+00 2.9979e+00
4.9724e+00 2.9986e+00
4.7355e+00 2.9991e+00
4.5026e+00 2.9994e+00
4.2736e+00 2.9996e+00
4.0482e+00 2.9997e+00
3.8262e+00 2.9998e+00
3.6073e+00 2.9999e+00
3.3915e+00 2.9999e+00
3.1784e+00 3.0000e+00
2.9679e+00 2.9831e+00
2.7598e+00 2.9190e+00
2.5539e+00 2.8112e+00
2.3500e+00 2.6655e+00
2.1479e+00 2.4888e+00
1.9475e+00 2.2888e+00
1.7486e+00 2.0740e+00
1.5510e+00 1.8527e+00
1.3546e+00 1.6337e+00
1.1592e+00 1.4249e+00
9.6469e-01 1.2341e+00
7.7089e-01 1.0680e+00
5.7767e-01 9.3254e-01
3.8487e-01 8.3230e-01
1.9237e-01 7.7075e-01
6.1232e-16 7.5000e-01
-1.9237e-01 7.7075e-01
-3.8487e-01 8.3230e-01
-5.7767e-01 9.3254e-01
-7.7089e-01 1.0680e+00
-9.6469e-01 1.2341e+00
-1.1592e+00 1.4249e+00
-1.3546e+00 1.6337e+00
-1.5510e+00 1.8527e+00
-1.7486e+00 2.0740e+00
-1.9475e+00 2.2888e+00
-2.1479e+00 2.4888e+00
-2.3500e+00 2.6655e+00
-2.5539e+00 2.8112e+00
-2.7598e+00 2.9190e+00
-2.9679e+00 2.9831e+00
-3.1784e+00 3.0000e+00
-3.3915e+00 2.9999e+00
-3.6073e+00 2.9999e+00
-3.8262e+00 2.9998e+00
-4.0482e+00 2.9997e+00
-4.2736e+00 2.9996e+00
-4.5026e+00 2.9994e+00
-4.7355e+00 2.9991e+00
-4.9724e+00 2.9986e+00
-5.2134e+00 2.9979e+00
-5.4588e+00 2.9970e+00
-5.7086e+00 2.9957e+00
-5.9627e+00 2.9940e+00
-6.2211e+00 2.9915e+00
-6.4835e+00 2.9881e+00
-6.7494e+00 2.9835e+00
-7.0179e+00 2.9773e+00
-7.2879e+00 2.9691e+00
-7.5577e+00 2.9581e+00
-7.8251e+00 2.9437e+00
-8.0874e+00 2.9252e+00
-8.3412e+00 2.9015e+00
-8.5829e+00 2.8720e+00
-8.8086e+00 2.8359e+00
-9.0150e+00 2.7926e+00
-9.1991e+00 2.7421e+00
-9.3594e+00 2.6846e+00
-9.4956e+00 2.6206e+00
-9.6085e+00 2.5512e+00
-9.7001e+00 2.4773e+00
-9.7729e+00 2.4002e+00
-9.8298e+00 2.3209e+00
-9.8736e+00 2.2404e+00
-9.9070e+00 2.1594e+00
-9.9535e+00 1.9839e+00
-9.9776e+00 1.8128e+00
-9.9896e+00 1.6476e+00
-9.9954e+00 1.4885e+00
-9.9981e+00 1.3354e+00
-9.9992e+00 1.1877e+00
-9.9997e+00 1.0447e+00
-9.9999e+00 9.0575e-01
-1.0000e+01 7.7027e-01
-1.0000e+01 6.3767e-01
-1.0000e+01 5.0741e-01
-1.0000e+01 3.7899e-01
-1.0000e+01 2.5192e-01
-1.0000e+01 1.2574e-01
-1.0000e+01 3.6739e-16
-1.0000e+01 -1.2574e-01
-1.0000e+01 -2.5192e-01
-1.0000e+01 -3.7899e-01
-1.0000e+01 -5.0741e-01
-1.0000e+01 -6.3767e-01
-1.0000e+01 -7.7027e-01
-9.9999e+00 -9.0575e-01
-9.9997e+00 -1.0447e+00
-9.9992e+00 -1.1877e+00
-9.9981e+00 -1.3354e+00
-9.9954e+00 -1.4885e+00
-9.9896e+00 -1.6476e+00
-9.9776e+00 -1.8128e+00
-9.9535e+00 -1.9839e+00
-9.9070e+00 -2.1594e+00
-9.8736e+00 -2.2404e+00
-9.8298e+00 -2.3209e+00
-9.7729e+00 -2.4002e+00
-9.7001e+00 -2.4773e+00
-9.6085e+00 -2.5512e+00
-9.4956e+00 -2.6206e+00
-9.3594e+00 -2.6846e+00
-9.1991e+00 -2.7421e+00
-9.0150e+00 -2.7926e+00
-8.8086e+00 -2.8359e+00
-8.5829e+00 -2.8720e+00
-8.3412e+00 -2.9015e+00
-8.0874e+00 -2.9252e+00
-7.8251e+00 -2.9437e+00
-7.5577e+00 -2.9581e+00
-7.2879e+00 -2.9691e+00
-7.0179e+00 -2.9773e+00
-6.7494e+00 -2.9835e+00
-6.4835e+00 -2.9881e+00
-6.2211e+00 -2.9915e+00
-5.9627e+00 -2.9940e+00
-5.7086e+00 -2.9957e+00
-5.4588e+00 -2.9970e+00
-5.2134e+00 -2.9979e+00
-4.9724e+00 -2.9986e+00
-4.7355e+00 -2.9991e+00
-4.5026e+00 -2.9994e+00
-4.2736e+00 -2.9996e+00
-4.0482e+00 -2.9997e+00
-3.8262e+00 -2.9998e+00
-3.6073e+00 -2.9999e+00
-3.3915e+00 -2.9999e+00
-3.1784e+00 -3.0000e+00
-2.9679e+00 -2.9831e+00
-2.7598e+00 -2.9190e+00
-2.5539e+00 -2.8112e+00
-2.3500e+00 -2.6655e+00
-2.1479e+00 -2.4888e+00
-1.9475e+00 -2.2888e+00
-1.7486e+00 -2.0740e+00
-1.5510e+00 -1.8527e+00
-1.3546e+00 -1.6337e+00
-1.1592e+00 -1.4249e+00
-9.6469e-01 -1.2341e+00
-7.7089e-01 -1.0680e+00
-5.7767e-01 -9.3254e-01
-3.8487e-01 -8.3230e-01
-1.9237e-01 -7.7075e-01
-1.8370e-15 -7.5000e-01
1.9237e-01 -7.7075e-01
3.8487e-01 -8.3230e-01
5.7767e-01 -9.3254e-01
7.7089e-01 -1.0680e+00
9.6469e-01 -1.2341e+00
1.1592e+00 -1.4249e+00
1.3546e+00 -1.6337e+00
1.5510e+00 -1.8527e+00
1.7486e+00 -2.0740e+00
1.9475e+00 -2.2888e+00
2.1479e+00 -2.4888e+00
2.3500e+00 -2.6655e+00
2.5539e+00 -2.8112e+00
2.7598e+00 -2.9190e+00
2.9679e+00 -2.9831e+00
3.1784e+00 -3.0000e+00
3.3915e+00 -2.9999e+00
3.6073e+00 -2.9999e+00
3.8262e+00 -2.9998e+00
4.0482e+00 -2.9997e+00
4.2736e+00 -2.9996e+00
4.5026e+00 -2.9994e+00
4.7355e+00 -2.9991e+00
4.9724e+00 -2.9986e+00
5.2134e+00 -2.9979e+00
5.4588e+00 -2.9970e+00
5.7086e+00 -2.9957e+00
5.9627e+00 -2.9940e+00
6.2211e+00 -2.9915e+00
6.4835e+00 -2.9881e+00
6.7494e+00 -2.9835e+00
7.0179e+00 -2.9773e+00
7.2879e+00 -2.9691e+00
7.5577e+00 -2.9581e+00
7.8251e+00 -2.9437e+00
8.0874e+00 -2.9252e+00
8.3412e+00 -2.9015e+00
8.5829e+00 -2.8720e+00
8.8086e+00 -2.8359e+00
9.0150e+00 -2.7926e+00
9.1991e+00 -2.7421e+00
9.3594e+00 -2.6846e+00
9.4956e+00 -2.6206e+00
9.6085e+00 -2.5512e+00
9.7001e+00 -2.4773e+00
9.7729e+00 -2.4002e+00
9.8298e+00 -2.3209e+00
9.8736e+00 -2.2404e+00
9.9070e+00 -2.1594e+00
9.9535e+00 -1.9839e+00
9.9776e+00 -1.8128e+00
9.9896e+00 -1.6476e+00
9.9954e+00 -1.4885e+00
9.9981e+00 -1.3354e+00
9.9992e+00 -1.1877e+00
9.9997e+00 -1.0447e+00
9.9999e+00 -9.0575e-01
1.0000e+01 -7.7027e-01
1.0000e+01 -6.3767e-01
1.0000e+01 -5.0741e-01
1.0000e+01 -3.7899e-01
1.0000e+01 -2.5192e-01
1.0000e+01 -1.2574e-01
1.0000e+01 0.0000e+00
};

\addplot [mark=none,red,line width=1.5] table{
-2.1130e+00 8.4143e-01
-2.1384e+00 8.5395e-01
-2.1643e+00 8.6672e-01
-2.1911e+00 8.7993e-01
-2.2194e+00 8.9374e-01
-2.2494e+00 9.0826e-01
-2.2816e+00 9.2352e-01
-2.3160e+00 9.3950e-01
-2.3530e+00 9.5606e-01
-2.3928e+00 9.7297e-01
-2.4357e+00 9.8985e-01
-2.4818e+00 1.0061e+00
-2.5314e+00 1.0210e+00
-2.5846e+00 1.0332e+00
-2.6414e+00 1.0413e+00
-2.7015e+00 1.0433e+00
-2.7639e+00 1.0366e+00
-2.8267e+00 1.0185e+00
-2.8865e+00 9.8675e-01
-2.9388e+00 9.4031e-01
-2.9786e+00 8.8026e-01
-3.0017e+00 8.1000e-01
-3.0062e+00 7.3438e-01
-2.9930e+00 6.5808e-01
-2.9648e+00 5.8429e-01
-2.9254e+00 5.1419e-01
-2.8788e+00 4.4712e-01
-2.8288e+00 3.8119e-01
-2.7791e+00 3.1389e-01
-2.7337e+00 2.4283e-01
-2.6967e+00 1.6647e-01
-2.6724e+00 8.4875e-02
-2.6638e+00 1.7263e-13
-2.6724e+00 -8.4875e-02
-2.6967e+00 -1.6647e-01
-2.7337e+00 -2.4283e-01
-2.7791e+00 -3.1389e-01
-2.8288e+00 -3.8119e-01
-2.8788e+00 -4.4712e-01
-2.9254e+00 -5.1419e-01
-2.9648e+00 -5.8429e-01
-2.9930e+00 -6.5808e-01
-3.0062e+00 -7.3438e-01
-3.0017e+00 -8.1000e-01
-2.9786e+00 -8.8026e-01
-2.9388e+00 -9.4031e-01
-2.8865e+00 -9.8675e-01
-2.8267e+00 -1.0185e+00
-2.7639e+00 -1.0366e+00
-2.7015e+00 -1.0433e+00
-2.6414e+00 -1.0413e+00
-2.5846e+00 -1.0332e+00
-2.5314e+00 -1.0210e+00
-2.4818e+00 -1.0061e+00
-2.4357e+00 -9.8985e-01
-2.3928e+00 -9.7297e-01
-2.3530e+00 -9.5606e-01
-2.3160e+00 -9.3950e-01
-2.2816e+00 -9.2352e-01
-2.2494e+00 -9.0826e-01
-2.2194e+00 -8.9374e-01
-2.1911e+00 -8.7993e-01
-2.1643e+00 -8.6672e-01
-2.1384e+00 -8.5395e-01
-2.1130e+00 -8.4143e-01
-2.0876e+00 -8.2893e-01
-2.0617e+00 -8.1622e-01
-2.0349e+00 -8.0311e-01
-2.0066e+00 -7.8941e-01
-1.9765e+00 -7.7499e-01
-1.9445e+00 -7.5977e-01
-1.9101e+00 -7.4370e-01
-1.8734e+00 -7.2675e-01
-1.8340e+00 -7.0895e-01
-1.7920e+00 -6.9032e-01
-1.7472e+00 -6.7093e-01
-1.6997e+00 -6.5086e-01
-1.6494e+00 -6.3019e-01
-1.5963e+00 -6.0901e-01
-1.5405e+00 -5.8741e-01
-1.4821e+00 -5.6548e-01
-1.4210e+00 -5.4330e-01
-1.3575e+00 -5.2095e-01
-1.2916e+00 -4.9848e-01
-1.2234e+00 -4.7589e-01
-1.1532e+00 -4.5318e-01
-1.0810e+00 -4.3027e-01
-1.0071e+00 -4.0696e-01
-9.3181e-01 -3.8293e-01
-8.5541e-01 -3.5764e-01
-7.7843e-01 -3.3015e-01
-7.0171e-01 -2.9899e-01
-6.2670e-01 -2.6187e-01
-5.5616e-01 -2.1561e-01
-4.9522e-01 -1.5668e-01
-4.5214e-01 -8.3469e-02
-4.3629e-01 1.0475e-12
-4.5214e-01 8.3469e-02
-4.9522e-01 1.5668e-01
-5.5616e-01 2.1561e-01
-6.2670e-01 2.6187e-01
-7.0171e-01 2.9899e-01
-7.7843e-01 3.3015e-01
-8.5541e-01 3.5764e-01
-9.3181e-01 3.8293e-01
-1.0071e+00 4.0696e-01
-1.0810e+00 4.3027e-01
-1.1532e+00 4.5318e-01
-1.2234e+00 4.7589e-01
-1.2916e+00 4.9848e-01
-1.3575e+00 5.2095e-01
-1.4210e+00 5.4330e-01
-1.4821e+00 5.6548e-01
-1.5405e+00 5.8741e-01
-1.5963e+00 6.0901e-01
-1.6494e+00 6.3019e-01
-1.6997e+00 6.5086e-01
-1.7472e+00 6.7093e-01
-1.7920e+00 6.9032e-01
-1.8340e+00 7.0895e-01
-1.8734e+00 7.2675e-01
-1.9101e+00 7.4370e-01
-1.9445e+00 7.5977e-01
-1.9765e+00 7.7499e-01
-2.0066e+00 7.8941e-01
-2.0349e+00 8.0311e-01
-2.0617e+00 8.1622e-01
-2.0876e+00 8.2893e-01
-2.1130e+00 8.4143e-01
};

\end{axis}

\end{tikzpicture} &
\begin{tikzpicture}[scale=0.5]

\begin{axis}[
  xmin = -11,
  xmax = 11,
  ymin = -3.2,
  ymax = 3.2,
  scale only axis,
  axis equal image,
  hide axis,
  title = {\Huge$t=9$}
  ]

\addplot [mark=none,black,line width=1.5] table{
1.0000e+01 0.0000e+00
1.0000e+01 1.2574e-01
1.0000e+01 2.5192e-01
1.0000e+01 3.7899e-01
1.0000e+01 5.0741e-01
1.0000e+01 6.3767e-01
1.0000e+01 7.7027e-01
9.9999e+00 9.0575e-01
9.9997e+00 1.0447e+00
9.9992e+00 1.1877e+00
9.9981e+00 1.3354e+00
9.9954e+00 1.4885e+00
9.9896e+00 1.6476e+00
9.9776e+00 1.8128e+00
9.9535e+00 1.9839e+00
9.9070e+00 2.1594e+00
9.8736e+00 2.2404e+00
9.8298e+00 2.3209e+00
9.7729e+00 2.4002e+00
9.7001e+00 2.4773e+00
9.6085e+00 2.5512e+00
9.4956e+00 2.6206e+00
9.3594e+00 2.6846e+00
9.1991e+00 2.7421e+00
9.0150e+00 2.7926e+00
8.8086e+00 2.8359e+00
8.5829e+00 2.8720e+00
8.3412e+00 2.9015e+00
8.0874e+00 2.9252e+00
7.8251e+00 2.9437e+00
7.5577e+00 2.9581e+00
7.2879e+00 2.9691e+00
7.0179e+00 2.9773e+00
6.7494e+00 2.9835e+00
6.4835e+00 2.9881e+00
6.2211e+00 2.9915e+00
5.9627e+00 2.9940e+00
5.7086e+00 2.9957e+00
5.4588e+00 2.9970e+00
5.2134e+00 2.9979e+00
4.9724e+00 2.9986e+00
4.7355e+00 2.9991e+00
4.5026e+00 2.9994e+00
4.2736e+00 2.9996e+00
4.0482e+00 2.9997e+00
3.8262e+00 2.9998e+00
3.6073e+00 2.9999e+00
3.3915e+00 2.9999e+00
3.1784e+00 3.0000e+00
2.9679e+00 2.9831e+00
2.7598e+00 2.9190e+00
2.5539e+00 2.8112e+00
2.3500e+00 2.6655e+00
2.1479e+00 2.4888e+00
1.9475e+00 2.2888e+00
1.7486e+00 2.0740e+00
1.5510e+00 1.8527e+00
1.3546e+00 1.6337e+00
1.1592e+00 1.4249e+00
9.6469e-01 1.2341e+00
7.7089e-01 1.0680e+00
5.7767e-01 9.3254e-01
3.8487e-01 8.3230e-01
1.9237e-01 7.7075e-01
6.1232e-16 7.5000e-01
-1.9237e-01 7.7075e-01
-3.8487e-01 8.3230e-01
-5.7767e-01 9.3254e-01
-7.7089e-01 1.0680e+00
-9.6469e-01 1.2341e+00
-1.1592e+00 1.4249e+00
-1.3546e+00 1.6337e+00
-1.5510e+00 1.8527e+00
-1.7486e+00 2.0740e+00
-1.9475e+00 2.2888e+00
-2.1479e+00 2.4888e+00
-2.3500e+00 2.6655e+00
-2.5539e+00 2.8112e+00
-2.7598e+00 2.9190e+00
-2.9679e+00 2.9831e+00
-3.1784e+00 3.0000e+00
-3.3915e+00 2.9999e+00
-3.6073e+00 2.9999e+00
-3.8262e+00 2.9998e+00
-4.0482e+00 2.9997e+00
-4.2736e+00 2.9996e+00
-4.5026e+00 2.9994e+00
-4.7355e+00 2.9991e+00
-4.9724e+00 2.9986e+00
-5.2134e+00 2.9979e+00
-5.4588e+00 2.9970e+00
-5.7086e+00 2.9957e+00
-5.9627e+00 2.9940e+00
-6.2211e+00 2.9915e+00
-6.4835e+00 2.9881e+00
-6.7494e+00 2.9835e+00
-7.0179e+00 2.9773e+00
-7.2879e+00 2.9691e+00
-7.5577e+00 2.9581e+00
-7.8251e+00 2.9437e+00
-8.0874e+00 2.9252e+00
-8.3412e+00 2.9015e+00
-8.5829e+00 2.8720e+00
-8.8086e+00 2.8359e+00
-9.0150e+00 2.7926e+00
-9.1991e+00 2.7421e+00
-9.3594e+00 2.6846e+00
-9.4956e+00 2.6206e+00
-9.6085e+00 2.5512e+00
-9.7001e+00 2.4773e+00
-9.7729e+00 2.4002e+00
-9.8298e+00 2.3209e+00
-9.8736e+00 2.2404e+00
-9.9070e+00 2.1594e+00
-9.9535e+00 1.9839e+00
-9.9776e+00 1.8128e+00
-9.9896e+00 1.6476e+00
-9.9954e+00 1.4885e+00
-9.9981e+00 1.3354e+00
-9.9992e+00 1.1877e+00
-9.9997e+00 1.0447e+00
-9.9999e+00 9.0575e-01
-1.0000e+01 7.7027e-01
-1.0000e+01 6.3767e-01
-1.0000e+01 5.0741e-01
-1.0000e+01 3.7899e-01
-1.0000e+01 2.5192e-01
-1.0000e+01 1.2574e-01
-1.0000e+01 3.6739e-16
-1.0000e+01 -1.2574e-01
-1.0000e+01 -2.5192e-01
-1.0000e+01 -3.7899e-01
-1.0000e+01 -5.0741e-01
-1.0000e+01 -6.3767e-01
-1.0000e+01 -7.7027e-01
-9.9999e+00 -9.0575e-01
-9.9997e+00 -1.0447e+00
-9.9992e+00 -1.1877e+00
-9.9981e+00 -1.3354e+00
-9.9954e+00 -1.4885e+00
-9.9896e+00 -1.6476e+00
-9.9776e+00 -1.8128e+00
-9.9535e+00 -1.9839e+00
-9.9070e+00 -2.1594e+00
-9.8736e+00 -2.2404e+00
-9.8298e+00 -2.3209e+00
-9.7729e+00 -2.4002e+00
-9.7001e+00 -2.4773e+00
-9.6085e+00 -2.5512e+00
-9.4956e+00 -2.6206e+00
-9.3594e+00 -2.6846e+00
-9.1991e+00 -2.7421e+00
-9.0150e+00 -2.7926e+00
-8.8086e+00 -2.8359e+00
-8.5829e+00 -2.8720e+00
-8.3412e+00 -2.9015e+00
-8.0874e+00 -2.9252e+00
-7.8251e+00 -2.9437e+00
-7.5577e+00 -2.9581e+00
-7.2879e+00 -2.9691e+00
-7.0179e+00 -2.9773e+00
-6.7494e+00 -2.9835e+00
-6.4835e+00 -2.9881e+00
-6.2211e+00 -2.9915e+00
-5.9627e+00 -2.9940e+00
-5.7086e+00 -2.9957e+00
-5.4588e+00 -2.9970e+00
-5.2134e+00 -2.9979e+00
-4.9724e+00 -2.9986e+00
-4.7355e+00 -2.9991e+00
-4.5026e+00 -2.9994e+00
-4.2736e+00 -2.9996e+00
-4.0482e+00 -2.9997e+00
-3.8262e+00 -2.9998e+00
-3.6073e+00 -2.9999e+00
-3.3915e+00 -2.9999e+00
-3.1784e+00 -3.0000e+00
-2.9679e+00 -2.9831e+00
-2.7598e+00 -2.9190e+00
-2.5539e+00 -2.8112e+00
-2.3500e+00 -2.6655e+00
-2.1479e+00 -2.4888e+00
-1.9475e+00 -2.2888e+00
-1.7486e+00 -2.0740e+00
-1.5510e+00 -1.8527e+00
-1.3546e+00 -1.6337e+00
-1.1592e+00 -1.4249e+00
-9.6469e-01 -1.2341e+00
-7.7089e-01 -1.0680e+00
-5.7767e-01 -9.3254e-01
-3.8487e-01 -8.3230e-01
-1.9237e-01 -7.7075e-01
-1.8370e-15 -7.5000e-01
1.9237e-01 -7.7075e-01
3.8487e-01 -8.3230e-01
5.7767e-01 -9.3254e-01
7.7089e-01 -1.0680e+00
9.6469e-01 -1.2341e+00
1.1592e+00 -1.4249e+00
1.3546e+00 -1.6337e+00
1.5510e+00 -1.8527e+00
1.7486e+00 -2.0740e+00
1.9475e+00 -2.2888e+00
2.1479e+00 -2.4888e+00
2.3500e+00 -2.6655e+00
2.5539e+00 -2.8112e+00
2.7598e+00 -2.9190e+00
2.9679e+00 -2.9831e+00
3.1784e+00 -3.0000e+00
3.3915e+00 -2.9999e+00
3.6073e+00 -2.9999e+00
3.8262e+00 -2.9998e+00
4.0482e+00 -2.9997e+00
4.2736e+00 -2.9996e+00
4.5026e+00 -2.9994e+00
4.7355e+00 -2.9991e+00
4.9724e+00 -2.9986e+00
5.2134e+00 -2.9979e+00
5.4588e+00 -2.9970e+00
5.7086e+00 -2.9957e+00
5.9627e+00 -2.9940e+00
6.2211e+00 -2.9915e+00
6.4835e+00 -2.9881e+00
6.7494e+00 -2.9835e+00
7.0179e+00 -2.9773e+00
7.2879e+00 -2.9691e+00
7.5577e+00 -2.9581e+00
7.8251e+00 -2.9437e+00
8.0874e+00 -2.9252e+00
8.3412e+00 -2.9015e+00
8.5829e+00 -2.8720e+00
8.8086e+00 -2.8359e+00
9.0150e+00 -2.7926e+00
9.1991e+00 -2.7421e+00
9.3594e+00 -2.6846e+00
9.4956e+00 -2.6206e+00
9.6085e+00 -2.5512e+00
9.7001e+00 -2.4773e+00
9.7729e+00 -2.4002e+00
9.8298e+00 -2.3209e+00
9.8736e+00 -2.2404e+00
9.9070e+00 -2.1594e+00
9.9535e+00 -1.9839e+00
9.9776e+00 -1.8128e+00
9.9896e+00 -1.6476e+00
9.9954e+00 -1.4885e+00
9.9981e+00 -1.3354e+00
9.9992e+00 -1.1877e+00
9.9997e+00 -1.0447e+00
9.9999e+00 -9.0575e-01
1.0000e+01 -7.7027e-01
1.0000e+01 -6.3767e-01
1.0000e+01 -5.0741e-01
1.0000e+01 -3.7899e-01
1.0000e+01 -2.5192e-01
1.0000e+01 -1.2574e-01
1.0000e+01 0.0000e+00
};

\addplot [mark=none,red,line width=1.5] table{
1.5293e-01 4.4305e-01
1.2485e-01 4.3899e-01
9.6140e-02 4.3546e-01
6.6288e-02 4.3245e-01
3.4862e-02 4.2995e-01
1.5008e-03 4.2801e-01
-3.4080e-02 4.2669e-01
-7.2096e-02 4.2608e-01
-1.1270e-01 4.2626e-01
-1.5600e-01 4.2732e-01
-2.0206e-01 4.2931e-01
-2.5091e-01 4.3229e-01
-3.0255e-01 4.3626e-01
-3.5698e-01 4.4121e-01
-4.1416e-01 4.4705e-01
-4.7408e-01 4.5364e-01
-5.3670e-01 4.6074e-01
-6.0200e-01 4.6801e-01
-6.6995e-01 4.7498e-01
-7.4047e-01 4.8101e-01
-8.1348e-01 4.8530e-01
-8.8879e-01 4.8683e-01
-9.6608e-01 4.8441e-01
-1.0448e+00 4.7664e-01
-1.1242e+00 4.6200e-01
-1.2031e+00 4.3892e-01
-1.2796e+00 4.0598e-01
-1.3516e+00 3.6212e-01
-1.4163e+00 3.0701e-01
-1.4705e+00 2.4121e-01
-1.5116e+00 1.6634e-01
-1.5373e+00 8.4889e-02
-1.5460e+00 -1.1010e-12
-1.5373e+00 -8.4889e-02
-1.5116e+00 -1.6634e-01
-1.4705e+00 -2.4121e-01
-1.4163e+00 -3.0701e-01
-1.3516e+00 -3.6212e-01
-1.2796e+00 -4.0598e-01
-1.2031e+00 -4.3892e-01
-1.1242e+00 -4.6200e-01
-1.0448e+00 -4.7664e-01
-9.6608e-01 -4.8441e-01
-8.8879e-01 -4.8683e-01
-8.1348e-01 -4.8530e-01
-7.4047e-01 -4.8101e-01
-6.6995e-01 -4.7498e-01
-6.0200e-01 -4.6801e-01
-5.3670e-01 -4.6074e-01
-4.7408e-01 -4.5364e-01
-4.1416e-01 -4.4705e-01
-3.5698e-01 -4.4121e-01
-3.0255e-01 -4.3626e-01
-2.5091e-01 -4.3229e-01
-2.0206e-01 -4.2931e-01
-1.5600e-01 -4.2732e-01
-1.1270e-01 -4.2626e-01
-7.2096e-02 -4.2608e-01
-3.4080e-02 -4.2669e-01
1.5008e-03 -4.2801e-01
3.4862e-02 -4.2995e-01
6.6288e-02 -4.3245e-01
9.6140e-02 -4.3546e-01
1.2485e-01 -4.3899e-01
1.5293e-01 -4.4305e-01
1.8090e-01 -4.4775e-01
2.0931e-01 -4.5318e-01
2.3862e-01 -4.5953e-01
2.6925e-01 -4.6698e-01
3.0150e-01 -4.7573e-01
3.3561e-01 -4.8599e-01
3.7173e-01 -4.9795e-01
4.0997e-01 -5.1176e-01
4.5039e-01 -5.2747e-01
4.9308e-01 -5.4504e-01
5.3813e-01 -5.6425e-01
5.8571e-01 -5.8469e-01
6.3605e-01 -6.0574e-01
6.8945e-01 -6.2655e-01
7.4617e-01 -6.4605e-01
8.0641e-01 -6.6300e-01
8.7017e-01 -6.7606e-01
9.3720e-01 -6.8392e-01
1.0070e+00 -6.8536e-01
1.0787e+00 -6.7936e-01
1.1513e+00 -6.6515e-01
1.2237e+00 -6.4226e-01
1.2946e+00 -6.1045e-01
1.3626e+00 -5.6978e-01
1.4265e+00 -5.2050e-01
1.4850e+00 -4.6310e-01
1.5368e+00 -3.9828e-01
1.5809e+00 -3.2692e-01
1.6163e+00 -2.5010e-01
1.6422e+00 -1.6908e-01
1.6579e+00 -8.5241e-02
1.6632e+00 3.5240e-13
1.6579e+00 8.5241e-02
1.6422e+00 1.6908e-01
1.6163e+00 2.5010e-01
1.5809e+00 3.2692e-01
1.5368e+00 3.9828e-01
1.4850e+00 4.6310e-01
1.4265e+00 5.2050e-01
1.3626e+00 5.6978e-01
1.2946e+00 6.1045e-01
1.2237e+00 6.4226e-01
1.1513e+00 6.6515e-01
1.0787e+00 6.7936e-01
1.0070e+00 6.8536e-01
9.3720e-01 6.8392e-01
8.7017e-01 6.7606e-01
8.0641e-01 6.6300e-01
7.4617e-01 6.4605e-01
6.8945e-01 6.2655e-01
6.3605e-01 6.0574e-01
5.8571e-01 5.8469e-01
5.3813e-01 5.6425e-01
4.9308e-01 5.4504e-01
4.5039e-01 5.2747e-01
4.0997e-01 5.1176e-01
3.7173e-01 4.9795e-01
3.3561e-01 4.8599e-01
3.0150e-01 4.7573e-01
2.6925e-01 4.6698e-01
2.3862e-01 4.5953e-01
2.0931e-01 4.5318e-01
1.8090e-01 4.4775e-01
1.5293e-01 4.4305e-01
};

\end{axis}

\end{tikzpicture} &
\begin{tikzpicture}[scale=0.5]

\begin{axis}[
  xmin = -11,
  xmax = 11,
  ymin = -3.2,
  ymax = 3.2,
  scale only axis,
  axis equal image,
  hide axis,
  title = {\Huge$t=10$}
  ]

\addplot [mark=none,black,line width=1.5] table{
1.0000e+01 0.0000e+00
1.0000e+01 1.2574e-01
1.0000e+01 2.5192e-01
1.0000e+01 3.7899e-01
1.0000e+01 5.0741e-01
1.0000e+01 6.3767e-01
1.0000e+01 7.7027e-01
9.9999e+00 9.0575e-01
9.9997e+00 1.0447e+00
9.9992e+00 1.1877e+00
9.9981e+00 1.3354e+00
9.9954e+00 1.4885e+00
9.9896e+00 1.6476e+00
9.9776e+00 1.8128e+00
9.9535e+00 1.9839e+00
9.9070e+00 2.1594e+00
9.8736e+00 2.2404e+00
9.8298e+00 2.3209e+00
9.7729e+00 2.4002e+00
9.7001e+00 2.4773e+00
9.6085e+00 2.5512e+00
9.4956e+00 2.6206e+00
9.3594e+00 2.6846e+00
9.1991e+00 2.7421e+00
9.0150e+00 2.7926e+00
8.8086e+00 2.8359e+00
8.5829e+00 2.8720e+00
8.3412e+00 2.9015e+00
8.0874e+00 2.9252e+00
7.8251e+00 2.9437e+00
7.5577e+00 2.9581e+00
7.2879e+00 2.9691e+00
7.0179e+00 2.9773e+00
6.7494e+00 2.9835e+00
6.4835e+00 2.9881e+00
6.2211e+00 2.9915e+00
5.9627e+00 2.9940e+00
5.7086e+00 2.9957e+00
5.4588e+00 2.9970e+00
5.2134e+00 2.9979e+00
4.9724e+00 2.9986e+00
4.7355e+00 2.9991e+00
4.5026e+00 2.9994e+00
4.2736e+00 2.9996e+00
4.0482e+00 2.9997e+00
3.8262e+00 2.9998e+00
3.6073e+00 2.9999e+00
3.3915e+00 2.9999e+00
3.1784e+00 3.0000e+00
2.9679e+00 2.9831e+00
2.7598e+00 2.9190e+00
2.5539e+00 2.8112e+00
2.3500e+00 2.6655e+00
2.1479e+00 2.4888e+00
1.9475e+00 2.2888e+00
1.7486e+00 2.0740e+00
1.5510e+00 1.8527e+00
1.3546e+00 1.6337e+00
1.1592e+00 1.4249e+00
9.6469e-01 1.2341e+00
7.7089e-01 1.0680e+00
5.7767e-01 9.3254e-01
3.8487e-01 8.3230e-01
1.9237e-01 7.7075e-01
6.1232e-16 7.5000e-01
-1.9237e-01 7.7075e-01
-3.8487e-01 8.3230e-01
-5.7767e-01 9.3254e-01
-7.7089e-01 1.0680e+00
-9.6469e-01 1.2341e+00
-1.1592e+00 1.4249e+00
-1.3546e+00 1.6337e+00
-1.5510e+00 1.8527e+00
-1.7486e+00 2.0740e+00
-1.9475e+00 2.2888e+00
-2.1479e+00 2.4888e+00
-2.3500e+00 2.6655e+00
-2.5539e+00 2.8112e+00
-2.7598e+00 2.9190e+00
-2.9679e+00 2.9831e+00
-3.1784e+00 3.0000e+00
-3.3915e+00 2.9999e+00
-3.6073e+00 2.9999e+00
-3.8262e+00 2.9998e+00
-4.0482e+00 2.9997e+00
-4.2736e+00 2.9996e+00
-4.5026e+00 2.9994e+00
-4.7355e+00 2.9991e+00
-4.9724e+00 2.9986e+00
-5.2134e+00 2.9979e+00
-5.4588e+00 2.9970e+00
-5.7086e+00 2.9957e+00
-5.9627e+00 2.9940e+00
-6.2211e+00 2.9915e+00
-6.4835e+00 2.9881e+00
-6.7494e+00 2.9835e+00
-7.0179e+00 2.9773e+00
-7.2879e+00 2.9691e+00
-7.5577e+00 2.9581e+00
-7.8251e+00 2.9437e+00
-8.0874e+00 2.9252e+00
-8.3412e+00 2.9015e+00
-8.5829e+00 2.8720e+00
-8.8086e+00 2.8359e+00
-9.0150e+00 2.7926e+00
-9.1991e+00 2.7421e+00
-9.3594e+00 2.6846e+00
-9.4956e+00 2.6206e+00
-9.6085e+00 2.5512e+00
-9.7001e+00 2.4773e+00
-9.7729e+00 2.4002e+00
-9.8298e+00 2.3209e+00
-9.8736e+00 2.2404e+00
-9.9070e+00 2.1594e+00
-9.9535e+00 1.9839e+00
-9.9776e+00 1.8128e+00
-9.9896e+00 1.6476e+00
-9.9954e+00 1.4885e+00
-9.9981e+00 1.3354e+00
-9.9992e+00 1.1877e+00
-9.9997e+00 1.0447e+00
-9.9999e+00 9.0575e-01
-1.0000e+01 7.7027e-01
-1.0000e+01 6.3767e-01
-1.0000e+01 5.0741e-01
-1.0000e+01 3.7899e-01
-1.0000e+01 2.5192e-01
-1.0000e+01 1.2574e-01
-1.0000e+01 3.6739e-16
-1.0000e+01 -1.2574e-01
-1.0000e+01 -2.5192e-01
-1.0000e+01 -3.7899e-01
-1.0000e+01 -5.0741e-01
-1.0000e+01 -6.3767e-01
-1.0000e+01 -7.7027e-01
-9.9999e+00 -9.0575e-01
-9.9997e+00 -1.0447e+00
-9.9992e+00 -1.1877e+00
-9.9981e+00 -1.3354e+00
-9.9954e+00 -1.4885e+00
-9.9896e+00 -1.6476e+00
-9.9776e+00 -1.8128e+00
-9.9535e+00 -1.9839e+00
-9.9070e+00 -2.1594e+00
-9.8736e+00 -2.2404e+00
-9.8298e+00 -2.3209e+00
-9.7729e+00 -2.4002e+00
-9.7001e+00 -2.4773e+00
-9.6085e+00 -2.5512e+00
-9.4956e+00 -2.6206e+00
-9.3594e+00 -2.6846e+00
-9.1991e+00 -2.7421e+00
-9.0150e+00 -2.7926e+00
-8.8086e+00 -2.8359e+00
-8.5829e+00 -2.8720e+00
-8.3412e+00 -2.9015e+00
-8.0874e+00 -2.9252e+00
-7.8251e+00 -2.9437e+00
-7.5577e+00 -2.9581e+00
-7.2879e+00 -2.9691e+00
-7.0179e+00 -2.9773e+00
-6.7494e+00 -2.9835e+00
-6.4835e+00 -2.9881e+00
-6.2211e+00 -2.9915e+00
-5.9627e+00 -2.9940e+00
-5.7086e+00 -2.9957e+00
-5.4588e+00 -2.9970e+00
-5.2134e+00 -2.9979e+00
-4.9724e+00 -2.9986e+00
-4.7355e+00 -2.9991e+00
-4.5026e+00 -2.9994e+00
-4.2736e+00 -2.9996e+00
-4.0482e+00 -2.9997e+00
-3.8262e+00 -2.9998e+00
-3.6073e+00 -2.9999e+00
-3.3915e+00 -2.9999e+00
-3.1784e+00 -3.0000e+00
-2.9679e+00 -2.9831e+00
-2.7598e+00 -2.9190e+00
-2.5539e+00 -2.8112e+00
-2.3500e+00 -2.6655e+00
-2.1479e+00 -2.4888e+00
-1.9475e+00 -2.2888e+00
-1.7486e+00 -2.0740e+00
-1.5510e+00 -1.8527e+00
-1.3546e+00 -1.6337e+00
-1.1592e+00 -1.4249e+00
-9.6469e-01 -1.2341e+00
-7.7089e-01 -1.0680e+00
-5.7767e-01 -9.3254e-01
-3.8487e-01 -8.3230e-01
-1.9237e-01 -7.7075e-01
-1.8370e-15 -7.5000e-01
1.9237e-01 -7.7075e-01
3.8487e-01 -8.3230e-01
5.7767e-01 -9.3254e-01
7.7089e-01 -1.0680e+00
9.6469e-01 -1.2341e+00
1.1592e+00 -1.4249e+00
1.3546e+00 -1.6337e+00
1.5510e+00 -1.8527e+00
1.7486e+00 -2.0740e+00
1.9475e+00 -2.2888e+00
2.1479e+00 -2.4888e+00
2.3500e+00 -2.6655e+00
2.5539e+00 -2.8112e+00
2.7598e+00 -2.9190e+00
2.9679e+00 -2.9831e+00
3.1784e+00 -3.0000e+00
3.3915e+00 -2.9999e+00
3.6073e+00 -2.9999e+00
3.8262e+00 -2.9998e+00
4.0482e+00 -2.9997e+00
4.2736e+00 -2.9996e+00
4.5026e+00 -2.9994e+00
4.7355e+00 -2.9991e+00
4.9724e+00 -2.9986e+00
5.2134e+00 -2.9979e+00
5.4588e+00 -2.9970e+00
5.7086e+00 -2.9957e+00
5.9627e+00 -2.9940e+00
6.2211e+00 -2.9915e+00
6.4835e+00 -2.9881e+00
6.7494e+00 -2.9835e+00
7.0179e+00 -2.9773e+00
7.2879e+00 -2.9691e+00
7.5577e+00 -2.9581e+00
7.8251e+00 -2.9437e+00
8.0874e+00 -2.9252e+00
8.3412e+00 -2.9015e+00
8.5829e+00 -2.8720e+00
8.8086e+00 -2.8359e+00
9.0150e+00 -2.7926e+00
9.1991e+00 -2.7421e+00
9.3594e+00 -2.6846e+00
9.4956e+00 -2.6206e+00
9.6085e+00 -2.5512e+00
9.7001e+00 -2.4773e+00
9.7729e+00 -2.4002e+00
9.8298e+00 -2.3209e+00
9.8736e+00 -2.2404e+00
9.9070e+00 -2.1594e+00
9.9535e+00 -1.9839e+00
9.9776e+00 -1.8128e+00
9.9896e+00 -1.6476e+00
9.9954e+00 -1.4885e+00
9.9981e+00 -1.3354e+00
9.9992e+00 -1.1877e+00
9.9997e+00 -1.0447e+00
9.9999e+00 -9.0575e-01
1.0000e+01 -7.7027e-01
1.0000e+01 -6.3767e-01
1.0000e+01 -5.0741e-01
1.0000e+01 -3.7899e-01
1.0000e+01 -2.5192e-01
1.0000e+01 -1.2574e-01
1.0000e+01 0.0000e+00
};

\addplot [mark=none,red,line width=1.5] table{
1.7149e+00 1.3256e+00
1.6969e+00 1.3027e+00
1.6801e+00 1.2782e+00
1.6639e+00 1.2520e+00
1.6478e+00 1.2241e+00
1.6311e+00 1.1945e+00
1.6131e+00 1.1632e+00
1.5930e+00 1.1305e+00
1.5702e+00 1.0965e+00
1.5443e+00 1.0615e+00
1.5154e+00 1.0251e+00
1.4840e+00 9.8679e-01
1.4521e+00 9.4510e-01
1.4223e+00 8.9825e-01
1.3987e+00 8.4487e-01
1.3858e+00 7.8515e-01
1.3866e+00 7.2123e-01
1.4015e+00 6.5597e-01
1.4271e+00 5.9081e-01
1.4570e+00 5.2441e-01
1.4819e+00 4.5360e-01
1.4896e+00 3.7737e-01
1.4673e+00 3.0297e-01
1.4113e+00 2.4779e-01
1.3339e+00 2.2703e-01
1.2526e+00 2.3829e-01
1.1730e+00 2.6431e-01
1.0908e+00 2.8387e-01
1.0062e+00 2.7722e-01
9.3220e-01 2.3535e-01
8.8155e-01 1.6647e-01
8.5515e-01 8.5057e-02
8.4733e-01 -1.7082e-12
8.5515e-01 -8.5057e-02
8.8155e-01 -1.6647e-01
9.3220e-01 -2.3535e-01
1.0062e+00 -2.7722e-01
1.0908e+00 -2.8387e-01
1.1730e+00 -2.6431e-01
1.2526e+00 -2.3829e-01
1.3339e+00 -2.2703e-01
1.4113e+00 -2.4779e-01
1.4673e+00 -3.0297e-01
1.4896e+00 -3.7737e-01
1.4819e+00 -4.5360e-01
1.4570e+00 -5.2441e-01
1.4271e+00 -5.9081e-01
1.4015e+00 -6.5597e-01
1.3866e+00 -7.2123e-01
1.3858e+00 -7.8515e-01
1.3987e+00 -8.4487e-01
1.4223e+00 -8.9825e-01
1.4521e+00 -9.4510e-01
1.4840e+00 -9.8679e-01
1.5154e+00 -1.0251e+00
1.5443e+00 -1.0615e+00
1.5702e+00 -1.0965e+00
1.5930e+00 -1.1305e+00
1.6131e+00 -1.1632e+00
1.6311e+00 -1.1945e+00
1.6478e+00 -1.2241e+00
1.6639e+00 -1.2520e+00
1.6801e+00 -1.2782e+00
1.6969e+00 -1.3027e+00
1.7149e+00 -1.3256e+00
1.7347e+00 -1.3470e+00
1.7569e+00 -1.3667e+00
1.7821e+00 -1.3844e+00
1.8107e+00 -1.3993e+00
1.8429e+00 -1.4107e+00
1.8786e+00 -1.4172e+00
1.9173e+00 -1.4180e+00
1.9581e+00 -1.4120e+00
2.0001e+00 -1.3989e+00
2.0422e+00 -1.3787e+00
2.0836e+00 -1.3516e+00
2.1241e+00 -1.3185e+00
2.1636e+00 -1.2800e+00
2.2023e+00 -1.2370e+00
2.2406e+00 -1.1901e+00
2.2789e+00 -1.1400e+00
2.3175e+00 -1.0870e+00
2.3567e+00 -1.0314e+00
2.3962e+00 -9.7343e-01
2.4360e+00 -9.1301e-01
2.4757e+00 -8.5013e-01
2.5148e+00 -7.8469e-01
2.5526e+00 -7.1662e-01
2.5887e+00 -6.4584e-01
2.6225e+00 -5.7237e-01
2.6532e+00 -4.9628e-01
2.6803e+00 -4.1774e-01
2.7034e+00 -3.3699e-01
2.7218e+00 -2.5438e-01
2.7352e+00 -1.7036e-01
2.7434e+00 -8.5401e-02
2.7462e+00 1.1813e-12
2.7434e+00 8.5401e-02
2.7352e+00 1.7036e-01
2.7218e+00 2.5438e-01
2.7034e+00 3.3699e-01
2.6803e+00 4.1774e-01
2.6532e+00 4.9628e-01
2.6225e+00 5.7237e-01
2.5887e+00 6.4584e-01
2.5526e+00 7.1662e-01
2.5148e+00 7.8469e-01
2.4757e+00 8.5013e-01
2.4360e+00 9.1301e-01
2.3962e+00 9.7343e-01
2.3567e+00 1.0314e+00
2.3175e+00 1.0870e+00
2.2789e+00 1.1400e+00
2.2406e+00 1.1901e+00
2.2023e+00 1.2370e+00
2.1636e+00 1.2800e+00
2.1241e+00 1.3185e+00
2.0836e+00 1.3516e+00
2.0422e+00 1.3787e+00
2.0001e+00 1.3989e+00
1.9581e+00 1.4120e+00
1.9173e+00 1.4180e+00
1.8786e+00 1.4172e+00
1.8429e+00 1.4107e+00
1.8107e+00 1.3993e+00
1.7821e+00 1.3844e+00
1.7569e+00 1.3667e+00
1.7347e+00 1.3470e+00
1.7149e+00 1.3256e+00
};

\end{axis}

\end{tikzpicture} \\
\begin{tikzpicture}[scale=0.5]

\begin{axis}[
  xmin = -11,
  xmax = 11,
  ymin = -3.2,
  ymax = 3.2,
  scale only axis,
  axis equal image,
  hide axis,
  title = {\Huge$t=11$}
  ]

\addplot [mark=none,black,line width=1.5] table{
1.0000e+01 0.0000e+00
1.0000e+01 1.2574e-01
1.0000e+01 2.5192e-01
1.0000e+01 3.7899e-01
1.0000e+01 5.0741e-01
1.0000e+01 6.3767e-01
1.0000e+01 7.7027e-01
9.9999e+00 9.0575e-01
9.9997e+00 1.0447e+00
9.9992e+00 1.1877e+00
9.9981e+00 1.3354e+00
9.9954e+00 1.4885e+00
9.9896e+00 1.6476e+00
9.9776e+00 1.8128e+00
9.9535e+00 1.9839e+00
9.9070e+00 2.1594e+00
9.8736e+00 2.2404e+00
9.8298e+00 2.3209e+00
9.7729e+00 2.4002e+00
9.7001e+00 2.4773e+00
9.6085e+00 2.5512e+00
9.4956e+00 2.6206e+00
9.3594e+00 2.6846e+00
9.1991e+00 2.7421e+00
9.0150e+00 2.7926e+00
8.8086e+00 2.8359e+00
8.5829e+00 2.8720e+00
8.3412e+00 2.9015e+00
8.0874e+00 2.9252e+00
7.8251e+00 2.9437e+00
7.5577e+00 2.9581e+00
7.2879e+00 2.9691e+00
7.0179e+00 2.9773e+00
6.7494e+00 2.9835e+00
6.4835e+00 2.9881e+00
6.2211e+00 2.9915e+00
5.9627e+00 2.9940e+00
5.7086e+00 2.9957e+00
5.4588e+00 2.9970e+00
5.2134e+00 2.9979e+00
4.9724e+00 2.9986e+00
4.7355e+00 2.9991e+00
4.5026e+00 2.9994e+00
4.2736e+00 2.9996e+00
4.0482e+00 2.9997e+00
3.8262e+00 2.9998e+00
3.6073e+00 2.9999e+00
3.3915e+00 2.9999e+00
3.1784e+00 3.0000e+00
2.9679e+00 2.9831e+00
2.7598e+00 2.9190e+00
2.5539e+00 2.8112e+00
2.3500e+00 2.6655e+00
2.1479e+00 2.4888e+00
1.9475e+00 2.2888e+00
1.7486e+00 2.0740e+00
1.5510e+00 1.8527e+00
1.3546e+00 1.6337e+00
1.1592e+00 1.4249e+00
9.6469e-01 1.2341e+00
7.7089e-01 1.0680e+00
5.7767e-01 9.3254e-01
3.8487e-01 8.3230e-01
1.9237e-01 7.7075e-01
6.1232e-16 7.5000e-01
-1.9237e-01 7.7075e-01
-3.8487e-01 8.3230e-01
-5.7767e-01 9.3254e-01
-7.7089e-01 1.0680e+00
-9.6469e-01 1.2341e+00
-1.1592e+00 1.4249e+00
-1.3546e+00 1.6337e+00
-1.5510e+00 1.8527e+00
-1.7486e+00 2.0740e+00
-1.9475e+00 2.2888e+00
-2.1479e+00 2.4888e+00
-2.3500e+00 2.6655e+00
-2.5539e+00 2.8112e+00
-2.7598e+00 2.9190e+00
-2.9679e+00 2.9831e+00
-3.1784e+00 3.0000e+00
-3.3915e+00 2.9999e+00
-3.6073e+00 2.9999e+00
-3.8262e+00 2.9998e+00
-4.0482e+00 2.9997e+00
-4.2736e+00 2.9996e+00
-4.5026e+00 2.9994e+00
-4.7355e+00 2.9991e+00
-4.9724e+00 2.9986e+00
-5.2134e+00 2.9979e+00
-5.4588e+00 2.9970e+00
-5.7086e+00 2.9957e+00
-5.9627e+00 2.9940e+00
-6.2211e+00 2.9915e+00
-6.4835e+00 2.9881e+00
-6.7494e+00 2.9835e+00
-7.0179e+00 2.9773e+00
-7.2879e+00 2.9691e+00
-7.5577e+00 2.9581e+00
-7.8251e+00 2.9437e+00
-8.0874e+00 2.9252e+00
-8.3412e+00 2.9015e+00
-8.5829e+00 2.8720e+00
-8.8086e+00 2.8359e+00
-9.0150e+00 2.7926e+00
-9.1991e+00 2.7421e+00
-9.3594e+00 2.6846e+00
-9.4956e+00 2.6206e+00
-9.6085e+00 2.5512e+00
-9.7001e+00 2.4773e+00
-9.7729e+00 2.4002e+00
-9.8298e+00 2.3209e+00
-9.8736e+00 2.2404e+00
-9.9070e+00 2.1594e+00
-9.9535e+00 1.9839e+00
-9.9776e+00 1.8128e+00
-9.9896e+00 1.6476e+00
-9.9954e+00 1.4885e+00
-9.9981e+00 1.3354e+00
-9.9992e+00 1.1877e+00
-9.9997e+00 1.0447e+00
-9.9999e+00 9.0575e-01
-1.0000e+01 7.7027e-01
-1.0000e+01 6.3767e-01
-1.0000e+01 5.0741e-01
-1.0000e+01 3.7899e-01
-1.0000e+01 2.5192e-01
-1.0000e+01 1.2574e-01
-1.0000e+01 3.6739e-16
-1.0000e+01 -1.2574e-01
-1.0000e+01 -2.5192e-01
-1.0000e+01 -3.7899e-01
-1.0000e+01 -5.0741e-01
-1.0000e+01 -6.3767e-01
-1.0000e+01 -7.7027e-01
-9.9999e+00 -9.0575e-01
-9.9997e+00 -1.0447e+00
-9.9992e+00 -1.1877e+00
-9.9981e+00 -1.3354e+00
-9.9954e+00 -1.4885e+00
-9.9896e+00 -1.6476e+00
-9.9776e+00 -1.8128e+00
-9.9535e+00 -1.9839e+00
-9.9070e+00 -2.1594e+00
-9.8736e+00 -2.2404e+00
-9.8298e+00 -2.3209e+00
-9.7729e+00 -2.4002e+00
-9.7001e+00 -2.4773e+00
-9.6085e+00 -2.5512e+00
-9.4956e+00 -2.6206e+00
-9.3594e+00 -2.6846e+00
-9.1991e+00 -2.7421e+00
-9.0150e+00 -2.7926e+00
-8.8086e+00 -2.8359e+00
-8.5829e+00 -2.8720e+00
-8.3412e+00 -2.9015e+00
-8.0874e+00 -2.9252e+00
-7.8251e+00 -2.9437e+00
-7.5577e+00 -2.9581e+00
-7.2879e+00 -2.9691e+00
-7.0179e+00 -2.9773e+00
-6.7494e+00 -2.9835e+00
-6.4835e+00 -2.9881e+00
-6.2211e+00 -2.9915e+00
-5.9627e+00 -2.9940e+00
-5.7086e+00 -2.9957e+00
-5.4588e+00 -2.9970e+00
-5.2134e+00 -2.9979e+00
-4.9724e+00 -2.9986e+00
-4.7355e+00 -2.9991e+00
-4.5026e+00 -2.9994e+00
-4.2736e+00 -2.9996e+00
-4.0482e+00 -2.9997e+00
-3.8262e+00 -2.9998e+00
-3.6073e+00 -2.9999e+00
-3.3915e+00 -2.9999e+00
-3.1784e+00 -3.0000e+00
-2.9679e+00 -2.9831e+00
-2.7598e+00 -2.9190e+00
-2.5539e+00 -2.8112e+00
-2.3500e+00 -2.6655e+00
-2.1479e+00 -2.4888e+00
-1.9475e+00 -2.2888e+00
-1.7486e+00 -2.0740e+00
-1.5510e+00 -1.8527e+00
-1.3546e+00 -1.6337e+00
-1.1592e+00 -1.4249e+00
-9.6469e-01 -1.2341e+00
-7.7089e-01 -1.0680e+00
-5.7767e-01 -9.3254e-01
-3.8487e-01 -8.3230e-01
-1.9237e-01 -7.7075e-01
-1.8370e-15 -7.5000e-01
1.9237e-01 -7.7075e-01
3.8487e-01 -8.3230e-01
5.7767e-01 -9.3254e-01
7.7089e-01 -1.0680e+00
9.6469e-01 -1.2341e+00
1.1592e+00 -1.4249e+00
1.3546e+00 -1.6337e+00
1.5510e+00 -1.8527e+00
1.7486e+00 -2.0740e+00
1.9475e+00 -2.2888e+00
2.1479e+00 -2.4888e+00
2.3500e+00 -2.6655e+00
2.5539e+00 -2.8112e+00
2.7598e+00 -2.9190e+00
2.9679e+00 -2.9831e+00
3.1784e+00 -3.0000e+00
3.3915e+00 -2.9999e+00
3.6073e+00 -2.9999e+00
3.8262e+00 -2.9998e+00
4.0482e+00 -2.9997e+00
4.2736e+00 -2.9996e+00
4.5026e+00 -2.9994e+00
4.7355e+00 -2.9991e+00
4.9724e+00 -2.9986e+00
5.2134e+00 -2.9979e+00
5.4588e+00 -2.9970e+00
5.7086e+00 -2.9957e+00
5.9627e+00 -2.9940e+00
6.2211e+00 -2.9915e+00
6.4835e+00 -2.9881e+00
6.7494e+00 -2.9835e+00
7.0179e+00 -2.9773e+00
7.2879e+00 -2.9691e+00
7.5577e+00 -2.9581e+00
7.8251e+00 -2.9437e+00
8.0874e+00 -2.9252e+00
8.3412e+00 -2.9015e+00
8.5829e+00 -2.8720e+00
8.8086e+00 -2.8359e+00
9.0150e+00 -2.7926e+00
9.1991e+00 -2.7421e+00
9.3594e+00 -2.6846e+00
9.4956e+00 -2.6206e+00
9.6085e+00 -2.5512e+00
9.7001e+00 -2.4773e+00
9.7729e+00 -2.4002e+00
9.8298e+00 -2.3209e+00
9.8736e+00 -2.2404e+00
9.9070e+00 -2.1594e+00
9.9535e+00 -1.9839e+00
9.9776e+00 -1.8128e+00
9.9896e+00 -1.6476e+00
9.9954e+00 -1.4885e+00
9.9981e+00 -1.3354e+00
9.9992e+00 -1.1877e+00
9.9997e+00 -1.0447e+00
9.9999e+00 -9.0575e-01
1.0000e+01 -7.7027e-01
1.0000e+01 -6.3767e-01
1.0000e+01 -5.0741e-01
1.0000e+01 -3.7899e-01
1.0000e+01 -2.5192e-01
1.0000e+01 -1.2574e-01
1.0000e+01 0.0000e+00
};

\addplot [mark=none,red,line width=1.5] table{
2.7299e+00 1.6121e+00
2.7023e+00 1.6225e+00
2.6736e+00 1.6315e+00
2.6434e+00 1.6387e+00
2.6113e+00 1.6438e+00
2.5771e+00 1.6461e+00
2.5408e+00 1.6447e+00
2.5028e+00 1.6384e+00
2.4638e+00 1.6258e+00
2.4251e+00 1.6056e+00
2.3887e+00 1.5767e+00
2.3570e+00 1.5384e+00
2.3331e+00 1.4913e+00
2.3191e+00 1.4371e+00
2.3160e+00 1.3785e+00
2.3234e+00 1.3174e+00
2.3396e+00 1.2553e+00
2.3625e+00 1.1922e+00
2.3897e+00 1.1276e+00
2.4186e+00 1.0608e+00
2.4467e+00 9.9102e-01
2.4714e+00 9.1811e-01
2.4909e+00 8.4205e-01
2.5039e+00 7.6305e-01
2.5098e+00 6.8143e-01
2.5088e+00 5.9774e-01
2.5019e+00 5.1282e-01
2.4909e+00 4.2753e-01
2.4779e+00 3.4226e-01
2.4653e+00 2.5685e-01
2.4548e+00 1.7119e-01
2.4479e+00 8.5518e-02
2.4456e+00 -1.4671e-12
2.4479e+00 -8.5518e-02
2.4548e+00 -1.7119e-01
2.4653e+00 -2.5685e-01
2.4779e+00 -3.4226e-01
2.4909e+00 -4.2753e-01
2.5019e+00 -5.1282e-01
2.5088e+00 -5.9774e-01
2.5098e+00 -6.8143e-01
2.5039e+00 -7.6305e-01
2.4909e+00 -8.4205e-01
2.4714e+00 -9.1811e-01
2.4467e+00 -9.9102e-01
2.4186e+00 -1.0608e+00
2.3897e+00 -1.1276e+00
2.3625e+00 -1.1922e+00
2.3396e+00 -1.2553e+00
2.3234e+00 -1.3174e+00
2.3160e+00 -1.3785e+00
2.3191e+00 -1.4371e+00
2.3331e+00 -1.4913e+00
2.3570e+00 -1.5384e+00
2.3887e+00 -1.5767e+00
2.4251e+00 -1.6056e+00
2.4638e+00 -1.6258e+00
2.5028e+00 -1.6384e+00
2.5408e+00 -1.6447e+00
2.5771e+00 -1.6461e+00
2.6113e+00 -1.6438e+00
2.6434e+00 -1.6387e+00
2.6736e+00 -1.6315e+00
2.7023e+00 -1.6225e+00
2.7299e+00 -1.6121e+00
2.7569e+00 -1.6003e+00
2.7839e+00 -1.5870e+00
2.8111e+00 -1.5720e+00
2.8390e+00 -1.5552e+00
2.8678e+00 -1.5365e+00
2.8976e+00 -1.5154e+00
2.9286e+00 -1.4920e+00
2.9609e+00 -1.4661e+00
2.9943e+00 -1.4375e+00
3.0290e+00 -1.4062e+00
3.0649e+00 -1.3720e+00
3.1018e+00 -1.3348e+00
3.1397e+00 -1.2947e+00
3.1783e+00 -1.2517e+00
3.2175e+00 -1.2055e+00
3.2571e+00 -1.1564e+00
3.2967e+00 -1.1041e+00
3.3362e+00 -1.0488e+00
3.3752e+00 -9.9037e-01
3.4133e+00 -9.2891e-01
3.4504e+00 -8.6442e-01
3.4859e+00 -7.9699e-01
3.5195e+00 -7.2669e-01
3.5508e+00 -6.5367e-01
3.5795e+00 -5.7809e-01
3.6052e+00 -5.0016e-01
3.6276e+00 -4.2012e-01
3.6464e+00 -3.3827e-01
3.6613e+00 -2.5494e-01
3.6721e+00 -1.7053e-01
3.6786e+00 -8.5423e-02
3.6808e+00 7.5402e-13
3.6786e+00 8.5423e-02
3.6721e+00 1.7053e-01
3.6613e+00 2.5494e-01
3.6464e+00 3.3827e-01
3.6276e+00 4.2012e-01
3.6052e+00 5.0016e-01
3.5795e+00 5.7809e-01
3.5508e+00 6.5367e-01
3.5195e+00 7.2669e-01
3.4859e+00 7.9699e-01
3.4504e+00 8.6442e-01
3.4133e+00 9.2891e-01
3.3752e+00 9.9037e-01
3.3362e+00 1.0488e+00
3.2967e+00 1.1041e+00
3.2571e+00 1.1564e+00
3.2175e+00 1.2055e+00
3.1783e+00 1.2517e+00
3.1397e+00 1.2947e+00
3.1018e+00 1.3348e+00
3.0649e+00 1.3720e+00
3.0290e+00 1.4062e+00
2.9943e+00 1.4375e+00
2.9609e+00 1.4661e+00
2.9286e+00 1.4920e+00
2.8976e+00 1.5154e+00
2.8678e+00 1.5365e+00
2.8390e+00 1.5552e+00
2.8111e+00 1.5720e+00
2.7839e+00 1.5870e+00
2.7569e+00 1.6003e+00
2.7299e+00 1.6121e+00
};

\end{axis}

\end{tikzpicture} &
\begin{tikzpicture}[scale=0.5]

\begin{axis}[
  xmin = -11,
  xmax = 11,
  ymin = -3.2,
  ymax = 3.2,
  scale only axis,
  axis equal image,
  hide axis,
  title = {\Huge$t=12$}
  ]

\addplot [mark=none,black,line width=1.5] table{
1.0000e+01 0.0000e+00
1.0000e+01 1.2574e-01
1.0000e+01 2.5192e-01
1.0000e+01 3.7899e-01
1.0000e+01 5.0741e-01
1.0000e+01 6.3767e-01
1.0000e+01 7.7027e-01
9.9999e+00 9.0575e-01
9.9997e+00 1.0447e+00
9.9992e+00 1.1877e+00
9.9981e+00 1.3354e+00
9.9954e+00 1.4885e+00
9.9896e+00 1.6476e+00
9.9776e+00 1.8128e+00
9.9535e+00 1.9839e+00
9.9070e+00 2.1594e+00
9.8736e+00 2.2404e+00
9.8298e+00 2.3209e+00
9.7729e+00 2.4002e+00
9.7001e+00 2.4773e+00
9.6085e+00 2.5512e+00
9.4956e+00 2.6206e+00
9.3594e+00 2.6846e+00
9.1991e+00 2.7421e+00
9.0150e+00 2.7926e+00
8.8086e+00 2.8359e+00
8.5829e+00 2.8720e+00
8.3412e+00 2.9015e+00
8.0874e+00 2.9252e+00
7.8251e+00 2.9437e+00
7.5577e+00 2.9581e+00
7.2879e+00 2.9691e+00
7.0179e+00 2.9773e+00
6.7494e+00 2.9835e+00
6.4835e+00 2.9881e+00
6.2211e+00 2.9915e+00
5.9627e+00 2.9940e+00
5.7086e+00 2.9957e+00
5.4588e+00 2.9970e+00
5.2134e+00 2.9979e+00
4.9724e+00 2.9986e+00
4.7355e+00 2.9991e+00
4.5026e+00 2.9994e+00
4.2736e+00 2.9996e+00
4.0482e+00 2.9997e+00
3.8262e+00 2.9998e+00
3.6073e+00 2.9999e+00
3.3915e+00 2.9999e+00
3.1784e+00 3.0000e+00
2.9679e+00 2.9831e+00
2.7598e+00 2.9190e+00
2.5539e+00 2.8112e+00
2.3500e+00 2.6655e+00
2.1479e+00 2.4888e+00
1.9475e+00 2.2888e+00
1.7486e+00 2.0740e+00
1.5510e+00 1.8527e+00
1.3546e+00 1.6337e+00
1.1592e+00 1.4249e+00
9.6469e-01 1.2341e+00
7.7089e-01 1.0680e+00
5.7767e-01 9.3254e-01
3.8487e-01 8.3230e-01
1.9237e-01 7.7075e-01
6.1232e-16 7.5000e-01
-1.9237e-01 7.7075e-01
-3.8487e-01 8.3230e-01
-5.7767e-01 9.3254e-01
-7.7089e-01 1.0680e+00
-9.6469e-01 1.2341e+00
-1.1592e+00 1.4249e+00
-1.3546e+00 1.6337e+00
-1.5510e+00 1.8527e+00
-1.7486e+00 2.0740e+00
-1.9475e+00 2.2888e+00
-2.1479e+00 2.4888e+00
-2.3500e+00 2.6655e+00
-2.5539e+00 2.8112e+00
-2.7598e+00 2.9190e+00
-2.9679e+00 2.9831e+00
-3.1784e+00 3.0000e+00
-3.3915e+00 2.9999e+00
-3.6073e+00 2.9999e+00
-3.8262e+00 2.9998e+00
-4.0482e+00 2.9997e+00
-4.2736e+00 2.9996e+00
-4.5026e+00 2.9994e+00
-4.7355e+00 2.9991e+00
-4.9724e+00 2.9986e+00
-5.2134e+00 2.9979e+00
-5.4588e+00 2.9970e+00
-5.7086e+00 2.9957e+00
-5.9627e+00 2.9940e+00
-6.2211e+00 2.9915e+00
-6.4835e+00 2.9881e+00
-6.7494e+00 2.9835e+00
-7.0179e+00 2.9773e+00
-7.2879e+00 2.9691e+00
-7.5577e+00 2.9581e+00
-7.8251e+00 2.9437e+00
-8.0874e+00 2.9252e+00
-8.3412e+00 2.9015e+00
-8.5829e+00 2.8720e+00
-8.8086e+00 2.8359e+00
-9.0150e+00 2.7926e+00
-9.1991e+00 2.7421e+00
-9.3594e+00 2.6846e+00
-9.4956e+00 2.6206e+00
-9.6085e+00 2.5512e+00
-9.7001e+00 2.4773e+00
-9.7729e+00 2.4002e+00
-9.8298e+00 2.3209e+00
-9.8736e+00 2.2404e+00
-9.9070e+00 2.1594e+00
-9.9535e+00 1.9839e+00
-9.9776e+00 1.8128e+00
-9.9896e+00 1.6476e+00
-9.9954e+00 1.4885e+00
-9.9981e+00 1.3354e+00
-9.9992e+00 1.1877e+00
-9.9997e+00 1.0447e+00
-9.9999e+00 9.0575e-01
-1.0000e+01 7.7027e-01
-1.0000e+01 6.3767e-01
-1.0000e+01 5.0741e-01
-1.0000e+01 3.7899e-01
-1.0000e+01 2.5192e-01
-1.0000e+01 1.2574e-01
-1.0000e+01 3.6739e-16
-1.0000e+01 -1.2574e-01
-1.0000e+01 -2.5192e-01
-1.0000e+01 -3.7899e-01
-1.0000e+01 -5.0741e-01
-1.0000e+01 -6.3767e-01
-1.0000e+01 -7.7027e-01
-9.9999e+00 -9.0575e-01
-9.9997e+00 -1.0447e+00
-9.9992e+00 -1.1877e+00
-9.9981e+00 -1.3354e+00
-9.9954e+00 -1.4885e+00
-9.9896e+00 -1.6476e+00
-9.9776e+00 -1.8128e+00
-9.9535e+00 -1.9839e+00
-9.9070e+00 -2.1594e+00
-9.8736e+00 -2.2404e+00
-9.8298e+00 -2.3209e+00
-9.7729e+00 -2.4002e+00
-9.7001e+00 -2.4773e+00
-9.6085e+00 -2.5512e+00
-9.4956e+00 -2.6206e+00
-9.3594e+00 -2.6846e+00
-9.1991e+00 -2.7421e+00
-9.0150e+00 -2.7926e+00
-8.8086e+00 -2.8359e+00
-8.5829e+00 -2.8720e+00
-8.3412e+00 -2.9015e+00
-8.0874e+00 -2.9252e+00
-7.8251e+00 -2.9437e+00
-7.5577e+00 -2.9581e+00
-7.2879e+00 -2.9691e+00
-7.0179e+00 -2.9773e+00
-6.7494e+00 -2.9835e+00
-6.4835e+00 -2.9881e+00
-6.2211e+00 -2.9915e+00
-5.9627e+00 -2.9940e+00
-5.7086e+00 -2.9957e+00
-5.4588e+00 -2.9970e+00
-5.2134e+00 -2.9979e+00
-4.9724e+00 -2.9986e+00
-4.7355e+00 -2.9991e+00
-4.5026e+00 -2.9994e+00
-4.2736e+00 -2.9996e+00
-4.0482e+00 -2.9997e+00
-3.8262e+00 -2.9998e+00
-3.6073e+00 -2.9999e+00
-3.3915e+00 -2.9999e+00
-3.1784e+00 -3.0000e+00
-2.9679e+00 -2.9831e+00
-2.7598e+00 -2.9190e+00
-2.5539e+00 -2.8112e+00
-2.3500e+00 -2.6655e+00
-2.1479e+00 -2.4888e+00
-1.9475e+00 -2.2888e+00
-1.7486e+00 -2.0740e+00
-1.5510e+00 -1.8527e+00
-1.3546e+00 -1.6337e+00
-1.1592e+00 -1.4249e+00
-9.6469e-01 -1.2341e+00
-7.7089e-01 -1.0680e+00
-5.7767e-01 -9.3254e-01
-3.8487e-01 -8.3230e-01
-1.9237e-01 -7.7075e-01
-1.8370e-15 -7.5000e-01
1.9237e-01 -7.7075e-01
3.8487e-01 -8.3230e-01
5.7767e-01 -9.3254e-01
7.7089e-01 -1.0680e+00
9.6469e-01 -1.2341e+00
1.1592e+00 -1.4249e+00
1.3546e+00 -1.6337e+00
1.5510e+00 -1.8527e+00
1.7486e+00 -2.0740e+00
1.9475e+00 -2.2888e+00
2.1479e+00 -2.4888e+00
2.3500e+00 -2.6655e+00
2.5539e+00 -2.8112e+00
2.7598e+00 -2.9190e+00
2.9679e+00 -2.9831e+00
3.1784e+00 -3.0000e+00
3.3915e+00 -2.9999e+00
3.6073e+00 -2.9999e+00
3.8262e+00 -2.9998e+00
4.0482e+00 -2.9997e+00
4.2736e+00 -2.9996e+00
4.5026e+00 -2.9994e+00
4.7355e+00 -2.9991e+00
4.9724e+00 -2.9986e+00
5.2134e+00 -2.9979e+00
5.4588e+00 -2.9970e+00
5.7086e+00 -2.9957e+00
5.9627e+00 -2.9940e+00
6.2211e+00 -2.9915e+00
6.4835e+00 -2.9881e+00
6.7494e+00 -2.9835e+00
7.0179e+00 -2.9773e+00
7.2879e+00 -2.9691e+00
7.5577e+00 -2.9581e+00
7.8251e+00 -2.9437e+00
8.0874e+00 -2.9252e+00
8.3412e+00 -2.9015e+00
8.5829e+00 -2.8720e+00
8.8086e+00 -2.8359e+00
9.0150e+00 -2.7926e+00
9.1991e+00 -2.7421e+00
9.3594e+00 -2.6846e+00
9.4956e+00 -2.6206e+00
9.6085e+00 -2.5512e+00
9.7001e+00 -2.4773e+00
9.7729e+00 -2.4002e+00
9.8298e+00 -2.3209e+00
9.8736e+00 -2.2404e+00
9.9070e+00 -2.1594e+00
9.9535e+00 -1.9839e+00
9.9776e+00 -1.8128e+00
9.9896e+00 -1.6476e+00
9.9954e+00 -1.4885e+00
9.9981e+00 -1.3354e+00
9.9992e+00 -1.1877e+00
9.9997e+00 -1.0447e+00
9.9999e+00 -9.0575e-01
1.0000e+01 -7.7027e-01
1.0000e+01 -6.3767e-01
1.0000e+01 -5.0741e-01
1.0000e+01 -3.7899e-01
1.0000e+01 -2.5192e-01
1.0000e+01 -1.2574e-01
1.0000e+01 0.0000e+00
};

\addplot [mark=none,red,line width=1.5] table{
3.6199e+00 1.6038e+00
3.5920e+00 1.6134e+00
3.5630e+00 1.6214e+00
3.5326e+00 1.6276e+00
3.5003e+00 1.6316e+00
3.4661e+00 1.6327e+00
3.4299e+00 1.6301e+00
3.3920e+00 1.6229e+00
3.3531e+00 1.6100e+00
3.3142e+00 1.5902e+00
3.2766e+00 1.5627e+00
3.2423e+00 1.5267e+00
3.2134e+00 1.4824e+00
3.1920e+00 1.4307e+00
3.1793e+00 1.3733e+00
3.1758e+00 1.3119e+00
3.1808e+00 1.2479e+00
3.1934e+00 1.1820e+00
3.2121e+00 1.1144e+00
3.2354e+00 1.0454e+00
3.2617e+00 9.7496e-01
3.2893e+00 9.0305e-01
3.3169e+00 8.2951e-01
3.3435e+00 7.5398e-01
3.3686e+00 6.7604e-01
3.3916e+00 5.9547e-01
3.4118e+00 5.1264e-01
3.4290e+00 4.2830e-01
3.4430e+00 3.4312e-01
3.4538e+00 2.5739e-01
3.4614e+00 1.7140e-01
3.4658e+00 8.5546e-02
3.4673e+00 -1.2474e-12
3.4658e+00 -8.5546e-02
3.4614e+00 -1.7140e-01
3.4538e+00 -2.5739e-01
3.4430e+00 -3.4312e-01
3.4290e+00 -4.2830e-01
3.4118e+00 -5.1264e-01
3.3916e+00 -5.9547e-01
3.3686e+00 -6.7604e-01
3.3435e+00 -7.5398e-01
3.3169e+00 -8.2951e-01
3.2893e+00 -9.0305e-01
3.2617e+00 -9.7496e-01
3.2354e+00 -1.0454e+00
3.2121e+00 -1.1144e+00
3.1934e+00 -1.1820e+00
3.1808e+00 -1.2479e+00
3.1758e+00 -1.3119e+00
3.1793e+00 -1.3733e+00
3.1920e+00 -1.4307e+00
3.2134e+00 -1.4824e+00
3.2423e+00 -1.5267e+00
3.2766e+00 -1.5627e+00
3.3142e+00 -1.5902e+00
3.3531e+00 -1.6100e+00
3.3920e+00 -1.6229e+00
3.4299e+00 -1.6301e+00
3.4661e+00 -1.6327e+00
3.5003e+00 -1.6316e+00
3.5326e+00 -1.6276e+00
3.5630e+00 -1.6214e+00
3.5920e+00 -1.6134e+00
3.6199e+00 -1.6038e+00
3.6472e+00 -1.5926e+00
3.6744e+00 -1.5798e+00
3.7018e+00 -1.5653e+00
3.7298e+00 -1.5487e+00
3.7587e+00 -1.5300e+00
3.7885e+00 -1.5090e+00
3.8193e+00 -1.4854e+00
3.8513e+00 -1.4591e+00
3.8844e+00 -1.4301e+00
3.9187e+00 -1.3983e+00
3.9540e+00 -1.3635e+00
3.9905e+00 -1.3259e+00
4.0278e+00 -1.2854e+00
4.0661e+00 -1.2420e+00
4.1051e+00 -1.1957e+00
4.1447e+00 -1.1465e+00
4.1847e+00 -1.0945e+00
4.2248e+00 -1.0397e+00
4.2648e+00 -9.8195e-01
4.3043e+00 -9.2133e-01
4.3430e+00 -8.5782e-01
4.3804e+00 -7.9144e-01
4.4162e+00 -7.2222e-01
4.4498e+00 -6.5023e-01
4.4809e+00 -5.7560e-01
4.5090e+00 -4.9847e-01
4.5336e+00 -4.1908e-01
4.5543e+00 -3.3771e-01
4.5708e+00 -2.5470e-01
4.5828e+00 -1.7045e-01
4.5901e+00 -8.5413e-02
4.5926e+00 7.6143e-13
4.5901e+00 8.5413e-02
4.5828e+00 1.7045e-01
4.5708e+00 2.5470e-01
4.5543e+00 3.3771e-01
4.5336e+00 4.1908e-01
4.5090e+00 4.9847e-01
4.4809e+00 5.7560e-01
4.4498e+00 6.5023e-01
4.4162e+00 7.2222e-01
4.3804e+00 7.9144e-01
4.3430e+00 8.5782e-01
4.3043e+00 9.2133e-01
4.2648e+00 9.8195e-01
4.2248e+00 1.0397e+00
4.1847e+00 1.0945e+00
4.1447e+00 1.1465e+00
4.1051e+00 1.1957e+00
4.0661e+00 1.2420e+00
4.0278e+00 1.2854e+00
3.9905e+00 1.3259e+00
3.9540e+00 1.3635e+00
3.9187e+00 1.3983e+00
3.8844e+00 1.4301e+00
3.8513e+00 1.4591e+00
3.8193e+00 1.4854e+00
3.7885e+00 1.5090e+00
3.7587e+00 1.5300e+00
3.7298e+00 1.5487e+00
3.7018e+00 1.5653e+00
3.6744e+00 1.5798e+00
3.6472e+00 1.5926e+00
3.6199e+00 1.6038e+00
};

\end{axis}

\end{tikzpicture} &
\begin{tikzpicture}[scale=0.5]

\begin{axis}[
  xmin = -11,
  xmax = 11,
  ymin = -3.2,
  ymax = 3.2,
  scale only axis,
  axis equal image,
  hide axis,
  title = {\Huge$t=13$}
  ]

\addplot [mark=none,black,line width=1.5] table{
1.0000e+01 0.0000e+00
1.0000e+01 1.2574e-01
1.0000e+01 2.5192e-01
1.0000e+01 3.7899e-01
1.0000e+01 5.0741e-01
1.0000e+01 6.3767e-01
1.0000e+01 7.7027e-01
9.9999e+00 9.0575e-01
9.9997e+00 1.0447e+00
9.9992e+00 1.1877e+00
9.9981e+00 1.3354e+00
9.9954e+00 1.4885e+00
9.9896e+00 1.6476e+00
9.9776e+00 1.8128e+00
9.9535e+00 1.9839e+00
9.9070e+00 2.1594e+00
9.8736e+00 2.2404e+00
9.8298e+00 2.3209e+00
9.7729e+00 2.4002e+00
9.7001e+00 2.4773e+00
9.6085e+00 2.5512e+00
9.4956e+00 2.6206e+00
9.3594e+00 2.6846e+00
9.1991e+00 2.7421e+00
9.0150e+00 2.7926e+00
8.8086e+00 2.8359e+00
8.5829e+00 2.8720e+00
8.3412e+00 2.9015e+00
8.0874e+00 2.9252e+00
7.8251e+00 2.9437e+00
7.5577e+00 2.9581e+00
7.2879e+00 2.9691e+00
7.0179e+00 2.9773e+00
6.7494e+00 2.9835e+00
6.4835e+00 2.9881e+00
6.2211e+00 2.9915e+00
5.9627e+00 2.9940e+00
5.7086e+00 2.9957e+00
5.4588e+00 2.9970e+00
5.2134e+00 2.9979e+00
4.9724e+00 2.9986e+00
4.7355e+00 2.9991e+00
4.5026e+00 2.9994e+00
4.2736e+00 2.9996e+00
4.0482e+00 2.9997e+00
3.8262e+00 2.9998e+00
3.6073e+00 2.9999e+00
3.3915e+00 2.9999e+00
3.1784e+00 3.0000e+00
2.9679e+00 2.9831e+00
2.7598e+00 2.9190e+00
2.5539e+00 2.8112e+00
2.3500e+00 2.6655e+00
2.1479e+00 2.4888e+00
1.9475e+00 2.2888e+00
1.7486e+00 2.0740e+00
1.5510e+00 1.8527e+00
1.3546e+00 1.6337e+00
1.1592e+00 1.4249e+00
9.6469e-01 1.2341e+00
7.7089e-01 1.0680e+00
5.7767e-01 9.3254e-01
3.8487e-01 8.3230e-01
1.9237e-01 7.7075e-01
6.1232e-16 7.5000e-01
-1.9237e-01 7.7075e-01
-3.8487e-01 8.3230e-01
-5.7767e-01 9.3254e-01
-7.7089e-01 1.0680e+00
-9.6469e-01 1.2341e+00
-1.1592e+00 1.4249e+00
-1.3546e+00 1.6337e+00
-1.5510e+00 1.8527e+00
-1.7486e+00 2.0740e+00
-1.9475e+00 2.2888e+00
-2.1479e+00 2.4888e+00
-2.3500e+00 2.6655e+00
-2.5539e+00 2.8112e+00
-2.7598e+00 2.9190e+00
-2.9679e+00 2.9831e+00
-3.1784e+00 3.0000e+00
-3.3915e+00 2.9999e+00
-3.6073e+00 2.9999e+00
-3.8262e+00 2.9998e+00
-4.0482e+00 2.9997e+00
-4.2736e+00 2.9996e+00
-4.5026e+00 2.9994e+00
-4.7355e+00 2.9991e+00
-4.9724e+00 2.9986e+00
-5.2134e+00 2.9979e+00
-5.4588e+00 2.9970e+00
-5.7086e+00 2.9957e+00
-5.9627e+00 2.9940e+00
-6.2211e+00 2.9915e+00
-6.4835e+00 2.9881e+00
-6.7494e+00 2.9835e+00
-7.0179e+00 2.9773e+00
-7.2879e+00 2.9691e+00
-7.5577e+00 2.9581e+00
-7.8251e+00 2.9437e+00
-8.0874e+00 2.9252e+00
-8.3412e+00 2.9015e+00
-8.5829e+00 2.8720e+00
-8.8086e+00 2.8359e+00
-9.0150e+00 2.7926e+00
-9.1991e+00 2.7421e+00
-9.3594e+00 2.6846e+00
-9.4956e+00 2.6206e+00
-9.6085e+00 2.5512e+00
-9.7001e+00 2.4773e+00
-9.7729e+00 2.4002e+00
-9.8298e+00 2.3209e+00
-9.8736e+00 2.2404e+00
-9.9070e+00 2.1594e+00
-9.9535e+00 1.9839e+00
-9.9776e+00 1.8128e+00
-9.9896e+00 1.6476e+00
-9.9954e+00 1.4885e+00
-9.9981e+00 1.3354e+00
-9.9992e+00 1.1877e+00
-9.9997e+00 1.0447e+00
-9.9999e+00 9.0575e-01
-1.0000e+01 7.7027e-01
-1.0000e+01 6.3767e-01
-1.0000e+01 5.0741e-01
-1.0000e+01 3.7899e-01
-1.0000e+01 2.5192e-01
-1.0000e+01 1.2574e-01
-1.0000e+01 3.6739e-16
-1.0000e+01 -1.2574e-01
-1.0000e+01 -2.5192e-01
-1.0000e+01 -3.7899e-01
-1.0000e+01 -5.0741e-01
-1.0000e+01 -6.3767e-01
-1.0000e+01 -7.7027e-01
-9.9999e+00 -9.0575e-01
-9.9997e+00 -1.0447e+00
-9.9992e+00 -1.1877e+00
-9.9981e+00 -1.3354e+00
-9.9954e+00 -1.4885e+00
-9.9896e+00 -1.6476e+00
-9.9776e+00 -1.8128e+00
-9.9535e+00 -1.9839e+00
-9.9070e+00 -2.1594e+00
-9.8736e+00 -2.2404e+00
-9.8298e+00 -2.3209e+00
-9.7729e+00 -2.4002e+00
-9.7001e+00 -2.4773e+00
-9.6085e+00 -2.5512e+00
-9.4956e+00 -2.6206e+00
-9.3594e+00 -2.6846e+00
-9.1991e+00 -2.7421e+00
-9.0150e+00 -2.7926e+00
-8.8086e+00 -2.8359e+00
-8.5829e+00 -2.8720e+00
-8.3412e+00 -2.9015e+00
-8.0874e+00 -2.9252e+00
-7.8251e+00 -2.9437e+00
-7.5577e+00 -2.9581e+00
-7.2879e+00 -2.9691e+00
-7.0179e+00 -2.9773e+00
-6.7494e+00 -2.9835e+00
-6.4835e+00 -2.9881e+00
-6.2211e+00 -2.9915e+00
-5.9627e+00 -2.9940e+00
-5.7086e+00 -2.9957e+00
-5.4588e+00 -2.9970e+00
-5.2134e+00 -2.9979e+00
-4.9724e+00 -2.9986e+00
-4.7355e+00 -2.9991e+00
-4.5026e+00 -2.9994e+00
-4.2736e+00 -2.9996e+00
-4.0482e+00 -2.9997e+00
-3.8262e+00 -2.9998e+00
-3.6073e+00 -2.9999e+00
-3.3915e+00 -2.9999e+00
-3.1784e+00 -3.0000e+00
-2.9679e+00 -2.9831e+00
-2.7598e+00 -2.9190e+00
-2.5539e+00 -2.8112e+00
-2.3500e+00 -2.6655e+00
-2.1479e+00 -2.4888e+00
-1.9475e+00 -2.2888e+00
-1.7486e+00 -2.0740e+00
-1.5510e+00 -1.8527e+00
-1.3546e+00 -1.6337e+00
-1.1592e+00 -1.4249e+00
-9.6469e-01 -1.2341e+00
-7.7089e-01 -1.0680e+00
-5.7767e-01 -9.3254e-01
-3.8487e-01 -8.3230e-01
-1.9237e-01 -7.7075e-01
-1.8370e-15 -7.5000e-01
1.9237e-01 -7.7075e-01
3.8487e-01 -8.3230e-01
5.7767e-01 -9.3254e-01
7.7089e-01 -1.0680e+00
9.6469e-01 -1.2341e+00
1.1592e+00 -1.4249e+00
1.3546e+00 -1.6337e+00
1.5510e+00 -1.8527e+00
1.7486e+00 -2.0740e+00
1.9475e+00 -2.2888e+00
2.1479e+00 -2.4888e+00
2.3500e+00 -2.6655e+00
2.5539e+00 -2.8112e+00
2.7598e+00 -2.9190e+00
2.9679e+00 -2.9831e+00
3.1784e+00 -3.0000e+00
3.3915e+00 -2.9999e+00
3.6073e+00 -2.9999e+00
3.8262e+00 -2.9998e+00
4.0482e+00 -2.9997e+00
4.2736e+00 -2.9996e+00
4.5026e+00 -2.9994e+00
4.7355e+00 -2.9991e+00
4.9724e+00 -2.9986e+00
5.2134e+00 -2.9979e+00
5.4588e+00 -2.9970e+00
5.7086e+00 -2.9957e+00
5.9627e+00 -2.9940e+00
6.2211e+00 -2.9915e+00
6.4835e+00 -2.9881e+00
6.7494e+00 -2.9835e+00
7.0179e+00 -2.9773e+00
7.2879e+00 -2.9691e+00
7.5577e+00 -2.9581e+00
7.8251e+00 -2.9437e+00
8.0874e+00 -2.9252e+00
8.3412e+00 -2.9015e+00
8.5829e+00 -2.8720e+00
8.8086e+00 -2.8359e+00
9.0150e+00 -2.7926e+00
9.1991e+00 -2.7421e+00
9.3594e+00 -2.6846e+00
9.4956e+00 -2.6206e+00
9.6085e+00 -2.5512e+00
9.7001e+00 -2.4773e+00
9.7729e+00 -2.4002e+00
9.8298e+00 -2.3209e+00
9.8736e+00 -2.2404e+00
9.9070e+00 -2.1594e+00
9.9535e+00 -1.9839e+00
9.9776e+00 -1.8128e+00
9.9896e+00 -1.6476e+00
9.9954e+00 -1.4885e+00
9.9981e+00 -1.3354e+00
9.9992e+00 -1.1877e+00
9.9997e+00 -1.0447e+00
9.9999e+00 -9.0575e-01
1.0000e+01 -7.7027e-01
1.0000e+01 -6.3767e-01
1.0000e+01 -5.0741e-01
1.0000e+01 -3.7899e-01
1.0000e+01 -2.5192e-01
1.0000e+01 -1.2574e-01
1.0000e+01 0.0000e+00
};

\addplot [mark=none,red,line width=1.5] table{
4.4391e+00 1.5866e+00
4.4110e+00 1.5953e+00
4.3818e+00 1.6024e+00
4.3511e+00 1.6077e+00
4.3188e+00 1.6107e+00
4.2845e+00 1.6110e+00
4.2483e+00 1.6078e+00
4.2106e+00 1.6002e+00
4.1717e+00 1.5872e+00
4.1325e+00 1.5678e+00
4.0944e+00 1.5411e+00
4.0588e+00 1.5063e+00
4.0279e+00 1.4634e+00
4.0037e+00 1.4130e+00
3.9877e+00 1.3564e+00
3.9808e+00 1.2953e+00
3.9828e+00 1.2311e+00
3.9932e+00 1.1648e+00
4.0110e+00 1.0971e+00
4.0348e+00 1.0282e+00
4.0631e+00 9.5857e-01
4.0941e+00 8.8807e-01
4.1264e+00 8.1646e-01
4.1588e+00 7.4320e-01
4.1903e+00 6.6763e-01
4.2202e+00 5.8936e-01
4.2474e+00 5.0855e-01
4.2712e+00 4.2582e-01
4.2910e+00 3.4181e-01
4.3066e+00 2.5683e-01
4.3177e+00 1.7123e-01
4.3244e+00 8.5525e-02
4.3266e+00 -1.2310e-12
4.3244e+00 -8.5525e-02
4.3177e+00 -1.7123e-01
4.3066e+00 -2.5683e-01
4.2910e+00 -3.4181e-01
4.2712e+00 -4.2582e-01
4.2474e+00 -5.0855e-01
4.2202e+00 -5.8936e-01
4.1903e+00 -6.6763e-01
4.1588e+00 -7.4320e-01
4.1264e+00 -8.1646e-01
4.0941e+00 -8.8807e-01
4.0631e+00 -9.5857e-01
4.0348e+00 -1.0282e+00
4.0110e+00 -1.0971e+00
3.9932e+00 -1.1648e+00
3.9828e+00 -1.2311e+00
3.9808e+00 -1.2953e+00
3.9877e+00 -1.3564e+00
4.0037e+00 -1.4130e+00
4.0279e+00 -1.4634e+00
4.0588e+00 -1.5063e+00
4.0944e+00 -1.5411e+00
4.1325e+00 -1.5678e+00
4.1717e+00 -1.5872e+00
4.2106e+00 -1.6002e+00
4.2483e+00 -1.6078e+00
4.2845e+00 -1.6110e+00
4.3188e+00 -1.6107e+00
4.3511e+00 -1.6077e+00
4.3818e+00 -1.6024e+00
4.4110e+00 -1.5953e+00
4.4391e+00 -1.5866e+00
4.4668e+00 -1.5763e+00
4.4943e+00 -1.5644e+00
4.5222e+00 -1.5506e+00
4.5507e+00 -1.5349e+00
4.5800e+00 -1.5170e+00
4.6103e+00 -1.4967e+00
4.6417e+00 -1.4738e+00
4.6742e+00 -1.4481e+00
4.7078e+00 -1.4196e+00
4.7424e+00 -1.3882e+00
4.7782e+00 -1.3539e+00
4.8149e+00 -1.3166e+00
4.8526e+00 -1.2763e+00
4.8911e+00 -1.2331e+00
4.9304e+00 -1.1871e+00
4.9704e+00 -1.1383e+00
5.0108e+00 -1.0866e+00
5.0515e+00 -1.0322e+00
5.0922e+00 -9.7495e-01
5.1327e+00 -9.1496e-01
5.1725e+00 -8.5217e-01
5.2113e+00 -7.8659e-01
5.2487e+00 -7.1821e-01
5.2841e+00 -6.4708e-01
5.3170e+00 -5.7325e-01
5.3470e+00 -4.9686e-01
5.3735e+00 -4.1807e-01
5.3960e+00 -3.3716e-01
5.4140e+00 -2.5445e-01
5.4271e+00 -1.7038e-01
5.4351e+00 -8.5404e-02
5.4378e+00 6.6671e-13
5.4351e+00 8.5404e-02
5.4271e+00 1.7038e-01
5.4140e+00 2.5445e-01
5.3960e+00 3.3716e-01
5.3735e+00 4.1807e-01
5.3470e+00 4.9686e-01
5.3170e+00 5.7325e-01
5.2841e+00 6.4708e-01
5.2487e+00 7.1821e-01
5.2113e+00 7.8659e-01
5.1725e+00 8.5217e-01
5.1327e+00 9.1496e-01
5.0922e+00 9.7495e-01
5.0515e+00 1.0322e+00
5.0108e+00 1.0866e+00
4.9704e+00 1.1383e+00
4.9304e+00 1.1871e+00
4.8911e+00 1.2331e+00
4.8526e+00 1.2763e+00
4.8149e+00 1.3166e+00
4.7782e+00 1.3539e+00
4.7424e+00 1.3882e+00
4.7078e+00 1.4196e+00
4.6742e+00 1.4481e+00
4.6417e+00 1.4738e+00
4.6103e+00 1.4967e+00
4.5800e+00 1.5170e+00
4.5507e+00 1.5349e+00
4.5222e+00 1.5506e+00
4.4943e+00 1.5644e+00
4.4668e+00 1.5763e+00
4.4391e+00 1.5866e+00
};

\end{axis}

\end{tikzpicture}
\fi
\end{tabular}
\mcaption{A single vesicle discretized with $N=128$ points passing
through a constricted tube discretized with $N_{\wall}=256$ points.  The
boundary condition at the intake and outtake has a parabolic-profile and
on the rest of the solid wall is zero.}{f:stenosisGeom}
\end{figure}

In Figure~\ref{f:stenosisErrors}, we plot the errors in area and length
using $n_{\sdc}=1$ and a constant time step size.  Unsurprisingly, the
errors increase when the vesicle passes through the constriction.  In
this interval, a smaller time step should be taken.  In
Tables~\ref{t:adaptiveFirstOrderStenosis}--\ref{t:adaptiveThirdOrderStenosis},
we report results for adaptive time stepping with different tolerances
and different numbers of SDC corrections.  In
Figure~\ref{f:stenosisSummary}, we plot the errors in area and length
using $n_{\sdc}=1$ with and without adaptive time stepping.  We also
plot the time step size and the location of the rejected time steps.  As
expected, a much smaller time step size is taken as the vesicle passes
through the constriction.

We observe similar behavior as we observed for the {\em extensional}
example.  In particular,
\begin{itemize}
  \item To achieve three digits of accuracy with $n_{\sdc}=0$, a fixed
  time step requires more than 6,000 time steps, and with adaptive time
  steps, over 30,000 time steps are required.  This sharp increase is
  due to a very a small time step that must be taken to keep the local
  truncation error below the required threshold as the vesicle passes
  through the constriction (for this example, 65\% of the time steps
  are smaller than $10^{-4}$).  Using the second-order integrator
  $n_{\sdc}=1$ with adaptive time stepping, only $172 \times 3$
  adaptive time steps are required\footnote{Recall that $p-1=3$
  additional time steps are required because of the intermediate
  Gauss-Lobatto quadrature points.} are required to achieve 3 digits of
  accuracy.  The resulting speedup is a factor of greater than 22.


  \item Comparing the two adaptive time stepping integrators
  $n_{\sdc}=1$ and $n_{\sdc}=2$, six digits of accuracy can be computed
  with 16\% less CPU time by using $n_{\sdc}=2$.  Again, this indicates
  that if smaller tolerances are desired, additional SDC corrections
  should be used to increase the accuracy of each time step.
\end{itemize}

\begin{table}[htps]
\begin{center}
\begin{tabular}{c|cccc|cccc}
 & \multicolumn{4}{c|}{$\boldnsdc{0}$} & 
   \multicolumn{4}{c}{$\boldnsdc{1}$} \\ 
$m$ & $e_{A}$ & $e_{L}$ & \# fmm & CPU & 
      $e_{A}$ & $e_{L}$ & \# fmm & CPU \\
\hline
$750$  & $1.16$E$-2$ & $4.24$E$-2$ & $1.90$E$4$ & $1.0$ 
       & $6.05$E$-5$ & $2.57$E$-5$ & $1.26$E$5$ & $6.9$ \\  
$1500$ & $5.96$E$-3$ & $2.12$E$-2$ & $3.79$E$4$ & $2.2$  
       & $1.50$E$-5$ & $4.48$E$-6$ & $2.50$E$5$ & $13$ \\ 
$3000$ & $3.03$E$-3$ & $1.06$E$-2$ & $7.55$E$4$ & $4.2$
       & $3.73$E$-6$ & $7.00$E$-7$ & $4.96$E$5$ & $26$ \\ 
$6000$ & $1.53$E$-3$ & $5.28$E$-3$ & $1.51$E$5$ & $8.7$
       & $9.30$E$-7$ & $1.03$E$-7$ & $9.83$E$5$ & $53$  
\end{tabular}
\mcaption{The errors in area and length and the CPU time for a single
vesicle in a constricted tube ({\bf stenosis}) with a {\bf constant}
time step size using $\boldnsdc{0}$ (left) and $\boldnsdc{1}$ (right).
The CPU times for
Tables~\ref{t:SDC01stenosis}--\ref{t:adaptiveThirdOrderStenosis} are
relative to the cheapest simulation ($m=750$ and $n_{\sdc}=0$) which
took approximately $2.30$E$3$ seconds.  We achieve the expected first-
and second-order convergence rates.}{t:SDC01stenosis}

\begin{tabular}{c|cccc|cccc}
 & \multicolumn{4}{c|}{$\boldnsdc{2}$} & 
   \multicolumn{4}{c}{$\boldnsdc{3}$} \\ 
$m$ & $e_{A}$ & $e_{L}$ & \# fmm & CPU & 
      $e_{A}$ & $e_{L}$ & \# fmm & CPU \\
\hline
$750$  & $2.74$E$-5$ & $1.08$E$-6$  & $1.94$E$5$ & $11$ 
       & $3.51$E$-7$ & $4.38$E$-8$  & $2.62$E$5$ & $13$ \\ 
$1500$ & $7.16$E$-6$ & $7.36$E$-8$  & $3.84$E$5$ & $21$
       & $5.97$E$-8$ & $1.30$E$-9$  & $5.19$E$5$ & $27$ \\ 
$3000$ & $1.83$E$-6$ & $5.31$E$-9$  & $7.60$E$5$ & $39$
       & $9.30$E$-9$ & $4.60$E$-11$ & $1.03$E$6$ & $56$ \\ 
$6000$ & $4.63$E$-7$ & $4.53$E$-10$ & $1.51$E$6$ & $81$
       & $1.35$E$-9$ & $2.02$E$-12$ & $2.03$E$6$ & $108$  
\end{tabular}
\mcaption{The errors in area and length and the CPU time for a single
vesicle in a constricted tube ({\bf stenosis}) with a {\bf constant}
time step size using $\boldnsdc{2}$ (left) and $\boldnsdc{3}$ (right).
We do not obtain the expected convergence rates, but we do see that each
SDC correction does result in an increase in the accuracy.  We could
attempt to reach the asymptotic convergence rates, but this most likely
would require temporal resolutions where other sources of error will
dominate.}{t:SDC23stenosis}

\begin{tabular}{ccccccc}
Tolerance & $e_{A}$ & $e_{L}$ & Accepts & Rejects & \# fmm & CPU \\
\hline
$1$E$-2$ & $2.57$E$-3$ & $7.87$E$-3$ 
         & $3397$  & $21$ & $1.08$E$5$ & $5.7$ \\
$1$E$-3$ & $2.95$E$-4$ & $9.44$E$-4$ 
         & $33554$ & $12$ & $1.06$E$6$ & $59$ \\
\end{tabular}
\mcaption{The errors in area and length, the CPU time, and the number
of accepted and rejected time steps for a single vesicle in a
constricted tube ({\bf stenosis}) with an {\bf adaptive} time step size
using $\boldnsdc{0}$.  Since the time integrator is first-order, small
time steps have to be taken.  The result is that the CPU time is too
large at the reported tolerances.}{t:adaptiveFirstOrderStenosis}

\begin{tabular}{ccccccc}
Tolerance & $e_{A}$ & $e_{L}$ & Accepts & Rejects & \# fmm & CPU \\
\hline
$1$E$-2$ & $2.60$E$-3$ & $8.03$E$-4$ 
         & $68$   & $30$ & $2.06$E$4$ & $1.2$ \\
$1$E$-3$ & $5.55$E$-4$ & $1.11$E$-4$ 
         & $172$  & $68$ & $4.80$E$4$ & $2.6$ \\
$1$E$-4$ & $6.82$E$-5$ & $1.18$E$-5$ 
         & $503$  & $63$ & $1.11$E$5$ & $5.9$ \\
$1$E$-5$ & $7.20$E$-6$ & $9.08$E$-7$ 
         & $1559$ & $54$ & $3.15$E$5$ & $17$ \\
$1$E$-6$ & $7.35$E$-7$ & $5.60$E$-8$ 
         & $4904$ & $40$ & $9.61$E$5$ & $51$ \\
\end{tabular}
\mcaption{The errors in area and length, the CPU time, and the number
of accepted and rejected time steps for a single vesicle in a
constricted tube ({\bf stenosis}) with an {\bf adaptive} time step size
using $\boldnsdc{1}$.  This second-order method is able to take much
larger time steps than the first-order
method.}{t:adaptiveSecondOrderStenosis}

\begin{tabular}{ccccccc}
Tolerance & $e_{A}$ & $e_{L}$ & Accepts & Rejects & \# fmm & CPU \\
\hline
$1$E$-5$ & $7.19$E$-6$ & $4.63$E$-8$ 
         & $828$  & $136$ & $2.96$E$5$ & $16$ \\
$1$E$-6$ & $7.54$E$-7$ & $2.49$E$-9$ 
         & $2558$ & $105$ & $8.21$E$5$ & $43$ \\
\end{tabular}
\mcaption{The errors in area and length, the CPU time, and the number
of accepted and rejected time steps for a single vesicle in a
constricted tube ({\bf stenosis}) with an {\bf adaptive} time step size
using $\boldnsdc{2}$.  We see that if small tolerances are desired,
additional SDC iterations should be used to achieve a higher-order time
integrator.}{t:adaptiveThirdOrderStenosis}

\end{center}
\end{table}

\begin{figure}[htps]
\centering
\begin{tabular}{cccc}
\ifTikz
\input{stenosisErrorsSDC750.tikz} &
\input{stenosisErrorsSDC1500.tikz} &
\input{stenosisErrorsSDC3000.tikz} &
\input{stenosisErrorsSDC6000.tikz}
\fi
\end{tabular}
\mcaption{The errors in area and length for a single vesicle in a
constricted tube ({\bf stenosis}) with a {\bf constant} time step size
using $\boldnsdc{1}$ (right side of Table~\ref{t:SDC01stenosis}).
Notice that the error increases sharply as the vesicle passes through
the constriction for all the time step sizes.  This indicates that a
smaller time step size should be taken as the vesicle passes through
the constriction.}{f:stenosisErrors}
\end{figure}

\begin{figure}[htps]
\begin{tabular}{cc}
\ifTikz
\input{stenosisAdaptiveDT.tikz} &
\input{stenosisAdaptiveErrors.tikz}
\fi
\end{tabular}

\mcaption{Results for the {\bf stenosis} flow using $\boldnsdc{1}$.
Left: The time step size using 6000 {\bf constant} time step sizes
(dashed) and an {\bf adaptive} time step size (solid).  The open circles
indicate the 40 times when the time step size is rejected.  Notice that
the time step size decreases as the vesicle passes through the
constriction.  Right: The errors in area and length using {\bf constant}
(dashed) and {\bf adaptive} (solid) time steps, and the desired error
(black) of the adaptive time step.  When using adaptive time stepping,
the error in area nearly achieves the desired error everywhere except
when the vesicle passes through the constriction, where the error in
area actually drops.  However, when using a constant time step size,
very little error is committed at the start of the simulation, and then
the majority of the error is committed as the vesicle passes through the
constriction.  Even though the CPU savings are negligible (it is about
4\%), the adaptive time stepping method did not require any trial and
error to find an appropriate time step size.  The user simply specifies
that they desire six digits accuracy and it is achieved.  In addition,
by using an additional SDC correction, the CPU savings is increased to
19\%.}{f:stenosisSummary}
\end{figure}

\subsection{Couette}
\label{s:couette}
We consider eight vesicles discretized with $N=192$ points in a couette
apparatus where each wall is discretized with $N_{\mathrm{wall}}=512$
points (Figure~\ref{f:couetteGeom}).  At this resolution, our FMM
implementation of the double-layer potential is faster than a direct
evaluation.  Therefore, we use the FMM to evaluate all the single- and
double-layer potentials.  We use $p=4$ Gauss-Lobatto quadrature points
and take the long time horizon $T=50$.  At the time horizon, the inner
boundary has made five complete rotations with a constant velocity, and
the outer boundary is stationary.  The inner and outer boundaries are
aligned so that there is a narrow region that the vesicles must pass
through.  As a vesicle passes through this region, a smaller time step
should be taken.  We check the convergence rates with a constant time
step size using zero and one SDC correction in
Tables~\ref{t:SDC0couette} and~\ref{t:SDC1couette}.  For both time
integrators, the expected first- and second-order convergence rates are
achieved.

\begin{figure}[htps]
\centering
\begin{tabular}{cccccc}
\ifTikz
\input{couetteSnapsTime1.tikz} &
\input{couetteSnapsTime2.tikz} &
\input{couetteSnapsTime3.tikz} &
\input{couetteSnapsTime4.tikz} &
\input{couetteSnapsTime5.tikz} &
\input{couetteSnapsTime6.tikz} \\
\input{couetteSnapsTime7.tikz} &
\input{couetteSnapsTime8.tikz} &
\input{couetteSnapsTime9.tikz} &
\input{couetteSnapsTime10.tikz} &
\input{couetteSnapsTime11.tikz} &
\input{couetteSnapsTime12.tikz}
\fi
\end{tabular}
\mcaption{Eight vesicles discretized with $N=192$ points in a couette
apparatus whose solid walls are each discretized with $N_{\wall}=512$
points.  The outer boundary is stationary and the inner boundary has
constant angular velocity and has completed five full rotations at
$T=50$.  A single vesicle is colored in blue to help with visualization.
The interactions between the vesicles and the solid walls complicate and
simplify throughout the simulation and the simulation benefits from
adaptive time stepping.}{f:couetteGeom}
\end{figure}

In Figure~\ref{f:couetteErrors}, we plot the errors in area and length
using $n_{\sdc}=1$ with a constant time step size.  If $500$ time steps
are taken, the error in length has a sharp jump caused by the blue
vesicle in Figure~\ref{f:couetteGeom}.  By reducing the time step size,
the errors decrease at the expected rate, but there is still a jump in
the errors near $t=2$.  In Table~\ref{t:adaptiveSecondOrderCouette}, we
report results for $n_{\sdc}=1$ with adaptive time step sizes and two
different tolerances.  In Figure~\ref{f:couetteSummary}, we plot the
errors in area and length with and without adaptive time stepping.  We
also plot the time step size and the locations of rejected time steps.
We see that a small time step is taken initially as the vesicles smooth,
and then much larger time steps can be taken once the vesicles are
smooth and separated from one another and the solid walls.

While the desired error is achieved, the final error is much smaller
than the tolerance when compared to the {\em extensional} and {\em
stenosis} examples.  This is a result of having more vesicles since the
local truncation error is estimated with the maximum change in error,
where the maximum is taken over all the vesicles.  In contrast to the
previous two examples, different vesicles commit the largest error at
different times.  The result is that the global error of our adaptive
time stepping method is much less than the desired tolerance.  A
possible solution to this problem is to adjust the allowable truncation
error at each time step.  In particular, if we again assume that the
time horizon is $T=1$ and the desired tolerance is $\epsilon$, then
condition~\eqref{e:errorBounds} is replaced with
\begin{align*}
  |A(t + \Delta t) - A(t)| \leq A(t) \frac{\Delta t}{1-t} 
    \left(\epsilon - \frac{|A(t) - A(0)|}{A(0)} \right),
\end{align*}
and a similar condition for the length.  Using this new strategy with
the time integrator from Table~\ref{t:adaptiveSecondOrderCouette}, with
a tolerance of $1$E$-1$, the final error is $5.14$E$-2$, the number of
accepted time steps is reduced to 222, but 84 time steps are rejected.
Using a tolerance of $1$E$-2$, the final error is $9.62$E$-3$ with 583
accepted time steps, and the number of rejected time steps is only
slightly increased to 74.  While this strategy does offer a speedup, we
expect that multiscale time integrators, where each vesicle's time step
size can differ from the others', will overall result in a more robust
solver.

We have also run the same simulation with a coarser discretization.  We
leave all the parameters unchanged with the exception of $N=92$ points
per vesicle and $N_{\wall}=192$ points per solid wall.   Using adaptive
time stepping with one SDC correction, a tolerance of $1$E$-1$ and
$5$E$-2$ are attainable.  However, with a tolerance of $1$E$-2$, the
error condition~\eqref{e:errorBounds} can not be satisfied at the
initial condition for all $\Delta t$.  Therefore, the spatial error is
too large for condition~\eqref{e:errorBounds} to be satisfied.  We also
ran this resolution at a later time when the vesicle shapes are
smoothed.  There, the desired local truncation error is attainable and
the simulation can successfully complete.  This indicates that spatial
adaptivity is required, and this is a current extension we are
developing.

Our conclusions from this example are:
\begin{itemize}

\item Considering the time integrator $n_{\sdc}=1$, to achieve one
digit of accuracy, adaptive time stepping is at least 33\% faster.
However, to achieve two digits accuracy, with our current limitations
(fixed spatial resolution, and each vesicle having the same time step
size), adaptive time stepping does not conclusively offer a speedup.
However, unlike the simulations with fixed time step size, no trial and
error procedure to choose a time step size that results in the desired
global error is required.

\item With multiple vesicles, it is harder to achieve the desired
global error.  By adjusting the allowable local truncation error
throughout the simulation, the final error is much closer to the
tolerance.  However, each vesicle must still have the same time step
size.  We plan to develop integrators that use different time step
sizes for each vesicle.

\item Without spatial adaptivity, excessively large resolutions are
required for the entirety of the simulation.

\end{itemize}

\begin{table}[htps]
\begin{center}

\begin{tabular}{c|cccc}
 \multicolumn{4}{c}{$\boldnsdc{0}$} \\ 
$m$ & $e_{A}$ & $e_{L}$ & \# fmm & CPU \\ 
\hline
2000 & $1.24$E$-1$ & $3.50$E$-1$ & $2.73$E$5$ & $1.0$ \\
4000 & $2.04$E$-2$ & $1.21$E$-1$ & $5.28$E$5$ & $2.1$ \\
8000 & $5.99$E$-3$ & $5.79$E$-2$ & $1.05$E$6$ & $3.8$
\end{tabular}
\mcaption{The errors in area and length and the CPU time for eight
vesicles in a {\bf couette} apparatus with a {\bf constant} time step
size using $\boldnsdc{0}$.  The CPU for
Tables~\ref{t:SDC0couette}--\ref{t:adaptiveSecondOrderCouette} are
relative to the cheapest simulation ($m=2000$ and $n_{\sdc}=0$) which
took approximately $1.44$E$5$ seconds.  We achieve the expected
first-order convergence.}{t:SDC0couette}

\begin{tabular}{c|cccc}
 & \multicolumn{4}{c}{$\boldnsdc{1}$} \\ 
$m$ & $e_{A}$ & $e_{L}$ & \# fmm & CPU \\
\hline
$500$  & $7.14$E$-3$ & $5.39$E$-2$ & $4.27$E$5$ & $1.5$ \\
$1000$ & $2.79$E$-3$ & $4.68$E$-5$ & $8.28$E$5$ & $2.9$ \\
$2000$ & $7.04$E$-4$ & $9.34$E$-6$ & $1.69$E$6$ & $6.0$  
\end{tabular}
\mcaption{The errors in area and length and the CPU time for eight
vesicles in a {\bf couette} apparatus with a {\bf constant} time step
size using $\boldnsdc{1}$.  At the reported resolutions, we are just
beginning to achieve second-order convergence.}{t:SDC1couette}

\begin{tabular}{ccccccc}
Tolerance & $e_{A}$ & $e_{L}$ & Accepts & Rejects & \# fmm & CPU \\
\hline
$1$E$-1$ & $2.65$E$-2$ & $2.07$E$-2$ 
         & $256$ & $68$ & $2.90$E$5$ & $1.0$ \\
$1$E$-2$ & $4.21$E$-3$ & $3.86$E$-4$ 
         & $780$ & $65$ & $6.97$E$5$ & $2.5$
\end{tabular}
\mcaption{The errors in area and length, the CPU time, and the number of
accepted and rejected time steps for eight vesicles in a {\bf couette}
apparatus with an {\bf adaptive} time step size using
$\boldnsdc{1}$.}{t:adaptiveSecondOrderCouette}
\end{center}
\end{table}

\begin{figure}[htps]
\centering
\begin{tabular}{ccc}
\ifTikz
\input{couetteErrorsSDC500.tikz} &
\input{couetteErrorsSDC1000.tikz} &
\input{couetteErrorsSDC2000.tikz}
\fi
\end{tabular}
\mcaption{The errors in area and length for eight vesicles in a {\bf
couette} apparatus with a {\bf constant} time step size using
$\boldnsdc{1}$.  For all the time step sizes, there is a sharp increase
in the errors around $t=2$ which is committed by the vesicle colored in
blue in Figure~\ref{f:couetteGeom}.  Adaptive time stepping
automatically resolves these sharp jumps and avoids the need to use
trial and error to find an appropriate time step size.}{f:couetteErrors}
\end{figure}

\begin{figure}[htps]
\begin{tabular}{cc}
\ifTikz
\input{couetteAdaptiveDT.tikz} &
\input{couetteAdaptiveErrors.tikz}
\fi
\end{tabular}
\mcaption{Results for the {\bf couette} apparatus using
$\boldnsdc{1}$.  Left: The time step size using a {\bf constant} time
step (dashed) and an {\bf adaptive} time step size (solid).  The open
circles indicate the 65 times when the time step size is rejected.  We
see that the time step size is smallest when the error is large (near
$t=2$) and is largest when the vesicles and walls are well-separated
(near $t=15$ and $t=30$).  Right: The errors in area and length using
{\bf constant} (dashed) and {\bf adaptive} (solid) time steps, and the
desired error (black) of the adaptive time step.  Since each vesicle
must use the same time step size, the actual error does not closely
follow the desired error.}{f:couetteSummary}
\end{figure}

\section{Conclusions\label{s:conclusions}} 
We have developed and tested adaptive time stepping methods for an
integro-differential equation used to simulate vesicle suspensions.  To
select the time step size, we followed the standard strategy of
estimating the local truncation error, and then adjusting the time step
size so that the same amount of error is committed per time step.  This
strategy does not require users to use a trial and error procedure to
find a time step size that achieves the desirable error.  Because of
the cost of forming multiple numerical solutions, we introduced a
technique where the error in area and length of the vesicles, which are
fixed by the physics, are used to estimate the local truncation error.
This method proved reliable, but it can only be applied to problems
with measurable invariants.

In order to allow for larger time steps and longer time horizons, we
introduced new high-order time integrators based on spectral deferred
corrections (SDC).  An advantage of these integrators is that they only
require the solution from the previous time step making them compatible
with adaptive time step sizes.  To the best of our knowledge, this work
is the first to apply IMEX methods coupled with SDC to a
two-dimensional integro-differential equation.  Therefore, we currently
do not have any asymptotic results and our convergence results are
entirely based on numerical results.

These contributions are a major step towards the development of a robust
solver for vesicle suspensions, but several other features need to be
introduced.
\begin{itemize}
  \item Theoretical results, in particular the convergence rates of SDC,
  need to be developed.  This will allow us to take the largest possible
  time step size given the user specified tolerance.

  \item Our SDC formulation does not support viscosity contrast.  The
  challenge is that vesicle suspensions with viscosity contrast
  are governed by $\dot{\xx} = F(\xx,\dot{\xx})$ which makes
  formulating a Picard integral difficult. 

  \item Many suspensions would benefit from spatial adaptivity.  This
  was demonstrated in the {\em couette} example where we saw that our
  adaptive time stepping method requires a sufficiently resolved spatial
  discretization, but this discretization was not required for the
  entirety of the simulation.

  \item All of our results are two-dimensional.  However, none of the
  methods that we have introduced in this paper are restricted to two
  dimensions.  Performing three-dimensional simulations only requires a
  careful implementation of the methods and the development of suitable
  preconditioners.

\end{itemize}

\begin{appendices}

\section{Complete SDC Formulation \label{a:appendix1}}
Here we write our SDC formulation for the unbounded vesicle formulation
that includes the tension $\sigma_{j}$.  We also do not use the notation
$\vv(\xx_{j};\xx_{k})$ to help clarrify the discretization of the
prediction and SDC correction steps.  Let $\txx{}$ and $\tsigma{}$ be a
provisional solution of
\begin{align}
  \xx_{j}(t) = \xx_{0} + \int_{0}^{t} \left(\vv_{\infty}(\xx_{j}) + 
    \sum_{k=1}^{M}\SS(\xx_{j},\xx_{k})(-\BB(\xx_{k})(\xx_{k}) + 
    \TT(\xx_{k})(\sigma_{k}))\right)d\tau,
  \label{e:picardFull}
\end{align}
and then compute the residual
\begin{align*}
  \rr_{j}(t;\txx{},\tsigma{}) = \xx_{0} - \txx{j} + 
    \int_{0}^{t} \left(\vv_{\infty}(\txx{j}) + 
    \sum_{k=1}^{M}\SS(\txx{j},\txx{k})(-\BB(\txx{k})(\txx{k}) + 
    \TT(\txx{k})(\tsigma{k}))\right)d\tau.
\end{align*}
The error in the position is $\exx{j} = \xx_{j} - \txx{j}$, and the
error in tension is $\esigma{j} = \sigma_{j} - \tsigma{j}$.
Substituting the errors into~\eqref{e:picardFull},
\begin{equation}
\begin{aligned}
  \exx{j}(t) = \rr_{j}(t;\txx{},\tsigma{}) &+
    \int_{0}^{t}(\vv_{\infty}(\txx{j}+\exx{j}) -
      \vv_{\infty}(\txx{j}))d\tau \\
    &-\int_{0}^{t}\sum_{k=1}^{M}
    \left(\SS(\txx{j}+\exx{j},\txx{k}+\exx{k})
      \BB(\txx{k}+\exx{k}) - \SS(\txx{j},\txx{k})\BB(\txx{k})
        \right)(\txx{k}) d\tau \\
    &+\int_{0}^{t}\sum_{k=1}^{M}
      \left(\SS(\txx{j}+\exx{j},\txx{k}+\exx{k})
      \TT(\txx{k}+\exx{k}) - \SS(\txx{j},\txx{k})\TT(\txx{k})
        \right)(\tsigma{k}) d\tau \\
    &+\int_{0}^{t}\sum_{k=1}^{M} \SS(\txx{j}+\exx{j},\txx{k}+\exx{k})
      (-\BB(\txx{k}+\exx{k})(\exx{k}) +
        \TT(\txx{k}+\exx{k})(\esigma{k}))d\tau.
\end{aligned}
\label{e:sdcUpdateFull}
\end{equation}
For confined flows, $\vv_{\infty}$ depends on yet another variable, a
density function defined on the solid walls $\Gamma$, which we denote
by $\eeta$.  The dependence is through a double-layer potential
$\DD(\txx{j},\Gamma)$.  In this case, the term involving $\vv_{\infty}$
in~\eqref{e:sdcUpdateFull} becomes
\begin{align*}
  \int_{0}^{t}(\DD(\txx{j}+\exx{j},\Gamma) - \DD(\txx{j},\Gamma))
      (\eeta) d\tau,
\end{align*}
and the final term in~\eqref{e:sdcUpdateFull} includes the term
$\DD(\txx{j}+\exx{j},\Gamma)(\ee_{\eeta_{k}})$.

The inextensibility constraint~\eqref{e:inextens} is
\begin{equation}
\begin{aligned}
  0&=\Div(\xx_{j})\left(\vv_{\infty}(\xx_{j})
    +\sum_{k=1}^{M} \SS(\xx_{j},\xx_{k})(
    -\BB(\xx_{k})(\xx_{k}) + \TT(\xx_{k})(\sigma_{k}))\right) \\
  &=\frac{d\xx_{j}}{ds} \cdot \frac{d}{ds}\left(\vv_{\infty}(\xx_{j})
    +\sum_{k=1}^{M} \SS(\xx_{j},\xx_{k})(
    -\BB(\xx_{k})(\xx_{k}) + \TT(\xx_{k})(\sigma_{k}))\right),
\end{aligned}
\label{e:inextensMethod1}
\end{equation}
and the errors satisfy
\begin{align}
  0=&\Div(\txx{j}+\exx{j})\left(\vphantom{\sum_{k=1}^{M}} 
    \vv_{\infty}(\xx_{j}) \right. \nonumber \\
  &\left. + \sum_{k=1}^{M} \SS(\txx{j}+\exx{j},\txx{k}+\exx{k})\left(
    -\BB(\txx{k}+\exx{k})(\txx{k}+\exx{k}) + 
    \TT(\txx{k}+\exx{k})(\tsigma{k}+\esigma{k})\right)\right).
  \label{e:inextensUpdateMethod1}
\end{align}

\subsection{Discretization}
The provisional solution is found by
discretizing~\eqref{e:picardFull} and~\eqref{e:inextensMethod1} as
\begin{align*}
  \xx_{j}^{n+1} = \xx_{j}^{n} + \Delta t_{n}
    \left(\vv_{\infty}(\xx^{n}_{j}) + 
    \sum_{k=1}^{M} \SS(\xx^{n}_{j},\xx^{n}_{k})(
    -\BB(\xx^{n}_{k})(\xx^{n+1}_{k}) +
    \TT(\xx^{n}_{k})(\sigma^{n+1}_{k}))\right),
\end{align*}
and
\begin{align*}
  \Div(\xx^{n}_{j})\left(\sum_{k=1}^{M} \SS(\xx^{n}_{j},\xx^{n}_{k})(
    -\BB(\xx^{n}_{k})(\xx^{n+1}_{k}) +
    \TT(\xx^{n}_{k})(\sigma^{n+1}_{k}))\right)=0.
\end{align*}
The SDC updates~\eqref{e:sdcUpdateFull}
and~\eqref{e:inextensUpdateMethod1} are discretized as
\begin{align*}
  \ee^{n+1}_{\xx_{j}} = \ee^{n}_{\xx_{j}} + 
  \rr_{j}^{n+1} - \rr_{j}^{n}
  +\Delta t_{n} \left(\sum_{k=1}^{M} 
    \SS(\txx{j}^{n+1},\txx{k}^{n+1})\left(
    -\BB(\txx{k}^{n+1})(\ee^{n+1}_{\xx_{k}}) +
    \TT(\txx{k}^{n+1})(\delta^{n+1}_{\sigma_{k}})
    \right)\right),
\end{align*}
and
\begin{align*}
  \Div(\txx{j}^{n+1})&\left(\sum_{k=1}^{M} 
    \SS(\txx{j}^{n+1},\txx{k}^{n+1})\left(
    -\BB(\txx{k}^{n+1})(\ee^{n+1}_{\xx_{k}}) +
    \TT(\txx{k}^{n+1})(\delta^{n+1}_{\sigma_{k}})
    \right)\right) \\
  &=-\Div(\txx{j}^{n+1})\left(\vv_{\infty}(\txx{j}^{n+1})
    +\sum_{k=1}^{M} \SS(\txx{j}^{n+1},\txx{k}^{n+1})\left(
    -\BB(\txx{k}^{n+1})(\txx{k}^{n+1}) + 
    \TT(\txx{k}^{n+1})(\tsigma{k}^{n+1})\right)\right).
\end{align*}
With this discretization of~\eqref{e:inextensUpdateMethod1}, we see that
if $\ee^{n}_{\xx_{k}}=0$ and $\delta^{n}_{\sigma_{k}}=0$, then
\begin{align*}
  \Div(\txx{j}^{n+1})\left(\vv_{\infty}(\txx{j}^{n+1})
    +\sum_{k=1}^{M} \SS(\txx{j}^{n+1},\txx{k}^{n+1})\left(
    -\BB(\txx{k}^{n+1})(\txx{k}^{n+1}) + 
    \TT(\txx{k}^{n+1})(\tsigma{k}^{n+1})\right)\right) = 0,
\end{align*}
which means that
\begin{align*}
  \Div(\txx{j}^{n})\left( \sum_{k=1}^{M}
    \vv(\txx{j}^{n};\txx{k}^{n})\right) = 0.
\end{align*}
In other words, by using the
discretization~\eqref{e:numericSDCImplicit}, the fixed point of the SDC
iteration exactly satisfies the inextensibility condition.

\end{appendices}

\bibliographystyle{plainnat} 
\biboptions{sort&compress}
\bibliography{refs}

\begin{thebibliography}{46}
\providecommand{\natexlab}[1]{#1}
\providecommand{\url}[1]{\texttt{#1}}
\expandafter\ifx\csname urlstyle\endcsname\relax
  \providecommand{\doi}[1]{doi: #1}\else
  \providecommand{\doi}{doi: \begingroup \urlstyle{rm}\Url}\fi

\bibitem[Alpert(1999)]{alp1999}
B.~K. Alpert.
\newblock {Hybrid Gauss-Trapezoidal Quadrature Rules}.
\newblock \emph{SIAM Journal on Scientific Computing}, 20:\penalty0 1551--1584,
  1999.

\bibitem[Ascher et~al.(1995)Ascher, Ruuth, and Wetton]{asc:ruu:wet1995}
U.M. Ascher, S.J. Ruuth, and B.T.R. Wetton.
\newblock Implicit-explicit methods for time-dependent partial differential
  equations.
\newblock \emph{SIAM Journal on Numerical Analysis}, 32:\penalty0 797--823,
  1995.

\bibitem[Biben et~al.(2005)Biben, Kassner, and Misbah]{bib:kas:mis2005}
Thierry Biben, Klaus Kassner, and Chaouqi Misbah.
\newblock Phase-field approach to three-dimensional vesicle dynamics.
\newblock \emph{Physical Review E}, 72\penalty0 (041921), 2005.

\bibitem[B\"{o}hmer and Stetter(1984)]{boh:ste1984}
K.~B\"{o}hmer and H.~J. Stetter.
\newblock \emph{{Defect Correction Methods, Theory and Applications}}.
\newblock Springer-Verlag, New York, 1984.

\bibitem[Bourlioux et~al.(2003)Bourlioux, Layton, and Minion]{bou:lay:min2003}
Anne Bourlioux, Anita~T. Layton, and Michael~L. Minion.
\newblock {High-Order Multi-Implicit Spectral Deferred Correction Methods for
  Problems of Reactive Flow}.
\newblock \emph{Journal of Computational Physics}, 189:\penalty0 651--675,
  2003.

\bibitem[Bouzarth and Minion(2010)]{bou:min2010}
Elizabeth~L. Bouzarth and Michael~L. Minion.
\newblock A multirate time integrator for regularized stokeslets.
\newblock \emph{Journal of Computational Physics}, 229\penalty0 (11):\penalty0
  4208--4224, 2010.

\bibitem[Bu et~al.(2012)Bu, Huang, and Minion]{bu:hua:min2012}
Sunyoung Bu, Jingfang Huang, and Michael~L. Minion.
\newblock {Semi-implicit Krylov Deferred Correction Methods for Differential
  Algebraic Equations}.
\newblock \emph{Mathematics of Computation}, 81:\penalty0 2127--2157, 2012.

\bibitem[Christlieb and Ong(2010)]{chr:ong2010}
Andrew Christlieb and Benjamin Ong.
\newblock {Implicit Parallel Time Integrators}.
\newblock \emph{Journal on Scientific Computing}, 49:\penalty0 167--179, 2010.

\bibitem[Christlieb et~al.(2009)Christlieb, Ong, and Qui]{chr:ong:qui2009}
Andrew Christlieb, Benjamin Ong, and Jing-Mei Qui.
\newblock Comments on high order integrators embedded within integral deferred
  correction methods.
\newblock \emph{Communications in Applied Mathematics and Computational
  Science}, 4:\penalty0 27--56, 2009.

\bibitem[Christlieb et~al.(2010)Christlieb, Macdonald, and
  Ong]{chr:mac:ong2010}
Andrew~J. Christlieb, Colin~B. Macdonald, and Benjamin~W. Ong.
\newblock Parallel high-order integrators.
\newblock \emph{SIAM Journal on Scientific Computing}, 32\penalty0
  (2):\penalty0 818--835, 2010.

\bibitem[Christlieb et~al.(2012)Christlieb, Haynes, and Ong]{chr:hay:ong2012}
Andrew~J. Christlieb, Ronald~D. Haynes, and Benjamin~W. Ong.
\newblock {A Parallel Space-Time Algorithm}.
\newblock \emph{SIAM Journal on Scientific Computing}, 34\penalty0
  (5):\penalty0 233--248, 2012.

\bibitem[Du and Zhang(2007)]{du:zha2007}
Qiang Du and Jian Zhang.
\newblock Adaptive finite element method for a phase field bending elasticity
  model of vesicle membrane deformations.
\newblock \emph{SIAM Journal On Scientific Computing}, 30\penalty0
  (3):\penalty0 1634--1657, 2007.

\bibitem[Dutt et~al.(2000)Dutt, Greengard, and Rokhlin]{dut:gre:rok2000}
A.~Dutt, L.~Greengard, and V.~Rokhlin.
\newblock Spectral {D}eferred {C}orrection {M}ethods for {O}rdinary
  {D}ifferential {E}quations.
\newblock \emph{BIT Numerical Mathematics}, 40:\penalty0 241--266, 2000.

\bibitem[Germund G.~Dahlquist and Nevanlinna(1983)]{dah:lin:nev1983}
Werner~Liniger Germund G.~Dahlquist and Olavi Nevanlinna.
\newblock {Stability of Two-Step Methods for Variable Integration Steps}.
\newblock \emph{SIAM Journal on Numerical Analysis}, 20\penalty0 (5):\penalty0
  1071--1--85, October 1983.

\bibitem[Ghigliotti et~al.(2011)Ghigliotti, Rahimian, Biros, and
  Misbah]{ghi:rah:bir:mis2011}
G.~Ghigliotti, A.~Rahimian, G.~Biros, and C.~Misbah.
\newblock {Vesicle Migration and Spatial Organization Driven by Flow Line
  Curvature}.
\newblock \emph{Physical Review Letters}, 106:\penalty0 028101, 2011.

\bibitem[Hairer et~al.(1993)Hairer, Wanner, and N$\o$rsett]{hai:nor:wan1993}
Ernest Hairer, Gerhard Wanner, and Syvert~Paul N$\o$rsett.
\newblock \emph{{Solving Ordinary Differential Equations I: Nonstiff
  Problems}}.
\newblock Springer, 1993.

\bibitem[Hansen and Strain(2011)]{han:str2011}
Anders~C. Hansen and John Strain.
\newblock On the order of deferred correction.
\newblock \emph{Applied Numerical Mathematics}, 61:\penalty0 961--973, 2011.

\bibitem[Huang et~al.(2006{\natexlab{a}})Huang, Jia, and
  Minion]{hua:jia:min2006}
Jingfang Huang, Jun Jia, and Michael Minion.
\newblock Accelerating the convergence of spectral deferred correction methods.
\newblock \emph{Journal of Computational Physics}, 214:\penalty0 633--656,
  2006{\natexlab{a}}.

\bibitem[Huang et~al.(2006{\natexlab{b}})Huang, Lai, and
  Xiang]{hua:lai:xia2006}
Jingfang Huang, Ming-Chih Lai, and Yang Xiang.
\newblock {An Integral Equation Method for Epitaxial Step-flow Growth
  Simulations}.
\newblock \emph{Journal of Computational Physics}, 216:\penalty0 724--743,
  2006{\natexlab{b}}.

\bibitem[Huang et~al.(2007)Huang, Jia, and Minion]{hua:jia:min2007}
Jingfang Huang, Jun Jia, and Michael Minion.
\newblock Arbitrary order {K}rylov deferred correction methods for differential
  algebraic equations.
\newblock \emph{Journal of Computational Physics}, 221\penalty0 (2):\penalty0
  739--760, 2007.

\bibitem[Jia and Huang(2008)]{jia:hua2008}
Jun Jia and Jingfang Huang.
\newblock {Krylov deferred correction accelerated method of lines transpose for
  parabolic equations}.
\newblock \emph{Journal of Computational Physics}, 227:\penalty0 1739--1753,
  2008.

\bibitem[Kaoui et~al.(2011)Kaoui, Tahiri, Biben, Ez-Zahraouy, Benyoussef,
  Biros, and Misbah]{kao:tah:bib:ezz:ben:bir:mis2011}
B.~Kaoui, N.~Tahiri, T.~Biben, H.~Ez-Zahraouy, A.~Benyoussef, G.~Biros, and
  C.~Misbah.
\newblock {What Dictates Red Blood Cell Shapes and Dynamics in the
  Microvasculature?}
\newblock \emph{Physical Review E}, pages 1--11, 2011.

\bibitem[Kim and Lai(2010)]{kim:lai2010}
Yongsam Kim and Ming-Chih Lai.
\newblock Simulating the dynamics of inextensible vesicles by the penalty
  immersed boundary method.
\newblock \emph{Journal of Computational Physics}, 229\penalty0 (12):\penalty0
  4840--4853, 2010.

\bibitem[Laadhari et~al.(2014)Laadhari, Saramito, and Misbah]{laa:sar:mis2014}
Aymen Laadhari, Pierre Saramito, and Chaouqi Misbah.
\newblock Computing the dynamics of biomembranes by combining conservative
  level set and adaptive finite element methods.
\newblock \emph{Journal of Computational Physics}, 263:\penalty0 328--352,
  2014.

\bibitem[Lindberg(1980)]{lin1980}
Bengt Lindberg.
\newblock Error estimation and iterative improvement for discretization
  algorithms.
\newblock \emph{BIT}, 20:\penalty0 486--500, 1980.

\bibitem[Minion(2003)]{min2003}
Michael~L. Minion.
\newblock {Semi-Implicit Spectral Deferred Correction Methods for Ordinary
  Differential Equations}.
\newblock \emph{Communications in Mathematical Sciences}, 1\penalty0
  (3):\penalty0 471--500, 2003.

\bibitem[Misbah(2006)]{mis2006}
C.~Misbah.
\newblock {Vacillating Breathing and Tumbling of Vesicles under Shear Flow}.
\newblock \emph{Physical Review Letters}, 96\penalty0 (2), 2006.

\bibitem[Noguchi and Gompper(2005)]{nog:gom2005}
H.~Noguchi and D.~G. Gompper.
\newblock Shape transitions of fluid vesicles and red blood cells in capillary
  flows.
\newblock \emph{Proceedings Of The National Academy Of Sciences Of The United
  States Of America}, 102:\penalty0 14159--14164, 2005.

\bibitem[Ong et~al.(2012)Ong, Melfi, and Christlieb]{ong:mel:chr2012}
Benjamin Ong, Andrew Melfi, and Andrew Christlieb.
\newblock {Parallel Semi-Implicit Time Integrators}.
\newblock \emph{arxiv}, 2012.
\newblock 1209.4297v1.

\bibitem[Pozrikidis(1990)]{poz1990}
C.~Pozrikidis.
\newblock {The Axisymmetric Deformation Of A Red Blood Cell In Uniaxial
  Straining Stokes Flow}.
\newblock \emph{Journal of Fluid Mechanics}, 216:\penalty0 231--254, 1990.

\bibitem[Quaife and Biros(2013)]{qua:bir2013b}
Bryan Quaife and George Biros.
\newblock {High-volume fraction simulations of two-dimensional vesicle
  suspensions}.
\newblock \emph{arxiv}, 2013.
\newblock 1309.1128.

\bibitem[Rahimian et~al.(2010{\natexlab{a}})Rahimian, Lashuk, Veerapaneni,
  Chandramowlishwaran, Malhotra, Moon, Sampath, Shringarpure, Vetter, Vuduc,
  Zorin, and Biros]{rah:las:vee:cha:mal:moo:sam:shr:vet:vud:zor:bir2010}
Abtin Rahimian, Ilya Lashuk, Shravan~K. Veerapaneni, Aparna
  Chandramowlishwaran, Dhairya Malhotra, Logan Moon, Rahul Sampath, Aashay
  Shringarpure, Jeffrey Vetter, Richard Vuduc, Denis Zorin, and George Biros.
\newblock Petascale direct numerical simulation of blood flow on 200k cores and
  heterogeneous architectures.
\newblock In \emph{SC '10: Proceedings of the 2010 ACM/IEEE conference on
  Supercomputing}, pages 1--12, Piscataway, NJ, USA, 2010{\natexlab{a}}. IEEE
  Press.

\bibitem[Rahimian et~al.(2010{\natexlab{b}})Rahimian, Veerapaneni, and
  Biros]{rah:vee:bir2010}
Abtin Rahimian, Shravan~Kumar Veerapaneni, and George Biros.
\newblock Dynamic simulation of locally inextensible vesicles suspended in an
  arbitrary two-dimensional domain, a boundary integral method.
\newblock \emph{Journal of Computational Physics}, 229:\penalty0 6466--6484,
  2010{\natexlab{b}}.

\bibitem[Rahimian et~al.(2012)Rahimian, Veerapaneni, Zorin, and
  Biros]{rah:vee:zor:bir2012}
Abtin Rahimian, Shravan~K. Veerapaneni, Denis Zorin, and George Biros.
\newblock {Fully discrete Galerkin Spectral Boundary Integral Method for
  Vesicle Flows with Viscosity Contrast}.
\newblock pages 1--28, 2012.
\newblock Submitted for publication.

\bibitem[Ramanujan(1998)]{ram:poz1998}
C.~Ramanujan, S.~and~Pozrikidis.
\newblock Deformation of liquid capsules enclosed by elastic membranes in
  simple shear flow: large deformations and the effect of fluid viscosities.
\newblock \emph{Journal of Fluid Mechanics}, 361:\penalty0 117--143, 1998.

\bibitem[Sackmann(1996)]{sac1996}
E.~Sackmann.
\newblock Supported membranes: Scientific and practical applications.
\newblock \emph{Science}, 271:\penalty0 43--48, 1996.

\bibitem[Skeel(1982)]{ske1982}
Robert~D. Skeel.
\newblock {A Theoretical Framework for Proving Accuracy Results for Deferred
  Corrections}.
\newblock \emph{SIAM Journal on Numerical Analysis}, 19\penalty0 (1):\penalty0
  171--196, 1982.

\bibitem[Sohn et~al.(2010)Sohn, Tseng, Li, Voigt, and
  Lowengrub]{soh:tse:li:voi:low2010}
J.S. Sohn, Y.-H. Tseng, S.~Li, A.~Voigt, and J.S. Lowengrub.
\newblock Dynamics of multicomponent vesicles in a viscous fluid.
\newblock \emph{Journal of Computational Physics}, 229:\penalty0 119--144,
  January 2010.

\bibitem[Speck et~al.(2014)Speck, Ruprecht, Emmett, Bolten, and
  Krause]{spe:rup:emm:bol:kra2013}
Robert Speck, Daniel Ruprecht, Matthew Emmett, Matthias Bolten, and Rolf
  Krause.
\newblock A space-time parallel solver for the three-dimensional heat equation.
\newblock In \emph{Parallel Computing: Accelerating Computational Science and
  Engineering (CSE)}, volume~25 of \emph{Advances in Parallel Computing}, pages
  263 -- 272. IOS Press, 2014.

\bibitem[Stetter(1973)]{ste1973}
Hans~J. Stetter.
\newblock \emph{Analysis of Discretization Methods for Ordinary Differential
  Equations}, volume~23.
\newblock Springer Tracts in Natural Philosophy, 1973.

\bibitem[Sukumaran and Seifert(2001)]{suk:sei2001}
S.~Sukumaran and U.~Seifert.
\newblock {Influence of shear flow on vesicles near a wall: A numerical study}.
\newblock \emph{Physical Review Letter E}, 64\penalty0 (011916), 2001.

\bibitem[Veerapaneni et~al.(2009)Veerapaneni, Gueyffier, Zorin, and
  Biros]{vee:gue:zor:bir2009}
S.~K. Veerapaneni, D.~Gueyffier, D.~Zorin, and G.~Biros.
\newblock A boundary integral method for simulating the dynamics of
  inextensible vesicles suspended in a viscous fluid in 2{D}.
\newblock \emph{Journal of Computational Physics}, 228\penalty0 (7):\penalty0
  2334--2353, 2009.

\bibitem[Veerapaneni et~al.(2011)Veerapaneni, Rahimian, Biros, and
  Zorin]{vee:rah:bir:zor2011}
Shravan~K. Veerapaneni, Abtin Rahimian, George Biros, and Denis Zorin.
\newblock A fast algorithm for simulating vesicle flows in three dimensions.
\newblock \emph{Journal or Computational Physics}, 230\penalty0 (14):\penalty0
  5610--5634, 2011.

\bibitem[Weiser(2013)]{wei2013}
Martin Weiser.
\newblock {Faster SDC convergence on non-equidistant grids by DIRK sweeps}.
\newblock Technical Report 13-30, ZIB, Takustr.7, 14195 Berlin, 2013.

\bibitem[Zhao and Shaqfeh(2011)]{zha:sha2011a}
Hong Zhao and Eric S.~G. Shaqfeh.
\newblock The dynamics of a vesicle in simple shear flow.
\newblock \emph{Journal of Fluid Mechanics}, 674:\penalty0 578--604, 2011.

\bibitem[Zhao and Shaqfeh(2013)]{zha:sha2013b}
Hong Zhao and Eric S.~G. Shaqfeh.
\newblock The dynamics of a non-dilute vesicle suspension in simple shear flow.
\newblock \emph{Journal of Fluid Mechanics}, 725:\penalty0 709--731, 2013.

\end{thebibliography}
\end{document}